\documentclass{amsart}
\usepackage{amssymb}

\usepackage{graphicx}

\usepackage[all]{xy}
 \newtheorem{theorem}{Theorem}[section]

 \newtheorem{proposition}[theorem]{Proposition}
 \newtheorem{lemma}[theorem]{Lemma}
 \newtheorem{corollary}[theorem]{Corollary}
 \newtheorem{remark}[theorem]{Remark}
 \newtheorem{definition}[theorem]{Definition}
 \newtheorem{notation}[theorem]{Notation}

 \newtheorem{claim}[theorem]{Claim}

 \numberwithin{equation}{section}
 \def\subrel#1#2{\mathrel{\mathop{#2}\limits_{#1}}}
\allowdisplaybreaks[1]

\begin{document}

\title
[Hochschild cohomology of  the quadratic monomial algebra ${\rm N}_m$
]
{
Hochschild cohomology of  the quadratic monomial algebra ${\rm N}_m$  
}
\author{Tomohiro Itagaki, Kazunori Nakamoto,  and Takeshi Torii}

\address{Department of Economics, 
Faculty of Economics, 
Takasaki City University of Economics, 
Gunma 370-0801, Japan} 
\email{titagaki@tcue.ac.jp}

\address{Center for Medical Education and Sciences,
Faculty of Medicine, 
University of Yamanashi,
Yamanashi 409--3898, Japan}
\email{nakamoto@yamanashi.ac.jp}

\address{Department of Mathematics, 
%Faculty of Science, 
Okayama University,
Okayama 700--8530, Japan}
\email{torii@math.okayama-u.ac.jp}
\thanks{The first author was partially supported by 
JSPS KAKENHI Grant Number JP17K14175. 
The second author was partially supported by 
JSPS KAKENHI Grant Number JP20K03509. 
The third author was partially supported by
JSPS KAKENHI Grant Numbers JP17K05253 and JP23K03113.}



\date{March 29, 2024\ ({\tt version~1.0.0})}

\begin{abstract}
Let ${\rm N}_m(R) = \{ (a_{ij}) \in {\rm M}_m(R) \mid 
a_{11} = a_{22} = \cdots = a_{mm} \mbox{ and } 
a_{ij} = 0  \mbox{ for any } i > j   \}$ for a commutative ring $R$.  
Then ${\rm N}_m(R)$ is a quadratic monomial algebra over $R$. We calculate %the Hochschild cohomology 
${\rm HH}^{\ast}({\rm N}_m(R), {\rm M}_m(R)/{\rm N}_m(R))$ as $R$-modules.  We also determine the $R$-algebra structure of  the Hochschild cohomology ring ${\rm HH}^{\ast}({\rm N}_m(R), {\rm N}_m(R))$.   
For $m \ge 3$, ${\rm HH}^{\ast}({\rm N}_m(R), {\rm N}_m(R))$ is an infinitely generated algebra over $R$ and has no 
Batalin-Vilkovisky algebra structure giving the Gerstenhaber bracket.  
\end{abstract}

\maketitle

\section{Introduction}
For a commutative ring $R$, set 
\[ 
{\rm N}_m(R) = \left\{
\left(
\begin{array}{ccccc}
a & \ast & \ast & \cdots & \ast \\
0 & a & \ast & \cdots & \ast \\ 
0 & 0 & a & \cdots & \ast \\
\vdots & \vdots & \ddots & \ddots & \vdots \\
0 & 0 & 0 & \cdots & a 
\end{array}
\right) \in {\rm M}_m(R) 
\right\}. 
\]
%Putting 
%\[
%x_1 = E_{12}, \ x_2 = E_{23}, \ \ldots, \ x_{m-1} = E_{m-1, m} \in A,   
%\]
%we have an isomorphism $A \cong R\langle x_1, x_2, \ldots, x_{m-1} \rangle/\langle x_i x_j \mid j \neq i+1 \rangle$ as $R$-algebras. 
In this paper, we calculate ${\rm HH}^{\ast}({\rm N}_m(R), {\rm M}_m(R)/{\rm N}_m(R))$ as $R$-modules.  
We also calculate the Hochschild cohomology ring ${\rm HH}^{\ast}({\rm N}_m(R), {\rm N}_m(R))$ as $R$-algebras. 
Moreover, we calculate the Gerstenhaber bracket on ${\rm HH}^{\ast}({\rm N}_m(R), {\rm N}_m(R))$ and show that ${\rm HH}^{\ast}({\rm N}_m(R), {\rm N}_m(R))$ has no Batalin-Vilkovisky algebra structure which gives the Gerstenhaber bracket. 
\if
We define the two-sided ideal ${\rm J}_m(R)$ of ${\rm N}_m(R)$ over $R$ by 
\[
{\rm J}_m(R) = \{ (a_{ij}) \in {\rm N}_m(R) \mid a_{11} = a_{22} = \cdots = a_{mm} = 0 \}.  
\] 
We set $A={\rm N}_m(R)$.
We denote by ${}_A{\rm BMod}_A(R)$
the category of $A$-$A$-bimodules over $R$.
Let ${\rm J}_m(R)$ be the Jacobson radical of $A$.
We define 
\[ F^rA={\rm J}_m(R)^r \]
for $r\ge 0$.
Note that $F^0A=A$ and $F^mA=0$.
We obtain a filtration 
\[ A=F^0A\supset F^1A\supset F^2A\supset\cdots
     \supset F^{m-1}A\supset F^mA=0\] 
of $A$ in ${}_A{\rm BMod}_A(R)$.
\fi

In \cite{Nakamoto-Torii:Hochschild}, we have calculated the Hochschild cohomology ${\rm HH}^{\ast}(A, {\rm M}_3(k)/A)$ for any $k$-subalgebra $A$ of ${\rm M}_3(k)$ over an algebraically closed field $k$. 
The $k$-subalgebra ${\rm N}_3(k)$ is one of the most difficult $k$-subalgebras of ${\rm M}_3(k)$ to calculate ${\rm HH}^{\ast}(A, {\rm M}_3(k)/A)$. Indeed, we could not calculate ${\rm HH}^{\ast}({\rm N}_3(k), {\rm M}_3(k)/{\rm N}_3(k))$ until we used spectral sequences.  Hence, it is a challenging task to calculate ${\rm HH}^{\ast}({\rm N}_{m}(R), {\rm M}_m(R)/{\rm N}_m(R))$ for $m \ge 4$. It is also a hard job to determine the $R$-algebra structure of  
${\rm HH}^{\ast}({\rm N}_{m}(R), {\rm N}_m(R))$ for $m \ge 3$.

Setting $x_1 =E_{1, 2}, x_2=E_{2, 3}, \ldots, x_{m-1} = E_{m-1, m} \in {\rm N}_m(R)$, 
we have an isomorphism as $R$-algebras  
\[ 
{\rm N}_m(R) \cong R\langle x_1, x_2, \ldots, x_{m-1} \rangle/\langle x_{i}x_{j} \mid j \neq i+1 \rangle, 
\]  
where $E_{i, j}$ is the $(i, j)$-th matrix unit in ${\rm M}_m(R)$. 
Note that ${\rm N}_m(R)$ is a quadratic monomial algebra over $R$ with degree $| x_i | =1$ ($1 \le i \le m-1$).   
When $R$ is a field,  
${\rm N}_m(R)$ is a Koszul algebra over $R$. 
The quadratic dual algebra ${\rm N}_m(R)^{!}$ of ${\rm N}_m(R)$ is isomorphic to $R\langle y_1, y_2, \ldots, y_{m-1} \rangle/\langle y_{i}y_{i+1} \mid 
1 \le i \le m-2 \rangle$. %Let $|y_i| = 1$ for $1 \le i \le m-1$.  
Then ${\rm N}_m(R)^{!}$ is also a graded $R$-algebra 
with degree $| y_i | =1$ ($1 \le i \le m-1$). 
Put 
\[
\varphi(d) = {\rm rank}_{R} {\rm N}_m(R)^{!}_d, 
\]
where ${\rm N}_m(R)^{!}_d$ is the homogeneous part of ${\rm N}_m(R)^{!}$ of degree $d$. 
The Poincar\'e series $\displaystyle f^{!}(t) = \sum_{d \ge 0} \varphi(d)t^d$ of ${\rm N}_m(R)^{!}$ 
can be calculated by 
$\displaystyle f^{!}(t) = 1/(1+\sum_{k=1}^{m-1} (-1)^k (m-k)t^k)$ (Proposition~\ref{prop:ff^{!}}). 
%\frac{1}{\displaystyle 1+\sum_{k=1}^{m-1} (-1)^k (m-k)t^k}$. 

The first main theorem is the following: 

\begin{theorem}[Theorem~\ref{th:HHMmNm} and Corollary~\ref{cor:HHMmNm}]
Let $m \ge 3$. The Hochschild cohomology ${\rm HH}^{n}({\rm N}_m(R), {\rm M}_m(R)/{\rm N}_m(R))$ is a free $R$-module for $n\ge 0$. The rank of ${\rm HH}^{n}({\rm N}_m(R), {\rm M}_m(R)/{\rm N}_m(R))$ for $n\ge 0$ is given by 
\[
{\rm rank}_{R} {\rm HH}^{n}({\rm N}_m(R), {\rm M}_m(R)/{\rm N}_m(R)) = \left\{ 
\begin{array}{cc}
m-1 & (n=0),  \\
(m-2)\varphi(n) & (n>0).  
\end{array}
\right. 
\]
\end{theorem}

The second main theorems are the following: 

\begin{theorem}[Theorem~\ref{th:mainthm}]
Let $m \ge 3$. The Hochschild cohomology ${\rm HH}^{n}({\rm N}_m(R), {\rm N}_m(R))$ is a free $R$-module for $n\ge 0$. The rank of ${\rm HH}^{n}({\rm N}_m(R), {\rm N}_m(R))$ is given by 
\begin{eqnarray*} 
& & {\rm rank}_{R} {\rm HH}^{n}({\rm N}_m(R), {\rm N}_m(R)) \\ 
& = & \left\{
\begin{array}{lc}
2 &  (n=0), \\
2m-4 & (n=1), \\
\displaystyle \varphi(n)+(m-4)\varphi(n-1) +(-1)^m\varphi(n-m+1) +\sum_{k=2}^{m-1} (-1)^k (k+1)\varphi(n-k) & (n\ge 2).  
\end{array}
\right.
\end{eqnarray*} 
\end{theorem}

\begin{theorem}[Theorems~\ref{th:productzerom=3}  and \ref{th:productzerom>=4} and Corollary~\ref{cor:infinitelygenerated}]
Let $m\ge3$. 
There is an augmentation map $\epsilon : {\rm HH}^{\ast}({\rm N}_m(R), {\rm N}_m(R)) \to R$ as an $R$-algebra homomorphism such that the Kernel $\overline{{\rm HH}^{\ast}}({\rm N}_m(R), {\rm N}_m(R))$ of $\epsilon$ satisfies 
\[\overline{{\rm HH}^{\ast}}({\rm N}_m(R), {\rm N}_m(R))  \cdot \overline{{\rm HH}^{\ast}}({\rm N}_m(R), {\rm N}_m(R)) =0.
\]
In particular,  ${\rm HH}^{\ast}({\rm N}_m(R), {\rm N}_m(R))$ is an infinitely generated algebra over $R$. 
\end{theorem}

\begin{theorem}[Theorem~\ref{th:noBVstructure}] 
For $m \ge 3$, ${\rm HH}^{\ast}({\rm N}_m(R), {\rm N}_m(R))$ has no Batalin-Vilkovisky algebra structure over $R$ which gives the Gerstenhaber bracket  $[\;, \;]$. 
\end{theorem}

\bigskip 

As an application of the main theorems, 
we can calculate the dimension of the tangent space of the moduli of subalgebras of 
${\rm M}_m$ over ${\Bbb Z}$ at ${\rm N}_m$. Set $\displaystyle d = {\rm rank}_{R} {\rm N}_m(R) = \frac{m^2-m+2}{2}$. 

\begin{theorem}[Theorem~\ref{th:dimtangentspace}] 
The dimension of the Zariski tangent space of the moduli of rank $d$ subalgebras of ${\rm M}_m$ over ${\Bbb Z}$ at ${\rm N}_m$ is 
\[
\dim T_{{\rm Mold}_{m, d}/{\Bbb Z}, {\rm N}_m} = \frac{3m^2 - 7m+4}{2} 
\]
for $m \ge 3$. 
\end{theorem}

\bigskip 

The organization of this paper is as follows: in Section 2, we review Hochschild cohomology. 
In Section 3, we introduce several results on spectral sequences. 
In Section 4, we show that ${\rm HH}^{\ast}({\rm N}_m(R), R) \cong {\rm N}_m(R)^{!}$ as $R$-algebras. 
We also describe the Poincar\'e series $f^{!}(t)$ of ${\rm N}_m(R)^{!}$ explicitly.  
In Section 5, we calculate ${\rm HH}^{\ast}({\rm N}_m(R), {\rm M}_m(R)/{\rm N}_m(R))$ as $R$-modules. 
We also calculate the dimension of the tangent space of the moduli of subalgebras of 
${\rm M}_m$ over ${\Bbb Z}$ at ${\rm N}_m$. 
In Section 6, we determine the $R$-module structure of ${\rm HH}^{\ast}({\rm N}_m(R), {\rm N}_m(R))$.  
In Section 7, we determine the product structure of ${\rm HH}^{\ast}({\rm N}_m(R), {\rm N}_m(R))$. 
In Section 8, we describe the Gerstenhaber bracket on ${\rm HH}^{\ast}({\rm N}_m(R), {\rm N}_m(R))$. 
We also show that ${\rm HH}^{\ast}({\rm N}_m(R), {\rm N}_m(R))$ has no Batalin-Vilkovisky algebra structure giving 
the Gerstenhaber bracket $[\;, \;]$ for $m\ge 3$.  
In Section 9, we deal with ${\rm HH}^{\ast}({\rm N}_2(R), {\rm M}_2(R)/{\rm N}_2(R))$ and 
${\rm HH}^{\ast}({\rm N}_2(R), {\rm N}_2(R))$ as an appendix. 

\bigskip 

Throughout this paper, $R$ denotes a commutative ring and $A$ an associative algebra over $R$. 
We denote by $E_{i,j} \in {\rm M}_m(R)$ the matrix with entry $1$ in the $(i, j)$-component and $0$ the other components. We also denote by $I_m$ the identity matrix in ${\rm M}_m(R)$. 
By a {\it module} $M$ over $A$, we mean a left module $M$ over $A$, unless stated otherwise. 
We set ${\rm B}_m(R) = \{ (a_{ij}) \in {\rm M}_m(R) \mid a_{ij} = 0 \mbox{ for } i>j \}$ and 
$J({\rm N}_m(R)) = \{ (a_{ij}) \in {\rm M}_m(R) \mid a_{ij} = 0 \mbox{ for } i\ge j \}$.   
For a subset $S$ of an $R$-module $M$, we denote by $R\{ S \}$ the $R$-submodule of $M$ generated by 
$S$. For a homogeneous element $x$ of a graded $R$-module $M=\oplus_{i \in {\Bbb Z}} M_i$ (or a graded 
$R$-algebra $A=\oplus_{i \in {\Bbb Z}} A_i$),  we denote by $|x|$ the degree of $x$.

\section{Preliminaries on Hochschild cohomology}
In this section, we make a brief survey of Hochschild cohomology  
%
% some preparations 
%for calculating % the Hochschild cohomology ring 
%${\rm HH}^{\ast}({\rm N}_m(R), {\rm N}_m(R))$ and 
%${\rm HH}^{\ast}({\rm N}_m(R), {\rm M}_m(R)/{\rm N}_m(R))$. 
%
%\subsection{Definition of Hochschild cohomology}
%
%In this subsection, 
%we give a review of Hochschild cohomology groups
(cf.~\cite{Gerstenhaber} and \cite{Witherspoon}). 
Throughout this section, 
$M$ denotes an $A$-bimodule over $R$. 
%$R$ denotes a commutative ring, $A$ an associative algebra over $R$, and $M$ an $A$-bimodule over $R$. 

\begin{definition}\label{def:hochschildcohomology}\rm 
Assume that $A$ is a projective module over $R$. 
Let $A^e = A\otimes_{R} A^{op}$ be the enveloping algebra of $A$. 
For $A$-bimodules $A$ and $M$ over $R$, we can regard them as $A^e$-modules.  
We define the $i$-th Hochschild cohomology group ${\rm HH}^i (A, M)$ as 
${\rm Ext}^{i}_{A^{e}} (A, M)$.  
\end{definition}

We denote by
$B_{\ast}(A,A,A)$
the bar resolution of $A$ as $A$-bimodules over $R$.
For $p\ge 0$,
we have
\[ B_p(A,A,A)=A\otimes_R \overbrace{A\otimes_R\cdots\otimes_RA}^p
              \otimes_R A.\]
For an $A$-bimodule $M$ over $R$,
we define a cochain complex $C^{\ast}(A, M)$ to be 
\[ {\rm Hom}_{A^e}(B_{\ast}(A,A,A),M). \]
We can identify
$C^p(A, M)$ with an $R$-module
\[ {\rm Hom}_R(\overbrace{A\otimes_R\cdots\otimes_RA}^{p},M). \]
Under this identification,
the coboundary map $d^p:C^p(A, M)\to C^{p+1}(A, M)$ is given by
\[ \begin{array}{rcl}
   d^p(f)(a_1\otimes\cdots\otimes a_{p+1})&=&
   a_1\cdot f(a_2\otimes \cdots\otimes a_{p+1})\\
  && \displaystyle +\sum_{i=1}^p(-1)^i f(a_1\otimes\cdots\otimes a_ia_{i+1}
                     \otimes\cdots\otimes a_{p+1})\\
  &&+(-1)^{p+1}f(a_1\otimes\cdots\otimes a_p)\cdot a_{p+1}\\
  \end{array}\]  
for $f\in C^p(A, M)$ $(p \ge 1)$ and 
\[
d^{0}(m)(a) = am-ma
\]
for $m \in C^{0}(A, M)=M$. 
The Hochschild cohomology group ${\rm HH}^{\ast}(A, M)$ of $A$ with coefficients 
in $M$ can be calculated by taking the cohomology of the cochain 
complex $C^{\ast}(A, M)$:
\[ {\rm HH}^{\ast}(A, M)=H^{\ast}(C^{\ast}(A, M)).\]     

\begin{remark}\rm
In Definition \ref{def:hochschildcohomology}, the assumption that $A$ is a projective module over $R$ is needed for 
${\rm Ext}^{i}_{A^{e}} (A, M) \cong H^{i}(C^{\ast}(A, M))$ for $i \ge 0$. 
\end{remark} 

Let $N$ be another $A$-bimodule over $R$.
We define a map 
\[ \cup: C^{\ast}(A, M)\times C^{\ast}(A, N)\longrightarrow
   C^{\ast}(A, M\otimes_A N)\]
by
\[ (f\cup g)(a_1\otimes\cdots\otimes a_p\otimes 
    b_1\otimes\cdots\otimes b_q)=
   f(a_1\otimes\cdots \otimes a_p)\otimes
   g(b_1\otimes \cdots\otimes b_q)\]
for $f\in C^p(A, M)$ and $g\in C^q(A, N)$.
The map $\cup$ is $R$-bilinear and satisfies
\[ d^{p+q}(f\cup g)=
   d^p(f)\cup g+ (-1)^pf\cup d^q(g).\]
Hence the map $\cup$ induces a map
\[ {\rm HH}^p(A, M)\otimes_R{\rm HH}^q(A, N)\longrightarrow
   {\rm HH}^{p+q}(A, M\otimes_A N)\]
of $R$-modules.

By the above construction,
we see that the Hochschild cohomology
${\rm HH}^{\ast}(A, -)$
defines a lax monoidal functor
from the monoidal category of $A$-bimodules over $R$
to the monoidal category of graded $R$-modules.
Hence, 
${\rm HH}^{\ast}(A, M)$ is a 
graded associative algebra over $R$ 
if $M$ is a monoid object in
the category of $A$-bimodules over $R$.

Suppose that the unit map 
$R\to A$ is a split monomorphism.
We set $\overline{A}=A/RI$,
where $I\in A$ is the image of $1\in R$ under the unit map.
Let $\overline{B}_{\ast}(A,A,A)$ be 
the reduced bar resolution of $A$ as $A$-bimodules over $R$.
We have
\[ \overline{B}_p(A,A,A)\cong
   A\otimes_R \overbrace{\overline{A}\otimes_R\cdots\otimes_R
   \overline{A}}^p\otimes_R A\]
for $p\ge 0$.
For an $A$-bimodule $M$ over $R$, 
we denote  the cochain complex
${\rm Hom}_{A^e}(\overline{B}_{\ast}(A,A,A),M)$
by $\overline{C}^{\ast}(A, M)$.
The cochain complex
$\overline{C}^{\ast}(A,M)$ is a subcomplex of $C^{\ast}(A, M)$. 
Recall that the reduced bar resolution $\overline{B}_{\ast}(A,A,A)$
is chain homotopy equivalent to
the bar resolution $B_{\ast}(A,A,A)$. 
Hence, the inclusion  
$\overline{C}^{\ast}(A, M)\to C^{\ast}(A, M)$ induces
an isomorphism 
\[ H^{\ast}(\overline{C}^{\ast}(A, M))\cong {\rm HH}^{\ast}(A, M). \] 
We observe that the map
$\cup: C^{\ast}(A, M)\times C^{\ast}(A, N)\to
C^{\ast}(A, M\otimes_A N)$ induces an $R$-bilinear map
\[ \cup: \overline{C}^{\ast}(A, M)\times \overline{C}^{\ast}(A, N)\longrightarrow
   \overline{C}^{\ast}(A, M\otimes_A N),\]
where $N$ is another $A$-bimodule over $R$.
Hence the map $\cup: \overline{C}^{\ast}(A, M)\times \overline{C}^{\ast}(A, N) \to 
\overline{C}^{\ast}(A, M\otimes_A N)$
induces the same map
${\rm HH}^p(A, M)\otimes_R{\rm HH}^q(A, N) \to 
   {\rm HH}^{p+q}(A, M\otimes_A N)$
of $R$-modules as before.

\section{Spectral sequences}
In this section
we recall the construction
of spectral sequences associated
to filtered cochain complexes.
In particular,
we construct a spectral sequence
by introducing a filtration on
the Hochschild cochain complex $C^*(A,M)$
by powers of a two-sided ideal (e.g. the Jacobson radical $J$) of $A$. 

%\subsection{Existence of spectral sequences} 
\subsection{Review on the construction 
of spectral sequences associated to
filtered cochain complexes}\label{subsection:constructionss}\label{subsection:reviewss} 

We can consider spectral sequences in 
an abelian category $\mathcal{A}$. 
In this subsection
we recall the construction of spectral sequences
associated to filtered cochain complexes in $\mathcal{A}$.

Let $(C^*,d)$ be a cochain complex  in $\mathcal{A}$
equipped with a filtration
\[ C^*=F^0C^*\supset F^1C^*\supset \cdots 
   \supset F^pC^*\supset\cdots \]
by subcomplexes.
Throughout this paper we assume that
there exists $t\in\mathbb{Z}_{>0}$
such that $F^tC^*=0$.
We say that $(C^*,d,\{F^pC^*\}_{p\ge 0})$ is 
a filtered differential graded module in $\mathcal{A}$.

For a filtered differential graded module
$(C^*,d,\{F^pC^*\}_{p\ge 0})$,
we can construct an associated spectral sequence 
\[ E_1^{p,q}(C^*)\Longrightarrow H^{p+q}(C^*) \]
with
\[ d_r^{p,q}: E_r^{p,q}(C^*)\longrightarrow E_r^{p+r,q-r+1}(C^*) \]
(see, for example, \cite[Theorem~2.6]{McCleary}).
Note that we have an isomorphism
\[ E_1^{p,q}(C^*)\cong H^{p+q}(F^pC^*/F^{p+1}C^*).\] 
The differential
$d_1: E_1^{p,q}(C^*)\to E_1^{p+1,q}(C^*)$
is identified with
the connecting homomorphism
\[ H^{p+q}(F^pC^*/F^{p+1}C^*)\longrightarrow
   H^{p+q+1}(F^{p+1}C^*/F^{p+2}C^*) \]
associated to the short exact sequence
\[ 0 \longrightarrow F^{p+1}C^*/F^{p+2}C^*\longrightarrow
        F^pC^*/F^{p+2}C^*\longrightarrow
        F^pC^*/F^{p+1}C^*\longrightarrow 0\]
of cochain complexes.
If $F^tC^*=0$,
then the spectral sequence collapses
from the $E_t$-page.

%Now, we suppose that $\mathcal{A}$
%is an abelian monoidal category, in which
%the tensor product $\otimes: \mathcal{A}\times\mathcal{A}
%\to\mathcal{A}$ is right exact separately
%in each variable.
%Let $(C^*,d)$ be 
%a differential graded algebra in $\mathcal{A}$
%equipped with a filtration
%$C^*=F^0C^*\supset F^1C^*\supset \cdots 
%   \supset F^pC^*\supset\cdots$
%by subcomplexes.
%We say that $(C^*,d,\{F^pC^*\}_{p\ge 0})$ is 
%a filtered differential graded algebra in $\mathcal{A}$
%if 
%\[ F^rC^*\cdot F^sC^*\subset F^{r+s}C^* \]
%for any $r,s\ge 0$. 
%To a filtered differential graded algebra 
%$(C^*,d,\{F^pC^*\}_{p\ge 0})$,
%we obtain 
%an associated  spectral sequence 
%$E_1^{p,q}(C^*)\Longrightarrow H^{p+q}(C^*)$ 
%of algebras
%(see, for example, \cite[Theorem~2.14]{McCleary}). 

%\bigskip 

Now, we suppose that $\mathcal{A}$
is an abelian monoidal category, 
in which the tensor product
$\otimes: \mathcal{A}\times\mathcal{A}\to\mathcal{A}$
is right exact separately in each variable.

\begin{definition}\rm
Let $(A^*,d)$ be a differential graded algebra in $\mathcal{A}$.
Suppose that we have a filtration
\[ A^*=F^0A^*\supset F^1A^*\supset\cdots\supset
       F^nA^*\supset\cdots \supset F^t A^* = 0.\]
A triple $(A^*,d,\{F^pA^*\}_{p\ge 0})$ is said to be
a filtered differential graded algebra
if it satisfies the following two conditions:
\begin{enumerate}
\item[(1)]
For any $p\ge 0$, $d(F^pA^*)\subset F^pA^*$.
\item[(2)]
For any $r,s\ge 0$, $F^rA^*\cdot F^sA^*\subset F^{r+s}A^*$.
\end{enumerate}
\end{definition}

\bigskip 

To a filtered differential graded algebra 
$(A^*,d,\{F^pA^*\}_{p\ge 0})$, 
there is a spectral sequence 
\[ E_1^{p,q}=H^{p+q}(F^pA/F^{p+1}A)
   \Longrightarrow H^{p+q}(A) \]
of algebras in $\mathcal{A}$,
which converges to $H^{p+q}(A)$
as an algebra (see, for example, \cite[Theorem~2.14]{McCleary}).

\subsection{Filtrations and spectral sequences
on Hochschild cochain complexes}

In this subsection 
we consider filtrations on Hochschild complexes 
and associated spectral sequences.

Let $R$ be a commutative ring and
let $A$ be an associative algebra over $R$.
We assume that $A$ is a projective module over $R$. 
For an $A$-bimodule $M$ over $R$, 
we denote by $C^*(A,M)$ the Hochschild cochain complex.

First, we suppose that there exists a filtration of 
$A$-bimodules over $R$: 
\[  M=F^{0}M \supset F^{1}M \supset
    \cdots \supset F^pM \supset \cdots 
     \supset F^tM=0. \]
We denote by ${\rm Gr}^{p}(M)$ 
the $p$-th associated graded module $F^{p}M/F^{p+1}M$. 
Using the filtration
$\{F^pM\}_{p\ge 0}$ on $M$,
we can introduce a filtration
$\{F^pC^*(A,M)\}_{p\ge 0}$ on $C^*(A,M)$ by
\[ F^pC^*(A,M)= C^*(A,F^pM).\]
Hence we obtain the following proposition.

\if0
Hence we obtain a spectral sequence  
\[ E_1^{p,q}={\rm HH}^{p+q}(A, {\rm Gr}^p(M))\Longrightarrow
   {\rm HH}^{p+q}(A, M)\]
of $R$-modules with
\[ d_r: E_r^{p,q}\longrightarrow E_r^{p+r,q-r+1} \]
for $r\ge 1$, where ${\rm Gr}^p(M) = F^{p}M/F^{p+1}M$. 
\fi

\begin{proposition}\label{prop:spectralsequence} 
For an $A$-bimodule $M$ over $R$
equipped with a filtration   
$M=F^{0}M \supset F^{1}M \supset
    \cdots \supset F^pM \supset \cdots
    \supset F^tM=0$,
there exists a spectral sequence  
\[ E_1^{p,q}(A,M)\cong{\rm HH}^{p+q}(A, {\rm Gr}^p(M))\Longrightarrow
   {\rm HH}^{p+q}(A, M)\]
of $R$-modules with
\[ d_r: E_r^{p,q}(A,M)\longrightarrow E_r^{p+r,q-r+1}(A,M) \]
for $r\ge 1$, where ${\rm Gr}^p(M) = F^{p}M/F^{p+1}M$. 
\end{proposition} 

\if0
Note that the differential 
$d_{1} : E_1^{p,q}\longrightarrow E_1^{p+1,q}$ is 
identified with the connecting homomorphism 
\[ {\rm HH}^{p+q}(A, {\rm Gr}^{p}(M)) \to 
   {\rm HH}^{p+q+1}(A, {\rm Gr}^{p+1}(M)) \] 
of the long exact sequence 
%\[ \cdots \to {\rm HH}^{\ast}(A, {\rm Gr}^{p+1}(M)) 
%   \to {\rm HH}^{\ast}(A, F^{p}/F^{p+2}) \to 
%   {\rm HH}^{\ast}(A, {\rm Gr}^{p}(M)) 
%   \to {\rm HH}^{\ast+1}(A, {\rm Gr}^{p+1}(M)) \to \cdots \]  
induced by the short exact sequence 
\[ 0 \to {\rm Gr}^{p+1}(M) \to F^{p}/F^{p+2} \to {\rm Gr}^{p}(M) \to 0. \] 
%Moreover, the spectral sequence collapses at the $E_m$-page. 
\fi

%\proof 
%See, for example, \cite[\S2.2]{McCleary} for construction
%of spectral sequences.
%\qed

%\subsection{Spectral sequences of algebras}

In particular,
we consider a filtration on the $A$-bimodule $A$
by powers of a two-sided ideal.
%We set $A={\rm N}_m(R)$.
%We denote by ${}_A{\rm BMod}_A(R)$
%the category of $A$-$A$-bimodules over $R$.
%Let ${\rm J}_m(R)$ be the Jacobson radical of $A$.
Let $J$ be a two-sided ideal of $A$. 
We assume that $J^t=0$ for some $t>0$.
%We define a filtration $\{F^pA\}_{p\ge 0}$ on $A$ by
%\[ F^pA=J^p \]
%for $p\ge 0$.
%Note that $F^0A=A$.
%We obtain a filtration 
%\[ A=F^0A\supset F^1A\supset F^2A\supset\cdots
%     \supset F^{m-1}A\supset F^mA=0\] 
%of $A$ in ${}_A{\rm BMod}_A(R)$.
\if0
Recall that we have a Hochschild cochain complex
$C^*(A,M)$ for $M\in{}_A{\rm BMod}_A(R)$.
We consider a functor
\[ C^*(A,-): {}_A{\rm BMod}_A(R)\longrightarrow {\rm CoCh}(R),\]
where ${\rm CoCh}(R)$ is the category of 
cochain complexes over $R$.
We set
\[ F^rC^*(A,A)=C^*(A,F^rA) \]
for $r\ge 0$.
Since $A$ is a finitely generated free module
over $R$,
we see that $F^rC^*(A,A)$
is a subcomplex of $C^*(A,A)$.
\fi
By setting  
\[ F^pA=J^p\]
for $p\ge 0$,
there is a filtration
\[ A=F^0A\supset F^1A\supset\cdots \supset F^pA\supset\cdots
     \supset F^tA=0 \]
of $A$ by $A$-bimodules over $R$.
From this filtration,
we obtain a filtration
$\{C^*(A,J^p)\}_{p\ge 0}$ 
on $C^*(A,A)$. 
Recall that
$C^*(A,A)$ is a differential graded algebra over $R$.
We can easily verify that
\[ C^*(A,J^r)\cdot C^*(A,J^s)\subset C^*(A,J^{r+s}). \]
Thus, the triple $(C^*(A,A),d,\{C^*(A,J^r)\}_{r\ge 0})$
is a filtered differential graded algebra
over $R$,
and we obtain the following proposition.

\begin{proposition}\label{prop:spectral-sequence-A-A-version}
There is a spectral sequence of $R$-algebras
\[ {}^JE_1^{p,q}(A,A)
   \Longrightarrow {\rm HH}^{p+q}(A,A), \]
where
\[ {}^JE_1^{p,q}(A,A)\cong {\rm HH}^{p+q}(A,J^p/J^{p+1}).\]
\end{proposition}

%\proof
%The corollary follows from \cite[Theorem~2.14]{McCleary}.
%\qed

Now, we consider the induced filtration
on $A$-bimodules over $R$.  
Let $M$ be an $A$-bimodule over $R$.
By setting
\begin{eqnarray} 
\overline{J}^pM=\sum_{a+b=p}J^aMJ^b \label{eq:JpM} 
\end{eqnarray}
for $p\ge 0$,
we obtain a filtration
\[ M=\overline{J}^0M\supset 
     \overline{J}^1M\supset\cdots \supset 
     \overline{J}^pM\supset\cdots
     \supset \overline{J}^{2t-1}M=0 \]
of $M$ by $A$-bimodules over $R$.
From this filtration,
we obtain a filtration
$\{C^*(A,\overline{J}^pM)\}_{p\ge 0}$ 
on $C^*(A,M)$
and the following proposition.

\begin{proposition}\label{prop:spectralseqJ}
There is a spectral sequence of $R$-modules
\[ {}^JE_1^{p,q}(A,M) \Longrightarrow {\rm HH}^{p+q}(A,M). \]
We have an isomorphism
\[ {}^JE_1^{p,q}(A,M)\cong
   {\rm HH}^{p+q}(A,{\rm Gr}^p_J(M)),\]
where 
\[  {\rm Gr}^p_J(M)=\overline{J}^pM/\overline{J}^{p+1}M. \]

\end{proposition}

\begin{remark}\rm
Since 
\[ C^*(A,J^a)\cdot C^*(A,\overline{J}^bM)\cdot C^*(A,J^c)
   \subset C^*(A,\overline{J}^{a+b+c}M), \]
the triple $(C^*(A,M),d,\{C^*(A,\overline{J}^pM)\}_{p\ge 0})$ 
is a differential graded bimodule 
over the differential graded algebra
$(C^*(A,A),d,\{C^*(A,J^p)\}_{p\ge 0})$.
Thus,
the spectral sequence
$\{{}^JE_r^{*,*}(A,M),d_r\}_{r\ge 1}$
is a bimodule over 
the spectral sequence
$\{{}^JE_r^{*,*}(A,A),d_r\}_{r\ge 1}$.
\end{remark}
%}

\if0
\subsection{Spectral sequence with degree}

Let $C^*$ be a cochain complex.
Suppose that we have a filtration
\[ C^*=F^0C^*\supset F^1C^*\supset \cdots 
   \supset F^pC^*\supset\cdots\]
of $C^*$ by subcomplexes.
We denote by
\[ E_1^{p,q}\Longrightarrow H^{p+q}(C^*) \]
ths associated spectral sequence with
\[ d_r^{p,q}: E_r^{p,q}\longrightarrow E_r^{p+r,q-r+1}. \]
Note that we have
\[ E_1^{p,q}\cong H^{p+q}(F^pC^*/F^{p+1}C^*).\] 

Now, suppose that $C^*$ has a decomposition
\[ C^*=\bigoplus_{s\in\mathbb{Z}}C^{*,s} \]
which is compatible with the filtration
\[ F^pC^*=\bigoplus_{s\in\mathbb{Z}}(F^pC^*\cap C^{*,s}).\]
Set $F^pC^{*,s} = F^pC^{*}\cap C^{*,s}$. 
In this case we have a filtration
\[ C^{*,s}=F^0C^{*,s}\supset F^1C^{*,s}\supset\cdots
   \supset F^pC^{*,s}\supset\cdots\]
which induces a spectral sequence
\[ E_1^{p,q,s}\Longrightarrow H^{p+q}(C^{*,s}).\]  
We have a decomposition
\[ E_r^{p,q}=\bigoplus_{s\in\mathbb{Z}}E_r^{p,q,s} \]
and 
\[ d_r^{p,q}=\bigoplus_{s\in\mathbb{Z}} d_r^{p,q,s}, \]
where 
\[ d_r^{p,q,s}: E_r^{p,q,s}\longrightarrow E_r^{p+r,q-r+1,s}. \]
We also have a decomposition
\[ F^pH^{p+q}(C^*)=\bigoplus_{s\in\mathbb{Z}}F^pH^{p+q}(C^{*,s})\]
and an isomorphism
\[ E_{\infty}^{p,q,s}\cong
   F^pH^{p+q}(C^{*,s})/F^{p+1}H^{p+q}(C^{*,s}).\]                   
\fi

\subsection{Gradings on spectral sequences}\label{subsection:gradingonss} 

Let $\mathcal{A}$ be an abelian category with countable coproducts.  
We let $\mathcal{A}^{\mathbb{Z}}$ be
the abelian category 
of $\mathbb{Z}$-graded objects of $\mathcal{A}$
and grading-preserving morphisms.
In this subsection
we consider spectral sequences in
$\mathcal{A}^{\mathbb{Z}}$.
In this case we have trigradings on
spectral sequences.

Let 
\[ C^*=\bigoplus_{s\in\mathbb{Z}}C^{*,s} \]
be a cochain complex of $\mathcal{A}^{\mathbb{Z}}$
equipped with a filtration
\[ C^*=F^0C^*\supset F^1C^*\supset \cdots 
   \supset F^pC^*\supset\cdots
   \supset F^tC^*=0 \]
by subcomplexes,
where $C^{*,s}$ is
the component of 
(cohomological) degree $s\in\mathbb{Z}$. 
\if0
Then $C^{*,s}$ is a cochain complex of $\mathcal{A}$
equipped with a filtration
\[ C^{*,s}=F^0C^{*,s}\supset F^1C^{*,s}\supset \cdots 
   \supset F^pC^{*,s}\supset\cdots, \]
where $F^pC^{*,s}=F^pC^*\cap C^{*,s}$, 
\fi
Set $F^{p}C^{n, s} = F^{p}C^{n} \cap C^{n, s}$.  
Note that any filtration by subcomplexes in $\mathcal{A}^{\mathbb{Z}}$ 
is assumed to satisfy $F^pC^*=\bigoplus_{s\in\mathbb{Z}}F^{p}C^{*,s}$ (in other words, the filtration is compatible with the grading).  
Then we have the following associated spectral sequence as in \S\ref{subsection:constructionss}. 
%Proposition~\ref{prop:spectralsequence}.  

\begin{proposition}\label{prop:sstrigrading} 
There is a spectral sequence 
\[ E_1^{p,q}(C^*)\Longrightarrow H^{p+q}(C^*) \]
in $\mathcal{A}^{\mathbb{Z}}$,
with
\[ E_1^{p,q}(C^*)\cong H^{p+q}(F^pC^*/F^{p+1}C^*).\]
More precisely, 
let  $E_r^{p,q,s}(C^*)$ be 
the degree $s$ component of 
the spectral sequence $\{E_r^{p,q}(C^*),d_r^{p, q}\}_{r\ge 1}$.
Then $\{E_r^{p,q,s}(C^*),d_r^{p, q, s}\}_{r\ge 1}$
is a spectral sequence of $\mathcal{A}$
with
\[ d_r^{p, q, s} : E_r^{p,q,s}(C^*)\longrightarrow
        E_r^{p+r,q-r+1,s}(C^*),\]
where 
\[
E_{r}^{p, q}(C^*) = \bigoplus_{s \in {\Bbb Z}} E_{r}^{p, q, s}(C^{\ast}),  \; 
d_{r}^{p, q} = \bigoplus_{s \in {\Bbb Z}} d_{r}^{p,q, s},  
\]       
and 
\[
E_{1}^{p, q, s}(C^{\ast}) \cong H^{p+q}(F^{p}C^{\ast, s}/F^{p+1}C^{\ast, s}) \Longrightarrow H^{p+q}(C^{*, s}). 
\]
\end{proposition} 

\begin{notation}\label{notation:Hn,s}\rm
%Proposition~\ref{prop:sstrigrading}
Under the situation above, we set 
\[
H^{n, s}(C^{\ast}) = H^{n}(C^{\ast, s}). 
\]
Note that 
\[
H^{n}(C^{\ast}) = \bigoplus_{s\in {\Bbb Z}} H^{n, s}(C^{\ast}). 
\]
\end{notation}

\bigskip 

%%%%%%%% deleted part 
\if0
Then we have an associated spectral sequence
\[ E_1^{p,q}(C^*)\Longrightarrow H^{p+q}(C^*) \]
in $\mathcal{A}^{\mathbb{Z}}$,
with
\[ E_1^{p,q}(C^*)\cong H^{p+q}(F^pC^*/F^{p+1}C^*).\]

We denote by
\[ E_r^{p,q,s}(C^*) \]
the degree $s$ component of 
the spectral sequence $\{E_r^{p,q}(C^*),d_r\}_{r\ge 1}$.
Then $\{E_r^{p,q,s}(C^*),d_r\}_{r\ge 1}$
is a spectral sequence of $\mathcal{A}$
with
\[ d_r: E_r^{p,q,s}(C^*)\longrightarrow
        E_r^{p+r,q-r+1,s}(C^*).\]
\fi
%%%%%%%%

We need the following lemma
in \S\ref{section:HH-N-M/N} below.

\begin{lemma}\label{lemma:ss-collapse-internal-degree}
If there exists an integer $i$ such that
\[ H^{p+q, s}(F^pC^{*}/F^{p+1}C^{*}) = H^{p+q}(F^pC^{*,s}/F^{p+1}C^{*,s})=0 \]
for $s\neq q+i$,
then
the spectral sequence
\[ E_1^{p,q}(C^*)\Longrightarrow H^{p+q}(C^*) \]
collapses from the $E_2$-page.
\end{lemma}

\proof
By the assumption,
$E_1^{p,q,s}(C^*)=0$ for $s\neq q+i$.
Since the differential of the spectral
sequence $\{E_r^{p,q,s}(C^*),d_r\}_{r\ge 1}$
has the form
$d_r^{p,q,s}: E_r^{p,q,s}(C^*)\to
               E_r^{p+r,q-r+1,s}(C^*)$,
we see that $d_r$ is trivial unless $r=1$.
Hence the spectral sequence collapses from the $E_2$-page.
%We have a filtration on ${\rm HH}^{n,s}(A,A)$
%and an isomorphism
%\[ \mathrm{Gr}^p({\rm HH}^{n,s}(A,A))\cong E_{\infty}^{p,n-p,s}.\]
%By Lemma~\ref{lemma:graded_ss_vanishing},
%we obtain that
%\[ {\rm HH}^{n,s}(A,A)\cong E_{\infty}^{n-s,s,s}.\]
%Since we have a decomposition
%\[ {\rm HH}^n(A,A)=
%   \bigoplus_{s\in\mathbb{Z}} {\rm HH}^{n,s}(A,A) ,\]
%there is no extension problem.
%
\if0
Recall that the filtration 
$\{F^p{\rm HH}^n(A;A)\}_p$ on ${\rm HH}^n(A;A)$
is compatible with the decomposition
\[ {\rm HH}^n(A;A)=
   \bigoplus_{s\in\mathbb{Z}} {\rm HH}^{n,s}(A;A), \]
that is, we have
\[ F^p{\rm HH}^n(A;A)=\bigoplus_{s\in\mathbb{Z}}
   F^{p}{\rm HH}^{n,s}(A;A) \]
for each $p$.
This implies that 
the exact sequence
\[ 0\to F^{p+1}{\rm HH}^n(A;A)\longrightarrow
        F^p{\rm HH}^n(A;A)\longrightarrow
        E_{\infty}^{p,n-p}\to 0\]
is canonically split
by Lemma~\ref{lemma:graded_ss_vanishing}.
%\[ \mathrm{Gr}^p({\rm HH}^{n,s}(A;A))\cong E_{\infty}^{p,n-p,s}.\]
Hence we obtain a canonical isomorphism
%\[ F^p{\rm HH}^n(A;A)\cong
%   F^{p+1}{\rm HH}^n(A;A)\oplus E_{\infty}^{p,n-p}\]
%for each $p$.
%In particular, 
%we have a canonical isomorphism 
\[ {\rm HH}^n(A;A)\cong\bigoplus_p E_{\infty}^{p,n-p} \]
and
there is no extension problem.
\fi
%}
\qed

%\section{Projective resolutions}
\section{${\rm HH}^{\ast}({\rm N}_m(R), R)$}\label{section:projective} 
In this section, we show that ${\rm HH}^{\ast}({\rm N}_m(R), R) \cong {\rm N}_m(R)^{!}$ as graded $R$-algebras. 
We also obtain several results on  $\varphi(n) = {\rm rank}_{R} {\rm N}_m(R)^{!}_{n}$.  
By using these results, it is possible to calculate  
${\rm HH}^{\ast}({\rm N}_m(R), {\rm M}_m(R)/{\rm N}_{m}(R))$ and 
 ${\rm HH}^{\ast}({\rm N}_m(R), {\rm N}_{m}(R))$ as $R$-modules. Indeed,  
${\rm M}_m(R)/{\rm N}_{m}(R)$ and ${\rm N}_{m}(R)$ have filtrations of ${\rm N}_m(R)$-bimodules over $R$ whose associated graded modules are isomorphic to direct sums of copies of ${\rm N}_m(R)$-bimodules $R$ over $R$. 
By calculating spectral sequences, we will determine the $R$-module structure of ${\rm HH}^{\ast}({\rm N}_m(R), {\rm M}_m(R)/{\rm N}_{m}(R))$ in \S\ref{section:HH-N-M/N} and that of  
${\rm HH}^{\ast}({\rm N}_m(R), {\rm N}_{m}(R))$ in \S\ref{section:HH-N-N}, respectively. 

In \S\ref{subsection:Quadraticmonomialalgebra}, we deal with quadratic monomial algebras $A$ 
over a commutative ring $R$. We show that ${\rm HH}^{\ast}(A, R) \cong A^{!}$ as graded $R$-algebras. 
In \S\ref{subsection:HHNmR}, we apply the results in \S\ref{subsection:Quadraticmonomialalgebra} to the case 
$A = {\rm N}_m(R)$ and determine the $R$-algebra structure of ${\rm HH}^{\ast}({\rm N}_m(R), R)$ for $m \ge 2$. In \S\ref{subsection:resultsonvarphi}, we obtain several results on  $\varphi(n) = {\rm rank}_{R} {\rm N}_m(R)^{!}_{n}$.  

\subsection{Quadratic monomial algebras}\label{subsection:Quadraticmonomialalgebra} 
In this subsection, we deal with quadratic monomial algebras over a commutative ring $R$ (cf. \cite[Chapter~1 \S2]{Quadratic}).  
Let $\{ e_1, \ldots, e_n \}$ be an $R$-basis of a free $R$-module $V$of rank $n$. 
Let $\displaystyle T(V) = \bigoplus_{i=0}^{\infty} V^{\otimes i}$ be the tensor algebra of $V$ over $R$.  
%We say that 
%For a subset ${\mathcal I} \subseteq [1, \ldots, n]\times [1, \ldots, n]$, set 
%\[
%S=S_{{\mathcal I}} = \{ e_i\otimes e_j \in V\otimes_{R} V \mid (i, j) \in {\mathcal I} \} 
%\mbox{ and } and 
%I _{S} = R\{ S \} \subseteq V\otimes_{R} V. 
%\] and 
For a subset $S$ of $\{ e_i\otimes e_j \in V\otimes_{R} V \mid 1 \le i, j \le n \}$, set 
$I _{S} = R\{ S \} \subseteq V\otimes_{R} V$. Then we say that $A_S = T(V)/\langle I_S \rangle$ is a {\it quadratic monomial algebra} over $R$, where $\langle I_S \rangle$ is the two-sided ideal of $T(V)$ generated by $I_S$.  
We also write $A_{S} = \{ V, I_S \}$ according to \cite[Chapter~1 \S2]{Quadratic}. 
Note that $A_{S}$ is a graded $R$-algebra with $| e_i | = 1$ ($1\le i \le n$). 
Denote by $A_{S, i}$ the homogenous part of $A_{S}$ of degree $i$.  
Then $A_{S} =R\oplus A_{S, +}$ is an augmented algebra over $R$ with augmentation map $\epsilon : A_{S} \to R$, 
where $A_{S, +} = \oplus_{i>0} A_{S, i}$ and $\epsilon(A_{S, +})=0$.     

Let $V^{\ast} = {\rm Hom}_{R}(V, R)$. 
Let $\{ e_1^{\ast}, \ldots, e_n^{\ast} \} \subset V^{\ast}$ be the dual basis of $\{ e_1, \ldots, e_n \}$.  Set 
\[
I_S^{\perp} = \{ f \in V^{\ast}\otimes_{R} V^{\ast} \mid f(v) = 0 \mbox{ for any } v \in I_S \} 
\]
and 
\[
S^{\perp} = I_S^{\perp} \cap \{ e_i^{\ast}\otimes e_j^{\ast} \in V^{\ast}\otimes_{R} V^{\ast} \mid 1 \le i, j \le n \}.  
\]
Then $I_S^{\perp} =  R\{ S^{\perp} \} \subseteq  V^{\ast}\otimes_{R} V^{\ast}$.  
We define the {\it quadratic dual algebra} $A_{S}^{!}$ of $A_{S}$ by 
\[
A_{S}^{!} = T(V^{\ast})/\langle I_S^{\perp} \rangle,     
\]
where $\langle I_S^{\perp} \rangle$ is the two-sided ideal of $T(V^{\ast})$ generated by $I_S^{\perp}$.  
The quadratic monomial algebra $A_{S}^{!} = \{ V^{\ast}, I_{S}^{\perp} \}$ over $R$ is a graded $R$-algebra with $| e_i^{\ast} | = 1$ ($1\le i \le n$). 

Put $A = A_{S}$, $A^{!} = A_{S}^{!}$, $I=I_S$, and $I^{\perp} = I_{S}^{\perp}$.  
%for a subset $S \subseteq \{ e_i\otimes e_j \in V\otimes_{R} V \mid 1 \le i, j \le n \}$. 
The degree $d$ part 
$A^{!}_d$ of $A^{!}$ can be described by 
\[
A^{!}_d = V^{\ast \otimes d}/(\sum_{i+j=d-2} V^{\ast \otimes i} \otimes_{R}I^{\perp}\otimes_{R} V^{\ast \otimes j}).  
\]
Thus, we can write the dual module $(A^{!}_d)^{\ast} = {\rm Hom}_{R}(A^{!}_d, R)$ of $A^{!}_d$ by   
\[
(A^{!}_d)^{\ast} =\bigcap_{i+j=d-2} V^{\otimes i} \otimes_{R}I\otimes_{R} V^{\otimes j} \subseteq V^{\otimes d}.  
\]
We define the complex $\widehat{K}_{\ast}(A)$ of graded free $A$-bimodules by 
\begin{eqnarray*}
\quad \cdots \longrightarrow A\otimes_{R} (A^{!}_3)^{\ast}\otimes_{R} A \stackrel{\widehat{d}_3}{\longrightarrow} 
A\otimes_{R} (A^{!}_2)^{\ast}\otimes_{R} A \stackrel{\widehat{d}_2}{\longrightarrow} A\otimes_{R} (A^{!}_1)^{\ast}\otimes_{R}A \stackrel{\widehat{d}_1}{\longrightarrow} A\otimes_{R}A \longrightarrow 0,   %\label{eq:widehatK} 
\end{eqnarray*}
where $\widehat{d}_i : A\otimes_{R} (A^{!}_i)^{\ast}\otimes_{R}A \to A\otimes_{R} (A^{!}_{i-1})^{\ast}\otimes_{R} A$  
is defined by 
\begin{multline} 
\sum_{\lambda} a\otimes a_{1}^{\lambda}\otimes a_{2}^{\lambda}\otimes \cdots \otimes a_{i}^{\lambda} \otimes b \label{eq:widehatKdifferential} \\ 
\longmapsto \sum_{\lambda} aa_{1}^{\lambda}\otimes a_{2}^{\lambda}\otimes \cdots \otimes a_{i}^{\lambda}\otimes b 
+(-1)^i \sum_{\lambda} a\otimes a_{1}^{\lambda}\otimes a_{2}^{\lambda}\otimes \cdots \otimes a_{i}^{\lambda} b.  
\end{multline} 
Note that $\widehat{K}_{\ast}(A)$ is a subcomplex of the reduced bar complex $\overline{B}_{\ast}(A,A,A)$. 
We also define the Koszul complex $K_{\ast}(A)$ of graded free $A$-modules by 
\[
\cdots \longrightarrow  
A\otimes_{R} (A^{!}_3)^{\ast} \stackrel{d_3}{\longrightarrow} 
A\otimes_{R} (A^{!}_2)^{\ast} \stackrel{d_2}{\longrightarrow} A\otimes_{R} (A^{!}_1)^{\ast} \stackrel{d_1}{\longrightarrow} A \longrightarrow 0,   
\] 
%\[
%\cdots \to A\otimes_{R} (A^{!}_3)^{\ast} \stackrel{d_3}{\to} A\otimes_{R} (A^{!}_2)^{\ast} \stackrel{d_2}{\to} A\otimes_{R} (A^{!}_1)^{\ast} \stackrel{d_1}{\to} A \stackrel{\epsilon}{\to}  R \to 0,   
%\] 
where $d_i : A\otimes_{R} (A^{!}_i)^{\ast} \to A\otimes_{R} (A^{!}_{i-1})^{\ast}$   is defined by  
\[\sum_{\lambda} a\otimes a_{1}^{\lambda}\otimes a_{2}^{\lambda}\otimes \cdots \otimes a_{n}^{\lambda} 
\longmapsto \sum_{\lambda} aa_{1}^{\lambda}\otimes a_{2}^{\lambda}\otimes \cdots \otimes a_{n}^{\lambda}.
\]   
Then $K_{\ast}(A) \cong \widehat{K}_{\ast}(A)\otimes_{A} R$. 

Using $\displaystyle A = \oplus_{i \ge 0} A_i$, we obtain the following decompositions:   
\[  
\widehat{K}_{d}(A) = A\otimes_{R} (A^{!}_d)^{\ast} \otimes_{R} A = \bigoplus_{l \ge 0} \widehat{K}_{d}(A)_{l}  \quad 
\mbox{ and } \quad 
K_{d}(A) = A\otimes_{R} (A^{!}_d)^{\ast} = \bigoplus_{l \ge 0} K_{d}(A)_{l}  
\]
for $d \ge 0$,   
where 
\[ 
\widehat{K}_{d}(A)_{l} = \bigoplus_{i+d+j=l} A_{i}\otimes_{R} (A^{!}_d)^{\ast} \otimes_{R} A_j   \quad \mbox{ and } \quad 
K_{d}(A)_{l} = \bigoplus_{i+d=l} A_{i}\otimes_{R} (A^{!}_d)^{\ast}.    
\] 
Here $(A^{!}_d)^{\ast} (\subseteq V^{\otimes d})$, $\widehat{K}_{d}(A)_{l}$,  and $K_{d}(A)_{l}$   
are equipped with internal degree $d$, $l$, $l$, respectively.  
Note that $\widehat{K}_{\ast}(A) = \oplus_{l \ge 0} \widehat{K}_{\ast}(A)_{l}$ and $K_{\ast}(A) = \oplus_{l \ge 0} K_{\ast}(A)_{l}$ are the direct sums of subcomplexes. 

By the augmentation map $\epsilon : A \to R$, $R$ can be considered as a left $A$-module (or an $A$-bimodule). The bar-resolution $\widetilde{Bar}_{\ast}(A, R)$ of the left $A$-module $R$ is 
\[
\cdots \stackrel{\partial_3}{\longrightarrow} A\otimes_{R} A_{+}^{\otimes 2} \stackrel{\partial_2}{\longrightarrow} A\otimes_{R} A_{+} \stackrel{\partial_1}{\longrightarrow} A \stackrel{\epsilon}{\longrightarrow} R \longrightarrow 0, 
\] 
where $\widetilde{Bar}_{i}(A, R) = A\otimes_{R} A_{+}^{\otimes i}$ and 
the differential $\partial_i : \widetilde{Bar}_{i}(A, R) \to \widetilde{Bar}_{i-1}(A, R)$ is given by 
\[
\partial_i(a_0\otimes a_1\otimes \cdots \otimes a_i) = \sum_{j=0}^{i-1} (-1)^j a_0\otimes \cdots \otimes a_{j}a_{j+1} \otimes \cdots \otimes a_i. 
\]
Let us consider the cochain complex $Cob^{\ast}(A) = {\rm Hom}_{A}(\widetilde{Bar}_{\ast}(A, R), R)$: for $i \ge 0$, we have 
\[
Cob^{i}(A) = {\rm Hom}_{A}(\widetilde{Bar}_{i}(A, R), R) \cong A_{+}^{\otimes i} \quad  (i>0),  \quad  Cob^{0}(A) = R, 
\]
\[
Cob^{i}(A) = \bigoplus_{j \ge i} Cob^{ij}(A),  \quad Cob^{ij}(A)= \bigoplus_{k_1+\cdots +k_i=j, k_s \ge 1} A_{k_1}\otimes_{R} \cdots \otimes_{R} A_{k_i} \quad (i>0), 
\]
\[
Cob^{0}(A) =  Cob^{00}(A)=R, \quad Cob^{0j}(A)=0 \quad (j>0).   
\]
Then ${\rm Ext}_{A}^i(R, R) = H^{i}(Cob^{\ast}(A)) = \oplus_{j}  {\rm Ext}_{A}^{ij}(R, R)$ and ${\rm Ext}_{A}^{ij}(R, R) = H^{i}(Cob^{\ast, j}(A))$.

\begin{proposition}\label{prop:quadraticmonomialalgebraKoszul} 
Let $A=\{V, I_S\}$ be the quadratic monomial algebra over a commutative ring $R$ associated to 
a subset $S$ of $\{ e_i\otimes e_j \in V\otimes_{R} V \mid 1 \le i, j \le n \}$.  Then 
\begin{enumerate}
\item\label{item:monomialquad-1} ${\rm Ext}^{ij}_{A}(R, R)=0$ for $i\neq j$.  
\item\label{item:monomialquad-2} ${\rm Ext}^{\ast}_{A}(R, R) \cong A^{!}$ as graded $R$-algebras.  
\item\label{item:monomialquad-3} $H_{i}(\widehat{K}_{\ast}(A)) = 0$ $(i>0)$ and $H_0(\widehat{K}_{\ast}(A)) = A$.    
\item\label{item:monomialquad-4} $H_{i}(K_{\ast}(A)) = 0$ $(i>0)$ and $H_0(K_{\ast}(A)) = R$.    
\item\label{item:monomialquad-5} $R$ admits a linear minimal graded free resolution as an $A$-module over $R$, in other words, there exists a graded free resolution over $R$ 
\[
\cdots \longrightarrow P_2 \stackrel{d_2}{\longrightarrow} P_1 \stackrel{d_1}{\longrightarrow} P_0 \stackrel{d_0}{\longrightarrow} R \longrightarrow 0
\] 
such that $d_i : P_i=A\otimes_{R} X_i \to P_{i-1}=A\otimes_{R}X_{i-1}$ can be described by a matrix whose entries are in $A_1$ with respect to $R$-bases of $X_i$ and $X_{i-1}$, where $X_i$ is a graded free $R$-module for each $i \ge 0$. 
\end{enumerate}
\end{proposition}

\proof 
Let $X_i = V^{\otimes i-1}\otimes_{R} I_S \otimes_{R} V^{\otimes d-i-1} \subset V^{\otimes d}$ for $1 \le i \le d-1$. 
Since $I_S = R\{ S \}$, the collection $(X_1, \ldots, X_{d-1})$ is distributive in 
the lattice $(V^{\otimes d}, \cap , \cup)$ consisting of free $R$-submodules by the same discussion in the proof of \cite[Chapter 1 Proposition~7.1]{Quadratic} $(c)\Longrightarrow (a)$. As in the proof of \cite[Chapter 2 Theorem~4.1]{Quadratic}, we see that (\ref{item:monomialquad-1}) ${\rm Ext}^{ij}_{A}(R, R)=0$ for $i\neq j$ and (\ref{item:monomialquad-4})  $H_{i}(K_{\ast}(A)) = 0$ $(i>0)$ and $H_0(K_{\ast}(A)) = R$ hold. 
We can also prove that (\ref{item:monomialquad-2})  ${\rm Ext}^{\ast}_{A}(R, R) \cong A^{!}$ as graded $R$-algebras 
in the same way as \cite[Chapter 1 Proposition~3.1]{Quadratic}. 
Since $K_{\ast}(A)$ gives a linear minimal graded free resolution of $R$, (\ref{item:monomialquad-5}) holds. 

Set ${\mathcal I} = \{ (i, j) \in  [1, \ldots, n]\times [1, \ldots, n] \mid e_i\otimes e_j \in S \}$.  
Let $\{ e'_1, \ldots, e'_n \}$ be a ${\Bbb Z}$-basis of  a free ${\Bbb Z}$-module $V_{{\Bbb Z}}$. 
Put $S_{{\Bbb Z}} = \{ e'_i\otimes e'_j \in V_{{\Bbb Z}}\otimes_{{\Bbb Z}} V_{{\Bbb Z}} \mid (i, j) \in {\mathcal I} \}$, $I_{S_{{\Bbb Z}}} = {\Bbb Z}\{ S_{{\Bbb Z}} \} \subseteq V_{{\Bbb Z}}\otimes_{{\Bbb Z}} V_{{\Bbb Z}}$, and $A_{{\Bbb Z}} =  T(V_{{\Bbb Z}})/\langle I_{S_{{\Bbb Z}}} \rangle$. 
Then $A_{{\Bbb Z}}$ is a quadratic monomial algebra over ${\Bbb Z}$ and  $A = A_{{\Bbb Z}}\otimes_{{\Bbb Z}} R$. 
Note that $\widehat{K}_{\ast}(A) \cong \widehat{K}_{\ast}(A_{{\Bbb Z}})\otimes_{{\Bbb Z}} R$ and $K_{\ast}(A) \cong K_{\ast}(A_{{\Bbb Z}})\otimes_{{\Bbb Z}} R$. 
To show that  (\ref{item:monomialquad-3}) holds for $\widehat{K}_{\ast}(A)$,  
it suffices to prove that  (\ref{item:monomialquad-3}) holds for $\widehat{K}_{\ast}(A_{{\Bbb Z}})$ by K\"{u}nneth theorem. %and $K_{\ast}(A_{{\Bbb Z}})$.  
For any field $k$, $A_{k} = A_{{\Bbb Z}}\otimes_{{\Bbb Z}} k$ is a quadratic monomial algebra over $k$. 
By \cite[Theorem~3.4.6]{Witherspoon},  (\ref{item:monomialquad-3}) holds for 
the complex $\widehat{K}_{\ast}(A_{k}) = \widehat{K}_{\ast}(A_{{\Bbb Z}})\otimes_{{\Bbb Z}} k$ of $A_{k}$-bimodules over $k$. Each internal degree component of $\widehat{K}_{\ast}(A_{{\Bbb Z}})$ is a complex of finitely generated free modules over ${\Bbb Z}$. Using Lemma~\ref{lemma:fieldacyclicisacyclic}, 
(\ref{item:monomialquad-3}) holds for $\widehat{K}_{\ast}(A_{{\Bbb Z}})$.  
(Note that we can prove  (\ref{item:monomialquad-4}) holds for $K_{\ast}(A)$ in the same way.)  
\qed

\bigskip 

The following lemma has been used in the proof of Proposition~\ref{prop:quadraticmonomialalgebraKoszul}. 

\begin{lemma}\label{lemma:fieldacyclicisacyclic} 
Let $C_{\ast}$ be a complex of finitely generated free modules over ${\Bbb Z}$. 
If $C_{\ast}\otimes_{\Bbb Z} k$ is acyclic for any field $k$, then $C_{\ast}$ is acyclic. 
\end{lemma}

\proof 
By the assumption, $H_{i}(C_{\ast}\otimes_{{\Bbb Z}} k) = 0$ for any $i \in {\Bbb Z}$. 
Using K\"{u}nneth theorem, we have an exact sequence
\[
0 \longrightarrow H_{i}(C_{\ast})\otimes_{{\Bbb Z}} k \longrightarrow H_{i}(C_{\ast}\otimes_{{\Bbb Z}} k) \longrightarrow {\rm Tor}_1(H_{i-1}(C_{\ast}), k) \longrightarrow 0,  
\]
which implies that $H_{i}(C_{\ast})\otimes_{{\Bbb Z}} k = 0$ for any field $k$. 
Since $H_{i}(C_{\ast})$ is a finitely generated module over ${\Bbb Z}$, 
$H_{i}(C_{\ast})$ must be $0$. 
Hence $C_{\ast}$ is acyclic.  
\qed 

%\begin{proposition}[{cf. \cite[Chapter 1 Proposition~4.2]{Quadratic}}]
%Let $A = \oplus_{n\ge 0} A_n$ be a graded algebra over a commutative ring $R$. 
%Assume that $A_0=R$ and $A_n$ $(n>0)$ is a free $R$-module of finite rank. 
%Every bounded above complex $F_{\ast}$ of free graded $A$-modules admits a decomposition 
%$F_{\ast} = P_{\ast}\oplus T_{\ast}$ into the direct sum of two subcomplexes of free graded $A$-modules, where 
%$P_{\ast}$ is minimal and $T_{\ast}$ is acyclic. 
%Here, we say that $F_{\ast}$ is {\it bounded above} if $F_m=0$ for $m\ll 0$ and that $P_{\ast}$ is {\it minimal} if $R%\otimes_{A} P_{i} \to R\otimes_{A} P_{i-1}$ is a zero morphism for any $i \in {\Bbb Z}$.  
%\end{proposition}
%
%\proof 
%Assume that 
%\[
%\cdots \to F_{i+1} \stackrel{d_{i+1}}{\to} F_{i} \stackrel{d_{i}}{\to} F_{i-1} \to \cdots 
%\]
%satisfies $\overline{d}_i = d_i \otimes_{A}R = 0$ for all $i\le n$ for some $n \in {\Bbb Z}$. 
%To prove the statement, it suffices to show that there exists a decomposition $F_{\ast} = F'_{\ast}\oplus T_{\ast}$ of into a direct sum of subcomplexes of free graded $A$-modules, where $T_{\ast}$ is an acyclic complex with $T_{i} = 0$ for $i\neq n, n+1$, and the differentials $d'_i$ in $F'_{\ast}$ satisfy $\overline{d'}_i = d'_i\otimes_{A} R=0$ for $i\le n+1$. 
%Set $F_i = X_i\otimes_{R} A$, where $X_i$ is a graded free $R$-module. 
%\qed

\begin{remark}\rm
In Proposition~\ref{prop:quadraticmonomialalgebraKoszul}, we can prove that 
 (\ref{item:monomialquad-4}) implies (\ref{item:monomialquad-3}) for 
$\widehat{K}_{\ast}(A_k)$ with a field $k$ in the following way.  
By \cite[Chapter 1 Proposition~4.2]{Quadratic}, there exists a decomposition $\widehat{K}_{\ast}(A_k) = P_{\ast} \oplus T_{\ast}$ into the direct sum of two subcomplexes of free graded $A_{k}$-bimodules, where $P_{\ast}$ is minimal and $T_{\ast}$ is acyclic. Here we say that $P_{\ast}$ is {\it minimal} if the induced map $P_{i+1}\otimes_{A_k} k \to P_{i}\otimes_{A_k} k$ vanish for any $i \in {\Bbb Z}$.  Using $\widehat{K}_{\ast}(A_k)\otimes_{A_k} k = K_{\ast}(A_k)$, we have 
\[ 
H_{i}(K_{\ast}(A_k)) = H_{i}(P_{\ast}\otimes_{A_k} k) \oplus H_{i}(T_{\ast}\otimes_{A_k} k) = 
\left\{
\begin{array}{cc}
0 & (i \neq 0), \\
k & (i=0). 
\end{array}
\right.
\]
Thus, we obtain 
\[
P_{i}\otimes_{k} k = \left\{
\begin{array}{cc}
0 & (i \neq 0), \\
k & (i=0). 
\end{array}
\right.
\]
Since $P_{\ast}$ has finite-dimensional grading components, $P_i=0$ for $i\neq 0$ by Nakayama's lemma for noncommutative graded algebras (\cite[Chapter 1 Lemma~4.1]{Quadratic}). 
Hence, $H_{i}(\widehat{K}_{\ast}(A_k)) = H_{i}(T_{\ast})= 0$ for $i>0$. We also see that $H_0(\widehat{K}_{\ast}(A_k)) = k$ 
directly. 
\end{remark}

\begin{remark}\rm
Let $\{ e^{\ast}_1, \ldots, e^{\ast}_n \}$ be the dual basis of $V^{\ast}$ of  an $R$-basis 
$\{ e_1, \ldots, e_n \}$ of a free $R$-module $V$.  
Let $\phi : V \to V^{\ast}$ be the $R$-isomorphism defined by 
$e_i \mapsto e_i^{\ast}$ ($1\le i \le n$). Set $A= \{V, I_S\}$ and $A^{!} = \{ V^{\ast}, I_S^{\perp} \}$. 
Then $\phi$ induces an $R$-isomorphism $(A^{!}_{d})^{\ast} \to A^{!}_{d}$ for $d \ge 0$ and an isomorphism of chain complexes of $A$-bimodules over $R$:  
\[
\begin{array}{cccccccccc} 
\cdots \longrightarrow & A\otimes_{R} (A^{!}_2)^{\ast}\otimes_{R} A & \stackrel{\widehat{d}_2}{\longrightarrow} & A\otimes_{R} (A^{!}_1)^{\ast}\otimes_{R}A & \stackrel{\widehat{d}_1}{\longrightarrow} & A\otimes_{R}A & \stackrel{\mu}{\longrightarrow} & A & \longrightarrow & 0   \\
 & \downarrow &  & \downarrow &  & \downarrow &  & \downarrow &  &    \\
\cdots \longrightarrow & A\otimes_{R} A^{!}_2\otimes_{R} A & \longrightarrow & A\otimes_{R} A^{!}_1\otimes_{R}A & \longrightarrow & A\otimes_{R}A & \stackrel{\mu}{\longrightarrow} & A & \longrightarrow & 0,   
\end{array}
\] 
where $\mu : A\otimes_{R} A \to A$ is defined by $\mu(a\otimes b) =ab$ and  the second exact row is given in 
\cite[Theorem~3]{Skoldberg}. In other words, $K_{\ast}(A)$ gives us the free resolution of $A$-bimodules of $A$ 
over $R$ which is isomorphic to the one in \cite[Theorem~3]{Skoldberg}. 
\end{remark}

\begin{theorem}\label{th:HochRisomA!}
Let $A=\{V, I_S\}$ be the monomial quadratic algebra over a commutative ring $R$ associated to 
a subset $S$ of $\{ e_i\otimes e_j \in V\otimes_{R} V \mid 1 \le i, j \le n \}$.  Then 
${\rm HH}^{\ast}(A, R) \cong A^{!}$ as graded $R$-algebras.  
\end{theorem}

\proof 
By Proposition~\ref{prop:quadraticmonomialalgebraKoszul} (\ref{item:monomialquad-3}), $\widehat{K}_{\ast}(A)$ gives us 
a graded free resolution of $A$ as $A$-bimodules over $R$.  By taking ${\rm Hom}_{A^e}(-, R)$ of $\widehat{K}_{\ast}(A)$, 
we obtain a cochain complex whose differentials are all $0$. Hence, ${\rm HH}^{i}(A, R) = H^{i}({\rm Hom}_{A^e}(\widehat{K}_{\ast}(A), R)) = {\rm Hom}_{A^e}(\widehat{K}_{i}(A), R) \cong {\rm Hom}_{R}((A^{!}_{i})^{\ast}, R) \cong A^{!}_{i}$. 

Let us show that ${\rm HH}^{\ast}(A, R) \cong A^{!}$ as graded $R$-algebras. 
Since $\widehat{K}_{\ast}(A)$ is a subcomplex of the reduced bar complex $\overline{B}_{\ast}(A,A,A)$, 
we have a morphism of chain complexes $f = (f_{i})$: 
\[
\begin{array}{cccccccccc} 
\cdots \longrightarrow & A\otimes_{R} (A^{!}_2)^{\ast}\otimes_{R} A & \stackrel{\widehat{d}_2}{\longrightarrow} & A\otimes_{R} (A^{!}_1)^{\ast}\otimes_{R}A & \stackrel{\widehat{d}_1}{\longrightarrow} & A\otimes_{R}A & \stackrel{\mu}{\longrightarrow} & A & \longrightarrow & 0   \\
 & \downarrow f_2 &  & \downarrow f_1 &  & \downarrow f_0 &  & \parallel &  &    \\
\cdots \longrightarrow & A\otimes_{R} \overline{A}^{\otimes 2} \otimes_{R} A & \longrightarrow & A\otimes_{R} \overline{A} \otimes_{R}A & \longrightarrow & A\otimes_{R}A & \stackrel{\mu}{\longrightarrow} & A & \longrightarrow & 0,     
\end{array}
\]
where $\overline{A} = A/R I$ and $I$ is the image of $1 \in R$ under the unit map $R \to A$.  
By taking ${\rm Hom}_{A^{e}}(-, R)$ of chain complexes, we have 
\[
\begin{array}{ccccccc}
0 & \longrightarrow & {\rm Hom}_{A^e}(A\otimes_{R} A, R) & \longrightarrow & {\rm Hom}_{A^e}(A\otimes_{R}\overline{A}\otimes_{R} A, R) & \longrightarrow & \cdots \\  
 & & \downarrow f_0^{\ast} & & \downarrow f_1^{\ast} & &  \\  
0 & \longrightarrow & {\rm Hom}_{A^e}(A\otimes_{R}A, R) & \stackrel{\widehat{d}_1^{\ast}}{\longrightarrow} & {\rm Hom}_{A^e}(A\otimes_{R}(A^{!}_1)^{\ast} \otimes_{R}A, R) & \stackrel{\widehat{d}_2^{\ast}}{\longrightarrow} & \cdots,  \\  
\end{array} 
\] 
which are isomorphic to 
\begin{eqnarray}
\begin{array}{ccccccccc}
0 & \longrightarrow & R & \stackrel{d^0}{\longrightarrow} & {\rm Hom}_{R}(\overline{A}, R) & \stackrel{d^1}{\longrightarrow} & {\rm Hom}_{R}(\overline{A}\otimes_{R} \overline{A}, R)  
& \stackrel{d^2}{\longrightarrow} & \cdots \\  
 & & \quad \downarrow id & & \quad \downarrow f^{\ast}_1& & \quad \downarrow f^{\ast}_2 & &   \label{eq:A!n} \\  
0 & \longrightarrow & R & \stackrel{0}{\longrightarrow} & A^{!}_1& \stackrel{0}{\longrightarrow} & A^{!}_2 & \stackrel{0}{\longrightarrow} &  \cdots.     \\  
\end{array} 
\end{eqnarray} 
This is a quasi-isomorphism of cochain complexes. 
For $g_m = e^{\ast}_{i_1} \otimes \cdots \otimes e^{\ast}_{i_m} \in {\rm Hom}_{R}(\overline{A}^{\otimes m}, R)$ and $g_{d} = e^{\ast}_{j_1} \otimes \cdots \otimes e^{\ast}_{j_{d}} \in {\rm Hom}_{R}(\overline{A}^{\otimes d}, R)$, we have  
$f^{\ast}_m(g_m) = e^{\ast}_{i_1}  \cdots e^{\ast}_{i_m} \in A^{!}_m$ and 
$f^{\ast}_{d}(g_{d}) = e^{\ast}_{j_1}  \cdots e^{\ast}_{j_{d}} \in A^{!}_{d}$.  
The restriction of the product $g_m \cdot g_{d} = e^{\ast}_{i_1} \otimes \cdots \otimes e^{\ast}_{i_m} \otimes e^{\ast}_{j_1} \otimes \cdots \otimes e^{\ast}_{j_{d}} \in {\rm Hom}_{R}(\overline{A}^{\otimes (m+d)}, R)$ to $(A_{m+d}^{!})^{\ast}$ is equal 
to $e^{\ast}_{i_1}  \cdots e^{\ast}_{i_m} e^{\ast}_{j_1}  \cdots e^{\ast}_{j_{d}} \in A^{!}_{m+d}$. 
This implies that ${\rm HH}^{\ast}(A, R) \cong A^{!}$ as graded $R$-algebras. 
\qed

\begin{proposition}\label{prop:HHAL=HHARtensorL}
Let $A=\{V, I_S\}$ be the monomial quadratic algebra over a commutative ring $R$. 
Let $L$ be an $A$-bimodule over $R$. Assume that $A_{+} L = L A_{+} = 0$, where $A_{+} = \oplus_{d>0} A_{d}$. Then ${\rm HH}^{i}(A, L) \cong {\rm HH}^{i}(A, R)\otimes_{R} L \cong A^{!}_{i} \otimes_{R} L$ for $i \ge 0$. 
\end{proposition}

\proof
By Proposition~\ref{prop:quadraticmonomialalgebraKoszul} (\ref{item:monomialquad-3}), $\widehat{K}_{\ast}(A)$ gives us 
a graded free resolution of $A$ as $A$-bimodules over $R$. 
Since $A_{+} L = L A_{+} = 0$, we obtain a cochain complex with zero differential 
\[
0 \longrightarrow {\rm Hom}_{A^e}(A\otimes_{R}A, L) \stackrel{0}{\longrightarrow} {\rm Hom}_{A^e}(A\otimes_{R} (A^{!}_1)^{\ast}\otimes_{R}A, L) \stackrel{0}{\longrightarrow} {\rm Hom}_{A^e}(A\otimes_{R} (A^{!}_2)^{\ast}\otimes_{R}A, L)  \stackrel{0}{\longrightarrow} \cdots   
\]
by taking ${\rm Hom}_{A^e}(-, L)$ of $\widehat{K}_{\ast}(A)$. 
This is isomorphic to 
\[
0 \longrightarrow L \stackrel{0}{\longrightarrow} A^{!}_1\otimes_{R} L \stackrel{0}{\longrightarrow} A^{!}_2\otimes_{R} L \stackrel{0}{\longrightarrow} \cdots.  
\]
Hence
\[
{\rm HH}^{i}(A, L) \cong A^{!}_{i}\otimes_{R} L \cong {\rm HH}^{i}(A, R)\otimes_{R} L
\]
for each $i \ge 0$. 
\qed 

\bigskip

Suppose that $A = \{V, I_S \}$ has finite rank over $R$. 
Set $\overline{A}=A/RI$,
where $I\in A$ is the image of $1\in R$ under the unit map $R\to A$. 
Denote by $\overline{A}_d$ the degree $d$ component of $\overline{A}$. Note that $\overline{A} = \bigoplus_{d\ge1} \overline{A}_{d}$. 
Let $\overline{B}_p(R, A, R) = \overbrace{\overline{A} \otimes_{R} \cdots \otimes_{R} \overline{A}}^{p}$ $(p>0)$ and 
$\overline{B}_0(R, A, R)=R$. For $p>0$, set 
\[
\displaystyle \overline{B}_p(R, A, R)_d = \bigoplus_{a_1+\cdots +a_p =d} \; \overline{A}_{a_1} \otimes_{R} \cdots \otimes_{R} \overline{A}_{a_p}.  
\]
For $p=0$, set 
\[
\overline{B}_0(R, A, R)_d = \left\{
\begin{array}{cc}
R & (d=0), \\
0 & (d\neq 0). 
\end{array}
\right.
\] 
Let $\overline{C}^{\ast}(A, R)$ and $\overline{C}^{\ast, d}(A, R)$ be the cochain complexes defined by  
$\overline{C}^{p}(A, R) = {\rm Hom}_{R}(\overline{B}_p(R, A, R), R)$ and $\overline{C}^{p, d}(A, R) = {\rm Hom}_{R}(\overline{B}_p(R, A, R)_d, R)$, respectively. 
The differential $d^p : \overline{C}^{p}(A, R) \to \overline{C}^{p+1}(A, R)$ (resp. $d^p : \overline{C}^{p, d}(A, R) \to \overline{C}^{p+1, d}(A, R)$) is defined by 
\[
d^{p}(f)(a_1\otimes \cdots \otimes a_{p+1}) = \sum_{j=1}^{p} (-1)^j f(a_1\otimes \cdots \otimes a_{j}a_{j+1}\otimes \cdots \otimes a_{p+1})
\]
for $f \in \overline{C}^{p}(A, R)$ (resp. $f \in \overline{C}^{p, d}(A, R)$). 
Since ${\rm rank}_{R} A < \infty$, we have 
\[
\overline{C}^{p}(A, R) = \bigoplus_{d \in {\Bbb Z}} \overline{C}^{p, d}(A, R). 
\]
Denoting by ${\rm HH}^{p, d}(A, R)$ the $p$-th cohomology of $\overline{C}^{\ast, d}(A, R)$ as in Notation~\ref{notation:Hn,s}, we have  
\[
{\rm HH}^{p}(A, R) = \bigoplus_{d \in {\Bbb Z}} {\rm HH}^{p, d}(A, R).   
\]

\begin{theorem}\label{th:HHndquadratic}  
Suppose that $A = \{V, I_S \}$ has finite rank over $R$. 
For $p\ge 0$, we have 
\[
{\rm HH}^{p, d}(A, R) \cong 
\left\{
\begin{array}{cc}
A^{!}_{p} & (d=p), \\
0 & (d\neq p). 
\end{array} 
\right.
\]
\end{theorem}

\proof
By Theorem~\ref{th:HochRisomA!}, ${\rm HH}^{p}(A, R) \cong A^{!}_{p}$ for $p \ge 0$. 
By (\ref{eq:A!n}), it can be verified that $A^{!}_{p} \subseteq {\rm HH}^{p, p}(A, R)$.
Hence, ${\rm HH}^{p, d}(A, R)=0$ for $d\neq p$ and ${\rm HH}^{p, p}(A, R) \cong A^{!}_{p}$. 
\qed

\subsection{The $R$-algebra structure of ${\rm HH}^{\ast}({\rm N}_m(R), R)$}\label{subsection:HHNmR}
In this subsection, we apply the results in \S\ref{subsection:Quadraticmonomialalgebra} to the case 
$A = {\rm N}_m(R)$ and determine the $R$-algebra structure of  ${\rm HH}^{\ast}({\rm N}_m(R), R)$ for $m\ge 2$.   
Let $R$ be a commutative ring. For $m \ge 2$, set 
\[ 
{\rm N}_m(R) = \left\{
\left(
\begin{array}{ccccc}
a & \ast & \ast & \cdots & \ast \\
0 & a & \ast & \cdots & \ast \\ 
0 & 0 & a & \cdots & \ast \\
\vdots & \vdots & \ddots & \ddots & \vdots \\
0 & 0 & 0 & \cdots & a 
\end{array}
\right) \in {\rm M}_m(R) 
\right\}. 
\]
Putting 
\[
 x_1 = E_{1, 2}, \ x_2 = E_{2, 3}, \ \ldots, \ x_{m-1} = E_{m-1, m} \in {\rm N}_m(R),      
\]
we have an isomorphism ${\rm N}_{m}(R) \cong R\langle x_1, x_2, \ldots, x_{m-1} \rangle/\langle x_i x_j \mid j \neq i+1 \rangle$. 
Let $V = Rx_1\oplus \cdots \oplus Rx_{m-1}$ be a free $R$-module of rank $m-1$. 
Set $S= \{ x_i \otimes x_j \mid j\neq i+1 \} \subset V\otimes V$. 
Then ${\rm N}_m(R) = \{ V, I_S \}$ is a quadratic monomial algebra over $R$ with $| x_i | =1$. 

We define the two-sided ideal $J({\rm N}_m(R))$ of ${\rm N}_m(R)$ over $R$ by 
\[
J({\rm N}_m(R)) = \{ (a_{ij}) \in {\rm N}_m(R) \mid a_{11} = a_{22} = \cdots = a_{nn} = 0 \}.  
\] 
Note that $J({\rm N}_m(R))= {\rm N}_m(R)_{+} = \oplus_{i > 0} {\rm N}_m(R)_i$. 
When $R$ is a field, $J({\rm N}_m(R))$ is the Jacobson radical of ${\rm N}_m(R)$. 
We denote $(I_m \bmod J({\rm N}_{m}(R)))\in {\rm N}_{m}(R)/J({\rm N}_{m}(R))$ by $e$.  Then 
we see that ${\rm N}_m(R)/J({\rm N}_m(R)) = Re \cong R$ as an ${\rm N}_m(R)$-bimodule over $R$. 
In the sequel, the ${\rm N}_m(R)$-bimodule $R$ over $R$ means ${\rm N}_m(R)/J({\rm N}_m(R)) = Re$. 

Put ${\rm N} = {\rm N}_m(R)$ and ${\rm J}=J({\rm N}_m(R))$.  Denote by ${\rm N}^{!} = {\rm N}_m(R)^{!}$ 
the quadratic dual algebra of ${\rm N}$. 
Note that $S^{\perp} =\{  x^{\ast}_i \otimes x^{\ast}_{i+1} \mid 1 \le i \le m-2 \}$, where 
$\{ x_1^{\ast}, \ldots, x_{m-1}^{\ast} \}$ is the dual basis of $V^{\ast} = {\rm Hom}_{R}(V, R)$ of the $R$-basis $\{ x_1, \ldots, x_{m-1} \}$ of $V$.  
Setting $y_i = x^{\ast}_i$ $(1\le i \le m-1)$, we can write ${\rm N}^{!} = R\langle y_1, y_2, \ldots, y_{m-1} \rangle / \langle y_i y_{i+1} \mid  1 \le i < m-1 \rangle$ with $| y_i | =1$. 
%By setting $| y_i | =1$, ${\rm N}^{!}$ is a graded $R$-algebra. 
Let us denote by ${\rm N}^{!}_{n}$ the homogeneous part of ${\rm N}^{!}$ of degree $n$.  We also denote by ${\mathcal B}({\rm N}^{!}_{n})$ the $R$-basis of ${\rm N}^{!}_{n}$  consisting of monomials of degree $n$ in $\{ y_1, \ldots, y_{m-1} \}$. (Set ${\mathcal B}({\rm N}^{!}_{0}) = \{ 1 \}$.)  Put ${\mathcal B}({\rm N}^{!}) = \cup_{n=0}^{\infty} {\mathcal B}({\rm N}^{!}_{n})$. 

\begin{theorem}\label{th:HHNR=N!}
We have an isomorphism 
\[
{\rm HH}^{\ast}({\rm N}, R) \cong {\rm N}^{!} \cong R\langle y_1, y_2, \ldots, y_{m-1} \rangle/\langle y_iy_{i+1} \mid 1 \le i \le m-2 \rangle  
\] 
of graded $R$-algebras, where $|y_i|=1$ for $1 \le i \le m-1$.  
\end{theorem}

\proof
The statement follows from Theorem~\ref{th:HochRisomA!}.  
\qed

\begin{proposition}\label{prop:JL=LJ=0case} 
Let $L$ be an ${\rm N}$-bimodule over $R$. Assume that ${\rm J} L = L {\rm J} = 0$. Then ${\rm HH}^{n}({\rm N}, L) \cong {\rm HH}^{n}({\rm N}, R)\otimes_{R} L \cong {\rm N}^{!}_{n} \otimes_{R} L$ for $n \ge 0$. 
\end{proposition}

\proof
The statement follows from Proposition~\ref{prop:HHAL=HHARtensorL}. 
\qed 

%\begin{corollary}\label{cor:JL=LJ=0case} 
%Let $L$ be an ${\rm N}$-bimodule over $R$. Assume that ${\rm J} L = L {\rm J} = 0$. Then ${\rm HH}^{n}({\rm N}, L) \cong {\rm HH}^{n}({\rm N}, R)\otimes_{R} L \cong {\rm N}^{!}_{n} \otimes_{R} L$ for $n \ge 0$. 
%\end{corollary}
%
%\proof 
%Let us consider the cochain complex 
%\[
%0 \to {\rm Hom}_{{\rm N}^{e}}({\rm N}\otimes_{R} {\rm N}, L) \stackrel{\delta'^0}{\to} {\rm Hom}_{{\rm N}^{e}}({\rm N}\otimes_{R}{\rm N}^{!}_1\otimes_{R} {\rm N}, L) \stackrel{\delta'^1}{\to} {\rm Hom}_{{\rm N}^{e}}({\rm N}\otimes_{R}{\rm N}^{!}_2\otimes_{R} {\rm N}, L) \to \cdots. 
%\] 
%In a similar way as Lemma~\ref{lemma:differential0}, we can prove $\delta'_n = 0$ for $n \ge 0$ since ${\rm J} L = L {\rm J} = 0$.   Hence, 
%\begin{eqnarray*}
%{\rm HH}^{n}({\rm N}, L) & = & 
%H^n({\rm Hom}_{{\rm N}^{e}}({\rm N}\otimes_{R}{\rm N}^{!}_{\ast}\otimes_{R} {\rm N}, L)) \\ 
%& \cong  & {\rm Hom}_{{\rm N}^e}({\rm N}\otimes_{R}{\rm N}^{!}_{n}\otimes_{R} {\rm N}, L). \\ 
%& \cong & {\rm Hom}_{R}({\rm N}^{!}_{n}, L) \\
%& \cong & {\rm N}^{!}_n \otimes_{R} L  \\
%& \cong & {\rm HH}^{n}({\rm N}, R)\otimes_{R} L. 
%\end{eqnarray*}
%This completes the proof. 
%\qed

\bigskip

Let $\overline{{\rm N}} = {\rm N}/RI_m$. Denote by $\overline{{\rm N}}_d$ the degree $d$ component of $\overline{{\rm N}}$. Note that $\overline{{\rm N}} = \bigoplus_{d=1}^{m-1} \overline{{\rm N}}_{d}$. Set
\[
\overline{B}_p({\rm N}, {\rm N}, {\rm N}) = {\rm N}\otimes_{R} \overbrace{\overline{{\rm N}} \otimes_{R} \cdots \otimes_{R} \overline{{\rm N}}}^{p} \otimes_{R} {\rm N} 
\]
for $p \ge 0$. 
Let $\overline{B}_p(R, {\rm N}, R) = \overbrace{\overline{{\rm N}} \otimes_{R} \cdots \otimes_{R} \overline{{\rm N}}}^{p}$ and $\displaystyle \overline{B}_p(R, {\rm N}, R)_d = \bigoplus_{a_1+\cdots +a_p =d} \; \overline{{\rm N}}_{a_1} \otimes_{R} \cdots \otimes_{R} \overline{{\rm N}}_{a_p}$.  Let $\overline{C}^{\ast}({\rm N}, R)$ and $\overline{C}^{\ast, d}({\rm N}, R)$ be the cochain complexes defined by  
$\overline{C}^{p}({\rm N}, R) = {\rm Hom}_{R}(\overline{B}_p(R, {\rm N}, R), R)$ and $\overline{C}^{p, d}({\rm N}, R) = {\rm Hom}_{R}(\overline{B}_p(R, {\rm N}, R)_d, R)$, respectively. 
Denoting by ${\rm HH}^{p, d}({\rm N}, R)$ the $p$-th cohomology of $\overline{C}^{\ast, d}({\rm N}, R)$ as in Notation~\ref{notation:Hn,s}, we have  
\[
{\rm HH}^{p}({\rm N}, R) = \bigoplus_{d=0}^{p(m-1)} {\rm HH}^{p, d}({\rm N}, R).   
\]

\begin{theorem}\label{th:HHndNR}
For $n\ge 0$, we have 
\[
{\rm HH}^{n, d}({\rm N}, R) \cong 
\left\{
\begin{array}{cc}
{\rm N}^{!}_n & (d=n), \\
0 & (d\neq n). 
\end{array} 
\right.
\]
\end{theorem}

\proof
The statement follows from Theorem~\ref{th:HHndquadratic}. 
\qed 

\subsection{Several results on $\varphi(n)$}\label{subsection:resultsonvarphi} In this subsection, we prove several results on $\varphi(n) = {\rm rank}_{R} {\rm N}^{!}_{n}$. These results will be used in \S\ref{section:HH-N-M/N} and \S\ref{section:HH-N-N} for describing the ranks of Hochschild cohomology over $R$. 

\begin{definition}\label{def:varphi}\rm 
For the $R$-algebra ${\rm N}^{!} = \oplus_{n=0}^{\infty} {\rm N}^{!}_n$, set  
\begin{eqnarray}
\varphi(n) & = & {\rm rank}_{R} {\rm N}^{!}_n = \sharp {\mathcal B}({\rm N}^{!}_{n}), \label{eq:varphi} \\ 
f^{!}(t) & = & \sum_{n=0}^{\infty} ({\rm rank}_{R} {\rm N}^{!}_n) \; t^n.  \label{eq:f^{!}} 
\end{eqnarray}
Note that $\varphi(n) =  0$ for $n<0$. 
\end{definition}

\begin{definition}\rm 
For $1 \le i \le m-1$, we define 
\begin{eqnarray*}
{\mathcal B}({\rm N}^{!}_n)(i) & = & 
\left\{ 
\begin{array}{lc}
\emptyset & (n=0), \\
\{ y_i f \in {\mathcal B}({\rm N}^{!}_n) \mid f \in {\mathcal B}({\rm N}^{!}_{n-1}) \} & (n>0),  
\end{array} 
\right. \\ 
\psi_{i}(n) & = & \sharp {\mathcal B}({\rm N}^{!}_n)(i) \qquad (n \ge 0).   
\end{eqnarray*}  
Note that $\psi_{i}(0) = 0$ and $\psi_{i}(1) = 1$ for $1 \le i \le m-1$ and that 
$\varphi(n) =  \sum_{i=1}^{m-1}\psi_{i}(n)$ for $n \ge 1$. 
\end{definition}

\begin{proposition}\label{prop:psisequence}
We have 
\[
\left(
\begin{array}{c}
\psi_{1}(n+1) \\
\psi_{2}(n+1) \\
\vdots \\
\psi_{m-1}(n+1) 
\end{array}
\right) = 
\left(
\begin{array}{ccccc}
1 & 0 & 1 & \cdots & 1\\
1 & 1 & 0 & \cdots & 1\\
%1 & 1 & 1 & \ddots & 1\\
\vdots & \vdots & \cdots & \ddots & \vdots \\
1 & 1 & 1 & \cdots & 0\\
1 & 1 & 1 & \cdots & 1\\
\end{array}
\right)
\left(
\begin{array}{c}
\psi_{1}(n) \\
\psi_{2}(n) \\
\vdots \\
\psi_{m-1}(n) 
\end{array}
\right)  \quad (n \ge 1) 
\] 
and 
\[  
\left(
\begin{array}{c}
\psi_{1}(1) \\
\psi_{2}(1) \\
\vdots \\
\psi_{m-1}(1) 
\end{array}
\right) = 
\left(
\begin{array}{c}
1 \\
1 \\
\vdots \\
1 
\end{array}
\right). 
\]
\end{proposition} 

\proof
Recall that $\psi_{i}(1) = 1$ for $1\le i \le m-1$.  Note that $y_iy_{i+1}=0$ in ${\rm N}^{!}$ for $1 \le i \le m-2$. 
Since 
\[{\mathcal B}({\rm N}^{!}_{n+1})(i) = y_i {\mathcal B}({\rm N}^{!}_n)(1) \coprod \cdots \coprod 
 y_i {\mathcal B}({\rm N}^{!}_n)(i) \coprod  y_i {\mathcal B}({\rm N}^{!}_n)(i+2) \coprod 
\cdots \coprod y_{i} {\mathcal B}({\rm N}^{!}_n)(m-1)
\] for  $1\le i \le m-1$ and $n \ge 1$, we can verify the statement. 
\qed

\begin{corollary}\label{cor:varphi}For $m \ge 2$, we have 
$\varphi(0)=1, \varphi(1)=m-1$, and  $\varphi(2)=m^2-3m+3$.  
In particular, $\varphi(n)>0$ for $n \ge 0$.  
\end{corollary}

\proof 
By Proposition~\ref{prop:psisequence}, we can calculate $\varphi(n)$ for $0 \le n \le 2$. 
(We can also calculate them directly.)  
We also easily see that $\psi_i(n) > 0$ for $1 \le i \le m-1$ and $n>0$ by induction.  
Hence, $\varphi(n) = \sum_{i=1}^{m-1} \psi_i(n) > 0$ for $n>0$.  
\qed 

\bigskip 

Let ${\rm N}_d$ be the degree $d$ component of ${\rm N}$. 
Let us define 
\begin{eqnarray}
f(t)  =  \sum_{n=0}^{\infty} \left({\rm rank}_{R} {\rm N}_{n}  \right) t^n. \label{eq:f} 
\end{eqnarray}  
Note that 
\begin{eqnarray}
f(t) = 1+ \sum_{k=1}^{m-1}(m-k)t^k = 1 + \sum_{k=1}^{m-1} k t^{m-k}.  \label{eq:explicitf} 
\end{eqnarray} 

By \cite[Chapter 2 Corollary~4.3]{Quadratic}, the quadratic monomial algebra ${\rm N}_m(R)$ over $R$ is Koszul when 
$R$ is a field. 
Hence, we have the following formula. 

\begin{proposition}\label{prop:ff^{!}} 
For $m \ge 2$, 
\[
f^{!}(t) f(-t) = 1. 
\]
In particular, 
\[
f^{!}(t) = \frac{1}{f(-t)} = \frac{1}{\displaystyle 1+\sum_{k=1}^{m-1}(-1)^{k} (m-k) t^{k}}. 
\]
\end{proposition} 
\proof 
By \cite[Chapter 2 Corollary~2.2]{Quadratic} or \cite[{Theorem~3.5.1}]{Loday-Vallette}, we see that $f^{!}(t) f(-t) = 1$ if $R$ is a field. The formula can be also proved for any commutative ring $R$ since $f(t)$ and $f^{!}(t)$ are common to any $R$.  
\qed 

\bigskip 

The following lemma will be used in \S\ref{section:HH-N-N}. 

\begin{lemma}\label{lemma:varphi-relation}
For $q>0$,
we have
\[ \varphi(q)=
   \sum_{r=1}^{m-1}(-1)^{r-1}(m-r)\varphi(q-r).\]
\end{lemma}

\proof
Since $f^{!}(x) f(-x)=1$ by Proposition~\ref{prop:ff^{!}}, 
we have
\[ \left(\sum_{q\ge 0}\varphi(q)x^q\right)\cdot
   \left(1+\sum_{k=1}^{m-1}(m-k)
   (-x)^k\right)=1.\]
Comparing the coefficients of $x^q$ in both sides,
we obtain
\[ \varphi(q)-(m-1)\varphi(q-1)+(m-2)\varphi(q-2)-\cdots
   +(-1)^{m-1}{\varphi(q-(m-1))}=0.\]
This completes the proof.    
\qed

\bigskip 

We have another formula for $\varphi(n)$. 

\begin{proposition}\label{prop:varphicomb} 
For $n \ge 0$, we have 
\begin{eqnarray}
\varphi(n) = (-1)^n \sum_{r \ge 0}(-1)^r \sum_{(a_1, \ldots, a_r)} (m-a_1)\cdots (m-a_r), \label{eq:varphicomb} 
\end{eqnarray}
where the second sum ranges over the $r$-tuples $(a_1, \ldots, a_r)$ of integers such that $1 \le a_i <m$ for $1\le i \le r$ and $a_1+\cdots + a_r = n$.  
\end{proposition}

\proof
Recall $\displaystyle \overline{B}_r(R, {\rm N}, R)_d = \oplus_{a_1+\cdots +a_r =d} \; \overline{{\rm N}}_{a_1} \otimes_{R} \cdots \otimes_{R} \overline{{\rm N}}_{a_r}$. Since ${\rm rank}_{R} \overline{{\rm N}}_{a} = m-a$ for $1 \le a \le m$ and $\overline{{\rm N}}_{0}=0$,  
\[
{\rm rank}_{R} \overline{B}_r(R, {\rm N}, R)_d = \sum (m-a_1)\cdots (m-a_r), 
\]
where the sum ranges over the $r$-tuples $(a_1, \ldots, a_r)$ of integers such that $1 \le a_i <m$ for $1 \le i \le r$ and $a_1+\cdots +a_r=d$. 
The rank of 
$\overline{C}^{r, d}({\rm N}, R) = {\rm Hom}_{R}(\overline{B}_r(R, {\rm N}, R)_d, R)$ is equal to 
${\rm rank}_{R} \overline{B}_r(R, {\rm N}, R)_d$.  
Note that the cochain complex 
\begin{eqnarray}
0 \longrightarrow \overline{C}^{0, d}({\rm N}, R) \longrightarrow \overline{C}^{1, d}({\rm N}, R) \longrightarrow \overline{C}^{2, d}({\rm N}, R) \longrightarrow \cdots 
\label{eq:barCrdNR} 
\end{eqnarray}
satisfies that $\overline{C}^{r, d}({\rm N}, R) = 0$ for $r>d$. By taking the Euler characteristic of 
(\ref{eq:barCrdNR}) and using ${\rm HH}^{r, d}({\rm N}, R) \cong H^{r}(\overline{C}^{\ast, d}({\rm N}, R))$, we obtain 
\[
\sum_{r \ge 0} (-1)^r {\rm rank}_{R} {\rm HH}^{r, d}({\rm N}, R) = \sum_{r \ge 0} (-1)^r {\rm rank}_{R} \overline{C}^{r, d}({\rm N}, R).    
\]
By Theorem~\ref{th:HHndNR}, ${\rm HH}^{r, d}({\rm N}, R)=0$ for $r\neq d$ and ${\rm HH}^{d, d}({\rm N}, R)\cong {\rm N}^{!}_{d}$.  Hence, 
\[
(-1)^d {\rm rank}_{R} {\rm N}^{!}_{d} = \sum_{r \ge 0} (-1)^r {\rm rank}_{R} \overline{C}^{r, d}({\rm N}, R).    
\]
Thus, we have
\begin{eqnarray*}
\varphi(n) & = & {\rm rank}_{R} {\rm N}^{!}_{n} \\ 
 & = & (-1)^n \sum_{r \ge 0} (-1)^r {\rm rank}_{R} \overline{C}^{r, n}({\rm N}, R)  \\
 & = &  (-1)^n \sum_{r \ge 0}(-1)^r \sum_{(a_1, \ldots, a_r)} (m-a_1)\cdots (m-a_r),  
\end{eqnarray*}
which is what we wanted to show. 
\qed 

\bigskip 

By taking account of the appearance of the term $(m-1)^{i_1}(m-2)^{i_2}\cdots (m-(m-1))^{i_{m-1}}$ in the right hand side of (\ref{eq:varphicomb}) in Proposition~\ref{prop:varphicomb}, we obtain the following corollary. 

\begin{corollary}\label{cor:varphicombination}
For $n \ge 0$, we have 
\begin{eqnarray} \quad\quad  
\varphi(n) = (-1)^n \sum \frac{(i_1+i_2+\cdots + i_{m-1})!}{i_1!i_2!\cdots i_{m-1}!} (1-m)^{i_1}(2-m)^{i_2}\cdots ((m-1)-m)^{i_{m-1}},   \label{eq:varphicomb2} 
\end{eqnarray}
where the sum ranges over the $(m-1)$-tuples $(i_1, i_2, \ldots, i_{m-1})$ of non-negative integers such that $i_1+2i_2+\cdots +(m-1)i_{m-1} =n$.   
\end{corollary}

\begin{remark}\rm 
Corollary~\ref{cor:varphicombination} can be also proved  by using Proposition~\ref{prop:ff^{!}}. Indeed, putting $y= - \sum_{k=1}^{m-1} (-1)^k (m-k) t^k = \sum_{k=1}^{m-1} (k-m)(-t)^k$, we have 
\[ 
f^{!}(t) = \frac{1}{f(-t)} = \frac{1}{1-y}  =  1+ y + y^2+ y^3+ \cdots 
\]  
by Proposition~\ref{prop:ff^{!}}. Then $f^{!}(t)$ equals to 
\begin{eqnarray*}
\sum_{r \ge 0} y^r & = & \sum_{r \ge 0} \left\{ \sum_{k=1}^{m-1} (k-m)(-t)^k \right\}^r \\
& = &  \sum_{n \ge 0} (-1)^n t^n \sum  \frac{(i_1+i_2+\cdots + i_{m-1})!}{i_1!i_2!\cdots i_{m-1}!} (1-m)^{i_1}(2-m)^{i_2}\cdots ((m-1)-m)^{i_{m-1}}, 
\end{eqnarray*}
where the second sum ranges over the $(m-1)$-tuples $(i_1, i_2, \ldots, i_{m-1})$ of non-negative integers such that $i_1+2i_2+\cdots +(m-1)i_{m-1} =n$.   Comparing the coefficients of $t^n$ on both sides, we obtain (\ref{eq:varphicomb2}). 
\end{remark}

\section{${\rm HH}^{\ast}({\rm N}_m(R), {\rm M}_m(R)/{\rm N}_m(R))$}
\label{section:HH-N-M/N}

In this section, we determine the $R$-module structure of  
${\rm HH}^{\ast}({\rm N}_m(R), {\rm M}_m(R)/{\rm N}_m(R))$ for $m \ge 3$. (The case $m=2$ will be discussed in \S\ref{section:m=2}.) In \S\ref{subsection:sssubquotient}, we construct a spectral sequence converging to ${\rm HH}^{\ast}({\rm N}_m(R), {\rm M}_m(R)/{\rm N}_m(R))$, which collapses from the $E_2$-page. In \S\ref{subsection:calHHM/N}, we show that ${\rm HH}^{\ast}({\rm N}_m(R), {\rm M}_m(R)/{\rm N}_m(R))$ is a free $R$-module by calculating $E_{2}^{p, q}$.  We also calculate the rank of the free $R$-module ${\rm HH}^{\ast}({\rm N}_m(R), {\rm M}_m(R)/{\rm N}_m(R))$ by using $\varphi(n)$ defined in Definition~\ref{def:varphi}.  

\subsection{Spectral sequences for subquotients of ${\rm M}_m(R)$}\label{subsection:sssubquotient} 

In this subsection
we introduce a $\mathbb{Z}$-grading 
on the matrix algebra ${\rm M}_m(R)$.
Using this grading,
we construct a spectral sequence
of $\mathbb{Z}$-graded $R$-modules
converging to the Hochschild cohomology
${\rm HH}^*({\rm N}_m(R),\widehat{M})$,
where $\widehat{M}$ is a subquotient
of the $\mathbb{Z}$-graded $R$-module
${\rm M}_m(R)$.
Furthermore,
we show that this spectral sequence
collapses from the $E_2$-page.

In this subsection
we work in the abelian category 
of $\mathbb{Z}$-graded $R$-modules.
Let $m \ge 2$ and ${\rm M}={\rm M}_m(R)$.
First, we introduce a grading on ${\rm M}$.
%Since ${\rm M}$ is a finitely generated free $R$-module,
We can choose a basis 
\[ \{E_{i,j}|\ 1\le i,j\le m\}\]
over $R$.
We define a (homological) degree $r$ component
of ${\rm M}$ by 
\[ {\rm M}_r=\bigoplus_{j-i=r}R\{E_{i,j}\}.\]
Then we can verify that
${\rm M}=\bigoplus_{r\in\mathbb{Z}}{\rm M}_r$ 
is a $\mathbb{Z}$-graded associative algebra over $R$.

Let ${\rm N}={\rm N}_m(R)$.
We can easily see that ${\rm N}$ is 
a $\mathbb{Z}$-graded subalgebra of ${\rm M}$.
For a $\mathbb{Z}$-graded
${\rm N}$-bimodule $L$ over $R$,
we let $C^*({\rm N},L)$ be the Hochschild cochain complex.
We have
\[ \begin{array}{rcl}
    C^p({\rm N},L)&\cong&
    {\rm Hom}_R({\rm N}^{\otimes p},L)\\[2mm]
    &\cong&
    ({\rm N}^*)^{\otimes p}\otimes_R L,\\
   \end{array} \]
%which is a finitely generated free module over $R$,
where ${\rm N}^*={\rm Hom}_R({\rm N},R)$.
We denote by
\[ C^{p,s}({\rm N},L) \]
the $R$-submodule of $C^p({\rm N},L)$ 
of (cohomological) degree $s$.
For example,
when $L={\rm N}$, 
we have
\[ E_{1,2}\in C^{0,-1}({\rm N},{\rm N}),\quad
   E_{2,3}^*\otimes I_m\in C^{1,1}({\rm N},{\rm N}),\quad
   E_{1,2}^*\otimes E_{1,3}^*\otimes E_{2,3}
   \in C^{2,2}({\rm N},{\rm N}),\]
where $\{I_m^*\} \cup \{ E_{i,j}^* \mid i<j \}$ is the dual basis of ${\rm N}^*$
with respect to the $R$-basis $\{I_m\} \cup \{E_{i,j} \mid i<j \}$ of ${\rm N}$.

Since the differential
$d: C^p({\rm N},L)\to C^{p+1}({\rm N},L)$
preserves the grading,
we have an isomorphism
\[ C^*({\rm N},L)\cong
   \bigoplus_{s\in\mathbb{Z}}C^{*,s}({\rm N},L)\]
of cochain complexes of $R$-modules.
Thus, we can regard $C^*({\rm N},L)$
as a cochain complex of $\mathbb{Z}$-graded $R$-modules. 
We set 
\begin{eqnarray}
{\rm HH}^{n, s}({\rm N}, L) = H^{n}(C^{\ast, s}({\rm N}, L))  \label{eq:defofHn,s} 
\end{eqnarray} 
as in Notation~\ref{notation:Hn,s}.  

Let 
\[ {\rm J}=\bigoplus_{j-i>0}R\{E_{i,j}\} \]
be the two-sided ideal of ${\rm N}$
consisting of upper triangular matrices
with zero diagonal entries. In other words, ${\rm J}$ coincides with $J({\rm N}_m(R))$ defined 
in \S\ref{subsection:HHNmR}. 
By (\ref{eq:JpM}), we have the filtration
$\{\overline{\rm J}^pL\}_{p\ge 0}$
on $M$,
which induces a 
filtration on the cochain complex
$C^*({\rm N},L)$.
Since
${\rm J}$ is a homogeneous two-sided ideal of ${\rm N}$,
we can verify that this filtration
is compatible with the grading.
Thus,
we obtain a spectral sequence
\[ {}^JE_1^{p,q}({\rm N},L)
   \Longrightarrow {\rm HH}^{p+q}({\rm N},L) \]
of $\mathbb{Z}$-graded $R$-modules by Proposition~\ref{prop:sstrigrading}, where 
\[
{}^JE_1^{p,q}({\rm N},L)
 \cong {\rm HH}^{p+q}({\rm N}, \overline{\rm J}^{p}L/\overline{\rm J}^{p+1}L).   
\]  
We note that
\[ {}^JE_r^{p,q}({\rm N},L)=0 \]
for $p\ge 2m-1$ since $\overline{\rm J}^{2m-1}L=0$.
%}

Next,
we consider the following situation.
Let $M''\subset M'\subset {\rm M}={\rm M}_m(R)$
be $\mathbb{Z}$-graded ${\rm N}$-sub-bimodules
over $R$.
We would like to construct another spectral
sequence converging to 
the Hochschild cohomology
%${\rm HH}^*({\rm N},\mathsf{M})$,
${\rm HH}^*({\rm N},\widehat{M})$,
where $\widehat{M}=M'/M''$.

For this purpose,
we define a filtration $\{F^p{\rm M}\}$
on the $\mathbb{Z}$-graded $R$-module
${\rm M}={\rm M}_m(R)$ by
reindexing the filtration
$\{\overline{\rm J}^p{\rm M}\}$
as follows
\[ F^{p-(m-1)}{\rm M}=\overline{{\rm J}}^p{\rm M} = \sum_{a+b=p} {\rm J}^a{\rm M}{\rm J}^b. \]
%where ${\rm J}={\rm J}_m(R)$ is the Jacobson ideal
%of ${\rm N}={\rm N}_m(R)$.
Then we have
\[ {\rm M}=F^{-(m-1)}{\rm M}\supset
           F^{-(m-1)+1}{\rm M}\supset\cdots
           \supset F^{m-1}{\rm M}\supset
           F^m{\rm M}=0.\]
Note that
\[ {\rm Gr}^p({\rm M})\cong
   R\{E_{i,j}|\ j-i=p\}.\]

Using this filtration on ${\rm M}$,
we define a filtration $\{F^p\widehat{M}\}$
on $\widehat{M}$ to be the induced filtration
\[ F^p\widehat{M}=((M'\cap F^p{\rm M})+M'')/M''.\]
Note that
${\rm Gr}^p(\widehat{M}) = F^p\widehat{M}/F^{p+1}\widehat{M}$ is a subquotient
of ${\rm Gr}^p({\rm M})$.
Using this filtration on $\widehat{M}$,
we obtain the following proposition.

\begin{proposition}\label{prop:spectral-sequence-matrix-filtration}
There is a spectral sequence
\[ {}^ME_1^{p,q}({\rm N},\widehat{M})\Longrightarrow
   H^{p+q}({\rm N},\widehat{M})\] 
of $\mathbb{Z}$-graded $R$-modules,
where
\[ {}^ME_1^{p,q}({\rm N},\widehat{M})\cong
    {\rm HH}^{p+q}({\rm N},{\rm Gr}^p(\widehat{M})).\]
We have 
${}^ME_r^{p,q}({\rm N},\widehat{M})=0$ 
unless $-(m-1)\le p\le m-1$.
\end{proposition}

We shall show that the spectral sequence
$\{{}^ME_r^{p,q}({\rm N},\widehat{M})\}_{r\ge 1}$
collapses from the $E_2$-page and 
that there is no extension problem. 
Recall the degree $s$ component ${}^ME_r^{p,q,s}({\rm N},\widehat{M})$ 
of $\{{}^ME_r^{p,q}({\rm N},\widehat{M})\}_{r\ge 1}$ in Proposition~\ref{prop:sstrigrading}.  

\begin{lemma}\label{lemma:graded_ss_vanishing-simple-version}
%If $s\neq q$, 
We have 
%then %${}^ME_1^{p,q,s}({\rm N},\widehat{M})=0$.
%\if0
\[ {}^ME_1^{p,q,s}({\rm N},\widehat{M})\cong
   \left\{\begin{array}{cl}
           {\rm HH}^{p+q}({\rm N},
            {\rm Gr}^p(\widehat{M})) & (s=q),\\[2mm]
           0                  & (s\neq q).\\
          \end{array}\right.
\]
%\fi
\end{lemma}

\proof
Since ${\rm J}\cdot {\rm Gr}^p(\widehat{M})
={\rm Gr}^p(\widehat{M})\cdot {\rm J}=0$
for each $p\ge 0$,
we have an isomorphism
\[ {}^ME_1^{p,q}({\rm N},\widehat{M})
    %\cong
    %{\rm HH}^{p+q}({\rm N},\mathrm{Gr}^p(\mathsf{M}))\\[2mm]
    \cong
    {\rm HH}^{p+q}({\rm N},R)\otimes_R \mathrm{Gr}^p(\widehat{M})%\\
\]   
by Proposition~\ref{prop:JL=LJ=0case}. %Corollary~\ref{cor:JL=LJ=0case}. 
Recall that 
\[
{\rm HH}^{n, d}({\rm N}, R) \cong 
\left\{
\begin{array}{cc}
{\rm N}^{!}_n & (d=n), \\
0 & (d\neq n) 
\end{array} 
\right.
\]
by Theorem~\ref{th:HHndNR}. 
%Recall that 
%\[ {\rm HH}^*({\rm N},R)\cong R\langle y_1,\ldots,y_{m-1}\rangle
%    /(y_iy_{i+1}|\ 1\le i\le m-2),\]
%where $y_i \in {\rm HH}^{1}({\rm N},R) = {\rm HH}^{1,1}({\rm N},R)$ 
%for $1\le i\le m-1$ by Theorems~\ref{th:HHNR=N!} and \ref{th:HHndNR}. 
The statement follows from the fact that 
$\mathrm{Gr}^p(\widehat{M})$
is a subquotient of ${\rm Gr}^p({\rm M})=R\{E_{i,j}|\ j-i=p\}$.
\qed

\begin{theorem}\label{theorem:collapse-E-2-term-no-extension}
The spectral sequence
${}^ME_1^{p,q}({\rm N},\widehat{M})
\Longrightarrow {\rm HH}^{p+q}({\rm N},\widehat{M})$
collapses from the $E_2$-page and
there is no extension problem.
\end{theorem}

\proof
By Lemmas~\ref{lemma:ss-collapse-internal-degree} and \ref{lemma:graded_ss_vanishing-simple-version}, 
we see that the spectral sequence collapses
from the $E_2$-page.

We shall show that there is no extension problem.
We have a filtration
$\{F^p{\rm HH}^n({\rm N},\widehat{M})\}$ 
on the $\mathbb{Z}$-graded $R$-module
${\rm HH}^n({\rm N},\widehat{M})$
given by
\[ F^p{\rm HH}^n({\rm N},\widehat{M})=
   {\rm Im}({\rm HH}^n({\rm N},F^p\widehat{M})\longrightarrow
            {\rm HH}^n({\rm N},\widehat{M})).\]
By Lemma~\ref{lemma:graded_ss_vanishing-simple-version},
${}^ME_{\infty}^{p,q,s}=0$ for $s\neq q$.
This implies that the exact sequence
\[ 0 \longrightarrow F^{p+1}{\rm HH}^n({\rm N},\widehat{M})\longrightarrow
        F^p{\rm HH}^n({\rm N},\widehat{M})\longrightarrow
        {}^ME_{\infty}^{p,n-p}({\rm N},\widehat{M})\longrightarrow 0\]
is canonically split.
%by Lemma~\ref{lemma:graded_ss_vanishing}.
We obtain a canonical isomorphism
\[ {\rm HH}^n({\rm N};\widehat{M})\cong
    \bigoplus_p {}^ME_{\infty}^{p,n-p}({\rm N},\widehat{M}) \]
and hence there is no extension problem.
\qed

\bigskip

In particular,
applying Theorem~\ref{theorem:collapse-E-2-term-no-extension}
to the case where
$\widehat{M}={\rm M}/{\rm N}$,
we obtain the following corollary.

\begin{corollary}\label{cor:ssM/N} 
The spectral sequence 
\[ {}^ME_1^{p,q}({\rm N},{\rm M}/{\rm N})
   \Longrightarrow 
   {\rm HH}^{p+q}({\rm N},{\rm M}/{\rm N})\]
of $\mathbb{Z}$-graded $R$-modules
collapses from the $E_2$-page.
There is an isomorphism
\[ {\rm HH}^{n}({\rm N}, {\rm M}/{\rm N}) = \bigoplus_{n,s}{\rm HH}^{n,s}({\rm N},{\rm M}/{\rm N})\cong
   \bigoplus_{n,s}{}^ME_{\infty}^{n-s,s}({\rm N},{\rm M}/{\rm N})\]
of bigraded $R$-modules.
\end{corollary}
%\qed

\if 
In this section, we calculate {${\rm HH}^{\ast}({\rm N}_m(R), {\rm M}_m(R)/{\rm N}_m(R))$. %for $m\ge 3$. 
Throughout this section, assume that $m \ge 3$. 
We set  $A={\rm N}_m(R)$ and $M={\rm M}_m(R)$. 
Let 
\begin{multline*}
{\rm M}_m(R) \supset L_1= \left\{
\left(
\begin{array}{ccccc}
\ast & \ast & \ast & \cdots & \ast \\
\vdots & \ddots & \ddots & \cdots & \vdots \\ 
\ast & \ast & \ast & \ddots & \ast \\
\ast & \ast & \ast & \ddots & \ast \\
0 & \ast & \ast & \cdots & \ast \\
\end{array}
\right)
\right\} 
\supset L_2= \left\{
\left(
\begin{array}{ccccc}
\ast & \ast & \ast & \cdots & \ast \\
\vdots & \ddots & \ddots & \cdots & \vdots \\ 
\ast & \ast & \ast & \ddots & \ast \\
0 & \ast & \ast & \ddots & \ast \\
0 & 0 & \ast & \cdots & \ast \\
\end{array}
\right)
\right\}  \\ 
\supset \cdots 
\supset L_{m-1} = {\rm B}_m(R)= \left\{
\left(
\begin{array}{ccccc}
\ast & \ast & \ast & \cdots & \ast \\
\vdots & \ddots & \ddots & \cdots & \vdots \\ 
0 & 0& \ast & \ddots & \ast \\
0 & 0 & 0 & \ddots & \ast \\
0 & 0 & 0 & \cdots & \ast \\
\end{array}
\right)
\right\}. 
\end{multline*} 
Set $F^{0} = M/A, \; F^{1}=L_1/A, \; F^{2}=L_2/A, \; \ldots, \; F^{m-1}={\rm B}_m(R)/A$ and $F^{m}=0$.   
We have a filtration 
\[
M/A = F^{0} \supset F^{1} \supset F^{2} \supset \cdots \supset F^{m-1} \supset F^{m} = 0  
%0 = F^{m} \subset F^{m-1} \subset \cdots \subset  F^{2} \subset F^{1} \subset F^{0} = M/A
\]
of $A$-bimodules over $R$. 
We denote by ${\rm Gr}^{p}(M/A)$ the $p$-th associated graded module $F^{p}/F^{p+1}$. 
By Proposition~\ref{prop:spectralsequence}, 
we obtain a spectral sequence 
\[
E_{1}^{p, q} = {\rm HH}^{p+q}(A, {\rm Gr}^{p}(M/A)) \Longrightarrow {\rm HH}^{p+q}(A, M/A)
\]
with 
\[
d_{r} : E_{r}^{p, q} \to E_{r}^{p+r, q-r+1} 
\] 
for $r \ge 1$. Note that $E_{1}^{p, q} = 0$ unless $0 \le p \le m-1$ and $p+q \ge 0$. 
%The spectral sequence collapses at $E_{m}$-stage. 

%Let us  consider the spectral sequence
%\[ E_1^{p,q}\Longrightarrow {\rm HH}^{p+q}(A,M/A),\]
%where $M={\rm M}_m(R)$.
The $A$-bimodules $A$ and  $M/A$ are finitely generated
free modules over $R$.
We have
\[ A \cong R\{ I_m \} \oplus R\{ E_{i,j}|\ i<j \} \mbox{ and } M/A\cong R\{E_{i,j}|\ i\ge j\}/R\{I_m\}.\]
By setting 
\[ |E_{i,j}|=j-i, \]
$A$ and $M/A$ are graded $A$-bimodules.  Then $A= \oplus_{a=0}^{m-1} A_{a}$ and $M/A=\oplus_{a=-(m-1)}^{0} (M/A)_{a}$, where 
$A_a$ and $(M/A)_{a}$ are the homogeneous parts of $A$ and $M/A$ of degree $a$, respectively. 
We have
\[ C^p(A,M/A)\cong {\rm Hom}_R(A^{\otimes p},M/A)\cong
     (A^*)^{\otimes p}\otimes_R (M/A). \]
In the same way as in the case $C^*(A,A)$,
we can introduce a grading on the cochain complex $C^*(A,M/A)$
and obtain a decomposition
\[ C^*(A,M/A)=\bigoplus_{s\in \mathbb{Z}}C^{*,s}(A,M/A),\]          
which is compatible with the filtration on $M/A$.
In the same way as in the case $C^*(A,A)$,
we obtain the following theorem.

\begin{theorem}
The spectral sequence
$E_1^{p,q}\Longrightarrow {\rm HH}^{p+q}(A,M/A)$
collapses at the $E_2$-page and
there is no extension problem.
\end{theorem}

\fi

\subsection{Calculation of ${\rm HH}^{\ast}({\rm N}_m(R), {\rm M}_m(R)/{\rm N}_m(R))$}\label{subsection:calHHM/N} 
In this subsection, we assume that $m \ge 3$. Let us calculate ${\rm HH}^{\ast}({\rm N}_m(R), {\rm M}_m(R)/{\rm N}_m(R))$. By Corollary~\ref{cor:ssM/N}, we only need to calculate $E_{2}^{p, q}({\rm N}, {\rm M}/{\rm N})$.  
Put $F^{p} = F^{p}({\rm M}/{\rm N})$ and $E_{r}^{p, q} = E_{r}^{p, q}({\rm N}, {\rm M}/{\rm N})$.  
Recall that 
\[
{\rm M}/{\rm N} = F^{-(m-1)} \supset F^{-(m-2)} \supset F^{-(m-3)} \supset \cdots \supset F^{0}={\rm B}/{\rm N} \supset F^{1} = 0,  
\]
where ${\rm B} = {\rm B}_m(R)  = \{ (a_{ij}) \in {\rm M}_m(R) \mid a_{ij} = 0 \mbox{ for } i>j \}$. 
It is easy to see that ${\rm Gr}^{p}({\rm M}/{\rm N}) = F^{p}/F^{p+1}$ is isomorphic to the direct sum of finitely many copies of $R$ as  an ${\rm N}$-bimodule over $R$. Hence we have 
\[ 
E_{1}^{p, q} \cong {\rm HH}^{p+q}({\rm N}, R)\otimes_{R}(F^{p}/F^{p+1}) \cong {\rm N}^{!}_{p+q}\otimes _{R}(F^{p}/F^{p+1}). 
\]  
Since $F^{p}/F^{p+1}$ is a free module over $R$, so is $E_{1}^{p, q}$. 
Note that ${\rm rank}_{R} (F^{-(m-1)}/F^{-(m-2)})= 1, \; {\rm rank}_{R} (F^{-(m-2)}/F^{-(m-3)})= 2, \; 
\ldots, \; {\rm rank}_{R} (F^{-1}/F^{0}) = m-1$ and ${\rm rank}_{R} (F^{0}/F^{1})= m-1$. Then we have 
\[\begin{array}{ccl}
{\rm rank}_{R} E_{1}^{-(m-1), q} & = & \varphi(q-m+1), \\ %\quad (q \ge 0) 
{\rm rank}_{R} E_{1}^{-(m-2), q} & = & 2 \varphi(q-m+2),  \\ %\quad (q \ge -1)
{\rm rank}_{R} E_{1}^{-(m-3), q} & = & 3 \varphi(q-m+3),  \\ %\quad (q \ge -2)
\cdots  & \cdots & \cdots \\ 
{\rm rank}_{R} E_{1}^{-1, q} & = & (m-1) \varphi(q-1),  \\ %\quad (q \ge -(m-2))  
{\rm rank}_{R} E_{1}^{0, q} & = & (m-1) \varphi(q),  \\ %\quad (q \ge -(m-1))
\end{array}
\] 
where $\varphi(n)= {\rm rank}_{R} {\rm N}^{!}_{n}$ in Definition~\ref{def:varphi}.  
The $R$-homomorphism $d_{1}^{p, q} : E_{1}^{p, q} \to E_{1}^{p+1, q}$ can be identified with the connecting homomorphism $\delta : {\rm HH}^{p+q}({\rm N}, F^{p}/F^{p+1}) \to {\rm HH}^{p+q+1}({\rm N}, F^{p+1}/F^{p+2})$ obtained by 
the short exact sequence $0 \to F^{p+1}/F^{p+2} \to F^{p}/F^{p+2} \to F^{p}/F^{p+1} \to 0$. 

We can write 
\[ 
{\rm Gr}^{p}({\rm M}/{\rm N}) = F^{p}/F^{p+1} = \left\{ 
\begin{array}{cl}
RE_{1-p, 1} \oplus RE_{2-p, 2} \oplus \cdots \oplus RE_{m, m+p} & (-(m-1) \le p \le -1), \\
(RE_{1,1}\oplus RE_{2,2} \oplus \cdots RE_{m,m})/RI_m & (p=0).  
\end{array}
\right. 
\] 
By Proposition~\ref{prop:quadraticmonomialalgebraKoszul} (\ref{item:monomialquad-3}),  
$\widehat{K}_{\ast}(A)$ gives us a graded free resolution of $A$ as $A$-bimodules over $R$ if $A={\rm N}$. 
Recall the differential $\widehat{d}_n : A\otimes_{R} (A^{!}_n)^{\ast}\otimes_{R}A \to A\otimes_{R} (A^{!}_{n-1})^{\ast}\otimes_{R} A$ is given by (\ref{eq:widehatKdifferential}). 
For $-(m-1) \le p \le -1$, 
\[
\begin{array}{cccl}
d_{1}^{p, q} : & {\rm N}^{!}_{p+q}\otimes {\rm Gr}^{p}({\rm M}/{\rm N}) & \longrightarrow & {\rm N}^{!}_{p+q+1}\otimes {\rm Gr}^{p+1}({\rm M}/{\rm N}) \\ 
 & f\otimes E_{i-p, i} & \longmapsto & y_{i-p-1} f \otimes E_{i-p-1, i} +(-1)^{p+q+1} fy_{i} \otimes E_{i-p, i+1} 
\end{array}
\]
for $1 \le i \le m+p$. 

Let us consider the complex 
\[
0 \longrightarrow E_{1}^{-(m-1), q} \stackrel{d_{1}^{-(m-1), q}}{\longrightarrow} E_{1}^{-(m-2), q} \stackrel{d_{1}^{-(m-2), q}}{\longrightarrow} E_{1}^{-(m-3), q} \stackrel{d_{1}^{-(m-3), q}}{\longrightarrow} \cdots \stackrel{d_{1}^{-2, q}}{\longrightarrow} E_{1}^{-1, q} \stackrel{d_{1}^{-1, q}}{\longrightarrow} E_{1}^{0, q} \longrightarrow 0. 
\]
We define another complex 
\[
0 \longrightarrow C_{1}^{-(m-1), q} \stackrel{\delta_{1}^{-(m-1), q}}{\longrightarrow} C_{1}^{-(m-2), q} \stackrel{\delta_{1}^{-(m-2), q}}{\longrightarrow} C_{1}^{-(m-3), q} \stackrel{\delta_{1}^{-(m-3), q}}{\longrightarrow} \cdots \stackrel{\delta_{1}^{-2, q}}{\longrightarrow} C_{1}^{-1, q} \stackrel{\delta_{1}^{-1, q}}{\longrightarrow} C_{1}^{0, q} \longrightarrow 0  
\]
by $C^{i, q} = E_{1}^{i, q}$ for $-(m-1) \le i \le -1$, $C^{0, q} = {\rm N}^{!}_{q}\otimes_{R} (RE_{1,1}\oplus RE_{2,2} \oplus \cdots \oplus RE_{m,m})$,  $\delta^{i, q} = d_{1}^{i, q}$ for $-(m-1) \le i \le -2$, and 
\[ 
\begin{array}{cccl}
\delta^{-1, q} : & C^{-1, q} = {\rm N}^{!}_{q-1}\otimes {\rm Gr}^{-1}({\rm M}/{\rm N}) & \longrightarrow & 
C^{0, q} = {\rm N}^{!}_{q}\otimes_{R} (\oplus_{k=1}^{m} RE_{k, k})  \\ 
 & f\otimes E_{i+1, i} & \longmapsto & y_{i} f \otimes E_{i, i} +(-1)^{q} fy_{i} \otimes E_{i+1, i+1} 
\end{array}
\] 
for $1\le i \le m-1$. Then we have a homomorphism of complexes 
\begin{eqnarray*}
\begin{array}{ccccccccccccc}
0 &\longrightarrow & C^{-(m-1), q} & \stackrel{\delta^{-(m-1), q}}{\longrightarrow} & C^{-(m-2), q}  & \stackrel{\delta^{-(m-2), q}}{\longrightarrow} & \cdots & \stackrel{\delta^{-2, q}}{\longrightarrow} & C^{-1, q} & \stackrel{\delta^{-1, q}}{\longrightarrow} & C^{0, q} & \longrightarrow & 0 \\ 
& & \downarrow \phi_{-(m-1)}  & & \downarrow \phi_{-(m-2)}  & & & & \downarrow \phi_{-1} & & \downarrow \phi_{0}  &  & \\ 
0 & \longrightarrow & E_{1}^{-(m-1), q} & \stackrel{d_{1}^{-(m-1), q}}{\longrightarrow} & E_{1}^{-(m-2), q} & \stackrel{d_{1}^{-(m-2), q}}{\longrightarrow} 
& \cdots & \stackrel{d_{1}^{-2, q}}{\longrightarrow} & E_{1}^{-1, q} & \stackrel{d_{1}^{-1, q}}{\longrightarrow} & E_{1}^{0, q} & \longrightarrow & 0, 
\end{array}
\end{eqnarray*} 
where $\phi_i = id_{E_{1}^{i, q}}$ for $-(m-1) \le i \le -1$ and $\phi_{0} : C^{0, q} = {\rm N}^{!}_{q}\otimes_{R} (\oplus_{k=1}^{m} RE_{k,k}) \to E_{1}^{0, q} = {\rm N}^{!}_{q}\otimes_{R} ((\oplus_{k=1}^{m} RE_{k, k})/RI_m)$ is the projection. 

\begin{remark}\rm 
We can give an interpretation of $C^{p, q}$ from a viewpoint of spectral sequences. 
%Let $M'={\rm M}_m(R)/{\rm J}_{m}(R)$. 
Let ${\rm J}=J({\rm N}_m(R))$. 
As in the case ${\rm M}/{\rm N}$, we put $F'^{p}=F^{p}({\rm M}/{\rm J})$. 
%$F'^{0}={\rm M}/{\rm J}, F^{1} = L_{1}/J,  F'^{2} =L_2/J, \ldots, F'^{m-1} = {\rm B}_{m}(R)/J$ and 
%$F'^{m}=0$. 
We have a filtration 
\[
{\rm M}/{\rm J} = F'^{-(m-1)} \supset F'^{-(m-2)} \supset \cdots \supset  F'^{1} \supset F'^{0} \supset F'^{-1} = 0   
\]
of $A$-bimodules over $R$. 
We denote by ${\rm Gr}^{p}({\rm M}/{\rm J})$ the $p$-th associated graded module $F'^{p}/F'^{p+1}$. 
By Proposition~\ref{prop:spectral-sequence-matrix-filtration}, we obtain a spectral sequence 
\[
E'^{p, q}_{1} = {\rm HH}^{p+q}({\rm N}, {\rm Gr}^{p}({\rm M}/{\rm J})) \Longrightarrow {\rm HH}^{p+q}({\rm N}, {\rm M}/{\rm J})
\]
with 
\[
d_{r} : E'^{p, q}_{r} \longrightarrow E'^{p+r, q-r+1}_{r} 
\] 
for $r \ge 1$. Then $\{ C^{p, q} \} \cong \{ E'^{p, q}_{1} \}$ as chain complexes. The homomorphism 
$\{ C^{p, q} \to E^{p, q}_{1} \}$ of chain complexes can be identified with the canonical map of spectral sequences $E'^{p, q}_{1} \to E^{p, q}_{1}$. 
\end{remark}

\bigskip 

Set $G^{p} = {\rm Gr}^p({\rm M}/{\rm N})$ for $-(m-1) \le p \le -1$ and $G^{0} = \oplus_{k=1}^{m} RE_{k, k}$. 
For $-(m-1) \le p \le -1$, we define the $R$-homomorphism $s^{p+1, q} : C^{p+1, q} \to C^{p, q}$ by 
\[
\begin{array}{ccl}
%s^{p+1, q} : & 
C^{p+1, q} = {\rm N}^{!}_{p+q+1}\otimes G^{p+1} & \stackrel{s^{p+1, q}}{\longrightarrow} & C^{p, q} = {\rm N}^{!}_{p+q}\otimes G^{p} \vspace*{1ex} \\ 
  f \otimes E_{i, i+p+1} & \longmapsto & 
\left\{ 
\begin{array}{ll}
 0 & \mbox{ if } f \notin y_{i} {\mathcal B}({\rm N}^{!}_{p+q}) \vspace*{1ex} \\ 
 f' \otimes E_{i+1, i+p+1}  & 
\begin{array}{l}
\mbox{ if }  f= y_{i}f'  \\ 
\mbox{ for } f' \in {\mathcal B}({\rm N}^{!}_{p+q})   
\end{array}  
\end{array} (-p \le i \le m-1), 
 \right. \vspace*{1ex} \\ 
  f\otimes E_{m, m+p+1} & \longmapsto &   0 \qquad (i=m)
\end{array}
\]
for $f \in {\mathcal B}({\rm N}^{!}_{p+q+1})$. % and $0 \le i \le p+1$.  

\begin{lemma}\label{lemma:homotopy} 
For $-(m-2) \le p \le -1$, $\delta^{p-1, q} \circ s^{p, q} + s^{p+1, q} \circ \delta^{p, q} = id_{C^{p, q}}$ and $s^{-(m-2), q} \circ \delta^{-(m-1), q} = id_{C^{-(m-1), q}}$: 
%\begin{eqnarray*}
%\begin{array}{ccccccccccccccc}
%0 & \to & C^{-(m-1), q} & \stackrel{\delta^{-(m-1), q}}{\to} & C^{-(m-2), q}  & \stackrel{\delta^{-(m-2), q}}{\to} & \cdots & 
% \stackrel{\delta^{-3, q}}{\to} & C^{-2, q} & \stackrel{\delta^{-2, q}}{\longrightarrow} & C^{-1, q} & \stackrel{\delta^{-1, q}}{\longrightarrow} & C^{0, q} & \to & 0 \\ 
%& & \downarrow_{id}  & \swarrow_{s^{-(m-2), q}} & \downarrow_{id}  & \swarrow_{s^{-(m-3), q}} & & \swarrow_{s^{-2, q}}  & \downarrow_{id}  & \swarrow_{s^{-1, q}} &\downarrow_{id} & \swarrow_{s^{0, q}} & \downarrow_{id}  &  & \\ 
%0 & \to & C^{-(m-1), q} & \stackrel{\delta^{-(m-1), q}}{\to} & C^{-(m-2), q}  & \stackrel{\delta^{-(m-2), q}}{\to} &  \cdots & 
%\stackrel{\delta^{-3, q}}{\to} & C^{-2, q} & \stackrel{\delta^{-2, q}}{\longrightarrow} & C^{-1, q} & \stackrel{\delta^{-1, q}}{\longrightarrow} & C^{0, q} & \to & 0. \\ 
%\end{array}
%\end{eqnarray*} 
\begin{eqnarray*}
\xymatrix{
	0 \ar[r]  & C^{-(m-1), q\quad} \ar[r]^{\hspace*{-2ex}\delta^{-(m-1),q}} \ar[d]_{id}&  C^{-(m-2), q \qquad} \ar[r]^{\qquad\delta^{-(m-2),q}}  \ar[d]_{id}\ar[ld]_{s^{-(m-2),q}} & \ar[ld]_{\hspace*{3ex}s^{-(m-3),q}} &\!\!\!\!\!\!\cdots\!\!\!\!\!\!&  \ar[r]^{\delta^{-2,q}} & C^{-1, q} \ar[ld]_{s^{-1,q}}\ar[r]^{\delta^{-1,q}}  \ar[d]_{id}& C^{0, q} \ar[r] \ar[d]_{id} \ar[ld]_{s^{0,q}}\ar[r] & 0 
	\\
	0 \ar[r]  & C^{-(m-1), q\quad} \ar[r]^{\hspace*{-1ex}\delta^{-(m-1),q}} &  C^{-(m-2), q\qquad } \ar[r]^{\qquad\delta^{-(m-2),q}} & &\!\!\!\!\!\!\cdots\!\!\!\!\!\!& \ar[r]^{\delta^{-2,q}} & C^{-1, q} \ar[r]^{\delta^{-1,q}} & C^{0, q} \ar[r]  & 0. 
	}
\end{eqnarray*} 
\end{lemma}
\proof
Let us prove that $s^{-(m-2), q} \circ \delta^{-(m-1), q} = id_{C^{-(m-1), q}}$. For $f \otimes E_{m, 1} \in C^{-(m-1), q} = {\rm N}^{!}_{q-(m-1)} \otimes_{R}  RE_{m, 1}$, 
\begin{eqnarray*}
s^{-(m-2), q} \circ \delta^{-(m-1), q}(f \otimes E_{m, 1}) & = & s^{-(m-2), q}(y_{m-1}f \otimes E_{m-1, 1} + (-1)^{q-m+2} fy_1 \otimes E_{m, 2}) \\ 
& = & f\otimes E_{m, 1}.  
\end{eqnarray*}
Hence, $s^{-(m-2), q} \circ \delta^{-(m-1), q} = id_{C^{-(m-1), q}}$. 

Let us show that $\delta^{p-1, q} \circ s^{p, q} + s^{p+1, q} \circ \delta^{p, q} = id_{C^{p, q}}$ for $-(m-2) \le p \le -1$. 
Note that $C^{p, q} = {\rm N}^{!}_{p+q} \otimes_{R} (RE_{1-p, 1} \oplus RE_{2-p, 2} \oplus \cdots \oplus RE_{m, m+p})$. 
It suffices to prove that $(\delta^{p-1, q} \circ s^{p, q} + s^{p+1, q} \circ \delta^{p, q})(f\otimes E_{i-p, i}) = f\otimes E_{i-p, i}$ for $-(m-2) \le p \le -1$, $1 \le i \le m+p$, and $f \in {\mathcal B}({\rm N}^{!}_{p+q})$. 
%For $f \in A^{!}_{p+q}\otimes RE_{m-p+i, i+1}$ with $0 \le i \le p$,  
Note that 
\begin{eqnarray*}
\delta^{p, q}(f\otimes E_{i-p, i}) & = & y_{i-p-1}f\otimes E_{i-p-1, i} + (-1)^{p+q+1} fy_{i} \otimes E_{i-p, i+1}.  
\end{eqnarray*}   
Assume that $1 \le i \le m+p-1$ and that $f \in {\mathcal B}({\rm N}^{!}_{p+q})$. Then  
\begin{eqnarray*}
s^{p+1, q} \circ \delta^{p, q} (f\otimes E_{i-p, i}) & = & 
\left\{
\begin{array}{cc}
f \otimes E_{i-p, i} & (f \notin y_{i-p} {\mathcal B}({\rm N}^{!}_{p+q-1})), \\ 
(-1)^{p+q+1} f' y_{i} \otimes E_{i-p+1, i+1} & (f = y_{i-p} f').  
\end{array}
\right. 
\end{eqnarray*}
Since 
\begin{eqnarray*}
s^{p, q}  (f\otimes E_{i-p, i}) & = & 
\left\{
\begin{array}{cc}
0 & (f \notin y_{i-p} {\mathcal B}({\rm N}^{!}_{p+q-1})), \\ 
f'  \otimes E_{i-p+1, i} & (f = y_{i-p} f'),   
\end{array}
\right. 
\end{eqnarray*}
we have 
\begin{eqnarray*}
& & \delta^{p-1, q} \circ s^{p, q}  (f\otimes E_{i-p, i})  \\ 
& = & 
\left\{
\begin{array}{cc}
0 & (f \notin y_{i-p} {\mathcal B}({\rm N}^{!}_{p+q-1})), \\ 
f  \otimes E_{i-p, i}+(-1)^{p+q} f' y_{i} E_{i-p+1, i+1} & (f = y_{i-p} f').    
\end{array}
\right. 
\end{eqnarray*}
Thereby, we obtain
\[
(\delta^{p-1, q} \circ s^{p, q} + s^{p+1, q} \circ \delta^{p, q})(f \otimes E_{i-p, i}) = f \otimes E_{i-p, i}.   
\]

Assume that $i=m+p$ and that $f \in {\mathcal B}({\rm N}^{!}_{p+q})$. Then 
$s^{p+1, q} \circ \delta^{p, q} (f\otimes E_{m, m+p}) =  
f \otimes E_{m, m+p}$.  
Since $s^{p, q}(f\otimes E_{m, m+p}) = 0$, $\delta^{p-1, q} \circ s^{p, q}(f\otimes E_{m, m+p}) = 0$. 
Hence $(\delta^{p-1, q} \circ s^{p, q} + s^{p+1, q} \circ \delta^{p, q})(f \otimes E_{m, m+p}) = f \otimes E_{m, m+p}$. 
This completes the proof. 
\qed 

\begin{lemma}\label{lemma:scalardelta0}
Let $K^{0, q} = {\rm N}^{!}_{q} \otimes_{R} RI_m \subseteq C^{0, q} = {\rm N}^{!}_{q}\otimes_{R} (\sum_{k=1}^m RE_{k, k})$.  
Then $\delta^{-1, q}(C^{-1, q}) \cap K^{0, q} = 0$. 
\end{lemma} 
\proof
Let $x = f_2\otimes E_{2, 1}+f_3\otimes E_{3, 2} + \cdots+ f_i \otimes E_{i, i-1} + \cdots + f_{m} \otimes E_{m, m-1} \in C^{-1, q}$, where $f_i \in {\rm N}^{!}_{q-1}$ ($2 \le i \le m$). Note that 
\begin{multline*}
\delta^{-1, q}(x)  =  y_1f_2\otimes E_{1, 1} + ( (-1)^q f_2 y_1+y_2f_3) E_{2, 2} + ( (-1)^q f_3 y_2+ y_3f_4) E_{3, 3} + 
\cdots  \\ 
 + ( (-1)^q f_i y_{i-1} + y_{i} f_{i+1}) E_{i, i} + \cdots + ( (-1)^q f_{m-1}y_{m-2} + y_{m-1} f_{m})E_{m-1, m-1} \\ + (-1)^q f_{m} y_{m-1} \otimes E_{m, m}.   
\end{multline*} 
Suppose that $\delta^{-1, q}(x) \in K^{0, q}$. Then 
\begin{align}
y_1f_2 & =  (-1)^q f_2y_1+y_2f_3, \label{eq:Kcapdelta-1} \\
(-1)^q f_2y_1+y_2f_3 & =  (-1)^q f_3y_2+y_3f_4, \label{eq:Kcapdelta-2} \\ 
 & \cdots  \nonumber \\  
(-1)^q f_{m-1}y_{m-2}+y_{m-1}f_m & =  (-1)^q f_my_{m-1}.  \label{eq:Kcapdelta-4}\tag{5.$m$}
\end{align}
We can write $f_2 = f'_2 + y_2 f''_2$ such that  $f'_2 \in y_1 {\rm N}^{!}_{q-2} \oplus y_3 {\rm N}^{!}_{q-2} \oplus \cdots \oplus y_{m-1} {\rm N}^{!}_{q-2}$ and $f''_{2} \in {\rm N}^{!}_{q-2}$. 
By (\ref{eq:Kcapdelta-1}), we obtain  
\begin{eqnarray*}
%x_1 f_2 -(-1)^a f_2x_1 & = & x_2f_3 \\ 
y_1 f'_2 & = & (-1)^q f'_2 y_1 + (-1)^q y_2 f''_2 y_1 +  y_2f_3,   
\end{eqnarray*}
and hence 
\begin{eqnarray}
y_1 f'_2 -(-1)^q f'_2 y_1 & = & y_2 ((-1)^q  f''_2 y_1+f_3).  \label{eq:Kcapdelta-f2} 
\end{eqnarray}
Since the right hand side of (\ref{eq:Kcapdelta-f2}) has the leading term $y_2$, we have $y_1f'_2 - (-1)^q f'_2y_1=0$.  
Hence $y_1f'_2 = (-1)^q f'_2y_1$. We see that $f'_2 = c y_1^{q-1}$ for some $c \in R$.  

Similarly, we can write $f_{m} = f'_{m} + f''_{m} y_{m-2}$ such that  $f'_{m} \in {\rm N}^{!}_{q-2}y_1 \oplus \cdots \oplus  {\rm N}^{!}_{q-2}y_{m-3} \oplus {\rm N}^{!}_{q-2}y_{m-1}$ and $f''_{m} \in {\rm N}^{!}_{q-2}$. 
By (\ref{eq:Kcapdelta-4}), we obtain  
\begin{eqnarray*}
(-1)^q f_{m-1}y_{m-2} + y_{m-1} f'_m +y_{m-1}f''_{m} y_{m-2}  & = & (-1)^q f'_m y_{m-1}, 
\end{eqnarray*}
and hence 
\begin{eqnarray}
((-1)^q f_{m-1}  +y_{m-1}f''_{m}) y_{m-2}  & = & (-1)^q f'_m y_{m-1} - y_{m-1} f'_m.  \label{eq:Kcapdelta-fm} 
\end{eqnarray}
Since the left hand side of (\ref{eq:Kcapdelta-fm}) has the last term $y_{m-2}$, we have $(-1)^q f'_m y_{m-1} - y_{m-1} f'_m = 0$.   
Hence $y_{m-1} f'_m = (-1)^q f'_m y_{m-1}$. We see that $f'_m = d y_{m-1}^{q-1}$ for some $d \in R$.  

By using (\ref{eq:Kcapdelta-1}), (\ref{eq:Kcapdelta-2}), $\ldots$, and (\ref{eq:Kcapdelta-4}), $y_1 f_2 = (-1)^q f_m y_{m-1}$. Then we obtain 
\begin{eqnarray*}
y_1 (f'_2+ y_{2} f''_2) & = &  (-1)^q (f'_{m}+f''_{m}y_{m-2})y_{m-1}, \\ 
y_1 f'_2 & = & (-1)^q f'_{m} y_{m-1}. 
\end{eqnarray*} 
Since $f'_2= c y_{1}^{q-1}$ and $f'_{m} = d y_{m-1}^{q-1}$, $cy_{1}^{q} =(-1)^q dy_{m-1}^{q}$, which implies that  $c=d=0$. 
Hence $y_1 f_2= y_1(f'_2+ y_2f''_2) = y_1 y_2f''_2 =0$.  Therefore $\delta^{-1, q}(x) = 0$ by  (\ref{eq:Kcapdelta-1}), (\ref{eq:Kcapdelta-2}), $\ldots$, and (\ref{eq:Kcapdelta-4}).  
This completes the proof. 
\qed 

\begin{proposition}\label{prop:E2pq=0unlessq=0}
For $-(m-1) \le p \le -1$, we have $E_{2}^{p, q}=0$ for all $q$. In particular, $E_{2}^{p, q}=0$ unless $p=0$. 
\end{proposition}

\proof
For $-(m-1) \le p \le -2$, $E_{2}^{p, q} \cong H^{p}(C^{\ast, q}) = 0$ by Lemma~\ref{lemma:homotopy}.  
For $p=-1$, let us consider the commutative diagram with columns exact:  
\[
\begin{array}{ccccccccc}
& & & &   & & 0\phantom{\phi_{1}} & & \\ 
& & & &   & & \downarrow\phantom{\phi_{0}} & & \\ 
& & 0\phantom{\phi_{-2}} & & 0\phantom{\phi_{-1}}  & & K^{0, q} & & \\ 
& & \downarrow\phantom{\phi_{-2}} & & \downarrow\phantom{\phi_{-1}} & & \downarrow\phantom{\phi_{0}} & & \\ 
\cdots &  \stackrel{\delta^{-3, q}}{\longrightarrow} & C^{-2, q} & \stackrel{\delta^{-2, q}}{\longrightarrow} & C^{-1, q} & \stackrel{\delta^{-1, q}}{\longrightarrow} & C^{0, q} & \longrightarrow & 0 \\ 
& & \downarrow\phi_{-2} & & \downarrow\phi_{-1} & & \downarrow\phi_{0} & & \\ 
\cdots &  \stackrel{d_{1}^{-3, q}}{\longrightarrow} & E_{1}^{-2, q} & \stackrel{d_{1}^{-2, q}}{\longrightarrow} & E_{1}^{-1, q} & \stackrel{d_{1}^{-1, q}}{\longrightarrow} & E_{1}^{0, q} & \longrightarrow & 0 \\
& & \downarrow\phantom{\phi_{-2}} & & \downarrow\phantom{\phi_{-1}}  & & \downarrow\phantom{\phi_{0}} & & \\ 
& & 0\phantom{\phi_{-2}} & & 0\phantom{\phi_{-1}}  & &  0 . \phantom{\phi_{0}} & &  \\ 
\end{array}
\]
It is easy to see that $E_{2}^{-1, q} \cong H^{-1}(C^{\ast, q}) = 0$ by Lemmas~\ref{lemma:homotopy} and \ref{lemma:scalardelta0}. 
The last statement follows from that $E_{1}^{p, q} = 0$ unless $-(m-1) \le p \le 0$.    
\qed 

%\bigskip 
%
%By the discussions above,  we have $E_{2}^{p, q} = 0$ unless $p=0$. 
%%Hence the spectral sequence collapses at $E_2$-stage. 

\bigskip 

We introduce the following claim. This is true for the case $R$ is a field. 

\begin{claim}\label{claim:free} 
The $R$-module $E_{2}^{0, q} \cong E_{1}^{0, q}/ {\rm Im} \; d_{1}^{-1, q}$ is free. 
\end{claim}

\bigskip

Under the hypothesis that Claim~\ref{claim:free} is true, let us calculate ${\rm rank}_{R} E_{2}^{0, q}$. 
By taking the Euler characteristic of the cochain complex 
\begin{eqnarray*}
0 \longrightarrow E_{1}^{-(m-1), q} \stackrel{d_{1}^{-(m-1), q}}{\longrightarrow} E_{1}^{-(m-2), q} \stackrel{d_{1}^{-(m-2), q}}{\longrightarrow} 
%E_{1}^{-(m-3), q} \stackrel{d_{1}^{-(m-3), q}}{\to} 
\cdots \stackrel{d_{1}^{-2, q}}{\longrightarrow} E_{1}^{-1, q} \stackrel{d_{1}^{-1, q}}{\longrightarrow} E_{1}^{0, q} \longrightarrow 0,  
\end{eqnarray*}
we obtain  
\begin{eqnarray*}
\sum_{k=-(m-1)}^{0} (-1)^k  {\rm rank}_{R} E_{2}^{k, q} & = & \sum_{k=-(m-1)}^{0} (-1)^k  {\rm rank}_{R} E_{1}^{k, q} \\
& = & (m-1)\varphi(q) + \sum_{k=-(m-1)}^{-1} (-1)^k (m+k) \varphi(q+k).  
\end{eqnarray*} 
Since $E_{2}^{k, q} = 0$ for $-(m-1) \le k \le -1$ by Proposition~\ref{prop:E2pq=0unlessq=0}, 
\begin{eqnarray*}
{\rm rank}_{R} E_{2}^{0, q} & = & (m-1)\varphi(q) + \sum_{k=-(m-1)}^{-1} (-1)^k (m+k) \varphi(q+k).  
\end{eqnarray*} 
Hence 
\begin{eqnarray}\label{eq:rankE2m-1q}
{\rm rank}_{R} E_{2}^{0, q}  & = &  (m-1)\varphi(q) + \sum_{k=1}^{m-1} (-1)^{m+k} k \varphi(q-m+k).  \label{eq:rankE2}
\end{eqnarray}

Now, let us prove Claim~\ref{claim:free}. 
To emphasize $R$, we denote by $d_{1}^{-1, q}(R) : E_{1}^{-1, q}(R) \to E_{1}^{0, q}(R)$ the differential $d_{1}^{-1, q} : E_{1}^{-1, q} \to E_{1}^{0, q}$.  We can ragard $d_{1}^{-1, q}(S)$ as $d_{1}^{-1, q}(R)\otimes_{R} S$ for any ring homomorphism $R \to S$. 

\begin{lemma}\label{lemma:E2free} 
The $R$-module $E_{2}^{0, q} \cong E_{1}^{0, q}/ {\rm Im} \; d_{1}^{-1, q}$ is a free module over $R$. The rank of $E_{2}^{0, q}$ over $R$ is given by (\ref{eq:rankE2m-1q}).  
\end{lemma}

\proof
When $R$ is a field, the statement is true. 
For $R={\Bbb Z}$, let us consider the exact sequence  
\begin{eqnarray}\label{eq:exactseqZ} 
E_{1}^{-1, q}({\Bbb Z}) \stackrel{d_{1}^{-1, q}({\Bbb Z})}{\longrightarrow} E_{1}^{0, q}({\Bbb Z}) \longrightarrow E_{2}^{0, q}({\Bbb Z}) \longrightarrow 0. 
\end{eqnarray}
Note that $E_{2}^{0, q}({\Bbb Z})$ is finitely generated over ${\Bbb Z}$. 
Suppose that $E_{2}^{0, q}({\Bbb Z})$ has a torsion element $x \neq 0$ such that $px = 0$ for a prime number $p$.  
By tensoring (\ref{eq:exactseqZ}) by ${\Bbb F}_{p}$,  we have an exact sequence 
\begin{eqnarray*}
E_{1}^{-1, q}({\Bbb F}_p) \stackrel{d_{1}^{-1, q}({\Bbb F}_p)}{\longrightarrow} E_{1}^{0, q}({\Bbb F}_p) \longrightarrow E_{2}^{0, q}({\Bbb Z})\otimes_{{\Bbb Z}} {\Bbb F}_{p} \longrightarrow 0,  
\end{eqnarray*}
which implies that $E_{2}^{0, q}({\Bbb F}_p) \cong E_{2}^{0, q}({\Bbb Z})\otimes_{{\Bbb Z}} {\Bbb F}_{p}$. 
By tensoring (\ref{eq:exactseqZ}) by ${\Bbb Q}$, we also have $E_{2}^{0, q}({\Bbb Q}) \cong E_{2}^{0, q}({\Bbb Z})\otimes_{{\Bbb Z}} {\Bbb Q}$. By the fundamental theorem of finitely generated abelian groups, 
$\dim_{{\Bbb Q}}  E_{2}^{0, q}({\Bbb Z})\otimes_{{\Bbb Z}} {\Bbb Q} < \dim_{{\Bbb F}_p} E_{2}^{0, q}({\Bbb Z})\otimes_{{\Bbb Z}} {\Bbb F}_{p}$, which contradicts to the fact that both $E_{2}^{0, q}({\Bbb Q})$ and $E_{2}^{0, q}({\Bbb F}_p)$ have the same rank (\ref{eq:rankE2m-1q}).  Hence, $E_{2}^{0, q}({\Bbb Z})$ has no torsion element.  By using the fundamental theorem of finitely generated abelian groups again, we see that $E_{2}^{0, q}({\Bbb Z})$ is a free module of rank (\ref{eq:rankE2m-1q}) over ${\Bbb Z}$. 

Let us consider the case that $R$ is an arbitrary commutative ring.  By tensoring (\ref{eq:exactseqZ}) by $R$, 
we have $E_{2}^{0, q}(R) \cong E_{2}^{0, q}({\Bbb Z})\otimes_{{\Bbb Z}} R$. 
Since $E_{2}^{0, q}({\Bbb Z})$ is a free module of rank (\ref{eq:rankE2m-1q}) over ${\Bbb Z}$, 
$E_{2}^{0, q}(R)$ is also a free module of rank (\ref{eq:rankE2m-1q}) over $R$. This completes the proof. 
%For the proof, let us consider the following exact sequence: 
%\begin{eqnarray}
%E_{1}^{m-2, q}(R) \stackrel{d_{1}^{m-2, q}(R)}{\longrightarrow} E_{1}^{m-1, q}(R) \to E_{2}^{m-1, q}(R) \to 0. 
%\end{eqnarray}
%Let us show the statement when $R={\Bbb Z}$. 
\qed

\begin{theorem}\label{th:HHMmNm} 
Let $m \ge 3$. 
The cohomology group ${\rm HH}^{n}({\rm N}_m(R), {\rm M}_m(R)/{\rm N}_m(R))$ is a free module over $R$ for $n \ge 0$. The rank of ${\rm HH}^{n}({\rm N}_m(R), {\rm M}_m(R)/{\rm N}_m(R))$ is given by  
\[
{\rm rank}_{R} {\rm HH}^{n}({\rm N}_m(R), {\rm M}_m(R)/{\rm N}_m(R)) =  (m-1)\varphi(n) + \sum_{k=1}^{m-1} (-1)^{m+k} k \varphi(n-m+k).
\]
\end{theorem}

\proof 
The spectral sequence collapses from the $E_2$-page and  there is no extension problem by Corollary~\ref{cor:ssM/N}.  By Lemma~\ref{lemma:E2free}, $E_{2}^{0, q}$ is a free module over $R$.  The statement follows from 
that 
\[{\rm HH}^{n}({\rm N}_m(R), {\rm M}_m(R)/{\rm N}_m(R)) \cong E_{2}^{0, n}
\] 
and (\ref{eq:rankE2}). 
\qed 

\bigskip 

%Recall that ${\rm N}_m(R) \cong R\langle x_1, x_2, \ldots, x_{m-1} \rangle/\langle x_i x_j \mid j \neq i+1 \rangle$. 
%Setting $\deg x_i =1$, we can regard ${\rm N}_m(R)$ as a graded $R$-algebra.  
%Let us denote by ${\rm N}_m(R)_{n}$ the homogeneous part of degree $n$ of ${\rm N}_m(R)$.  
%Let us define 
Recall 
\begin{eqnarray*}
f(t) & = & \sum_{n=0}^{\infty} \left({\rm rank}_{R} {\rm N}_m(R)_{n} \right) t^n = 1+ \sum_{k=1}^{m-1} kt^{m-k},   \\
f^{!}(t) & = & \sum_{n=0}^{\infty} \left({\rm rank}_{R} {\rm N}_m(R)^{!}_{n}  \right) t^n = \frac{1}{f(-t)}=\frac{1}{\displaystyle 1+ \sum_{k=1}^{m-1} (-1)^{m-k}kt^{m-k}}  
\end{eqnarray*} 
in (\ref{eq:f}), (\ref{eq:explicitf}), (\ref{eq:f^{!}}), and Proposition~\ref{prop:ff^{!}}. Let us define 
\begin{eqnarray*}
h(t) & = & \sum_{n=0}^{\infty} \left({\rm rank}_{R} {\rm HH}^{n}({\rm N}_m(R), {\rm M}_m(R)/{\rm N}_m(R))\right) t^n. \\
\end{eqnarray*} 
%Note that 
%\[
%h(t) = 1+ \sum_{i=1}^{m-1}(m-i)t^i = 1 + \sum_{i=1}^{m-1} i t^{m-i}.  
%\]
%By \cite[{Theorem~3.5.1}]{Loday-Vallette}, we obtain 
%\[
%f^{{\rm N}_{m}^{!}(R)}(t) h(-t) = 1. 
%\]
%Hence, 
%\[
%f^{{\rm N}_{m}^{!}(R)}(t) = \frac{1}{h(-t)} = \frac{1}{\displaystyle 1+\sum_{i=1}^{m-1}(-1)^{i} (m-i) t^{i}}. 
%\]

\begin{theorem}\label{th:generating} 
Let $m \ge 3$. 
The generating function $h(t)$ is given by 
\[
h(t) = 1+(m-2) f^{!}(t). 
\]
\end{theorem}

\proof
By Theorem~\ref{th:HHMmNm}, 
\begin{eqnarray*}
h(t) & = & \left\{ (m-1) + \sum_{k=1}^{m-1} (-1)^{m+k} k t^{m-k} \right\} f^{!}(t) \\  
 & = & \frac{\displaystyle (m-1) + \sum_{k=1}^{m-1} (-1)^{m-k} k t^{m-k}}{\displaystyle 1+ \sum_{k=1}^{m-1} (-1)^{m-k}kt^{m-k}} \\  
%& = & \frac{\displaystyle (m-1) + \sum_{i=0}^{m-2} (-1)^{m+i-1} (i+1) t^{m-i-1}}{h(-t)} \\
%& = & \frac{\displaystyle (m-1) + \sum_{i=1}^{m-1} (-1)^{m-i} i t^{m-i}}{\displaystyle 1 + \sum_{i=1}^{m-1} (-1)^{m-i} i t^{m-i}} \\
& = & 1+ \frac{m-2}{\displaystyle 1 + \sum_{k=1}^{m-1} (-1)^{m-k} k t^{m-k}} \\
& = & 1+ (m-2) f^{!}(t). 
\end{eqnarray*} 
This completes the proof. 
\qed 

\begin{corollary}\label{cor:HHMmNm} 
Let $m \ge 3$. The rank of ${\rm rank}_{R} {\rm HH}^{n}({\rm N}_m(R), {\rm M}_m(R)/{\rm N}_m(R))$ is given by
\[
{\rm rank}_{R} {\rm HH}^{n}({\rm N}_m(R), {\rm M}_m(R)/{\rm N}_m(R)) = \left\{ 
\begin{array}{cc}
m-1 & (n=0), \\
(m-2)\varphi(n) & (n>0).  
\end{array}
\right. 
\]
\end{corollary}

\proof
Note that the constant term of $f^{!}(t)$ is $1$. 
The statement follows from Theorem~\ref{th:generating}.  
\qed 

\begin{remark}\rm 
If $m=3$, then $\varphi(n)=n+1$ for $n\ge 0$. We can easily check that the result of 
Corollary~\ref{cor:HHMmNm} is compatible with \cite[Theorem~5.4]{Nakamoto-Torii:Hochschild}:  
\[
{\rm rank}_{R} {\rm HH}^{n}({\rm N}_3(R), {\rm M}_3(R)/{\rm N}_3(R)) = \left\{ 
\begin{array}{lc}
2 & (n=0), \\
n+1 & (n>0).  
\end{array}
\right. 
\]
\end{remark}

\bigskip 

Recall that ${\rm HH}^{n, s}({\rm N}_m(R), {\rm M}_m(R)/{\rm N}_m(R)) = H^{n}(C^{\ast, s}({\rm N}_m(R), {\rm M}_m(R)/{\rm N}_m(R)))$ (cf. (\ref{eq:defofHn,s})). 
By the result above, we have: 

\begin{theorem}
Let $m \ge 3$. For each $n \ge 0$ and $s \in {\Bbb Z}$, 
${\rm HH}^{n, s}({\rm N}_m(R), {\rm M}_m(R)/{\rm N}_m(R))$ is a free $R$-module. 
The rank is given by 
\[
{\rm rank}_{R} {\rm HH}^{n, s}({\rm N}_m(R), {\rm M}_m(R)/{\rm N}_m(R)) = \left\{ 
\begin{array}{ll}
0 & (n\neq s), \\
m-1 & (n=s=0), \\
(m-2)\varphi(n) & (n=s>0).   
\end{array}
\right. 
\]
\end{theorem}

\proof
As in Corollary~\ref{cor:ssM/N}, we have
\begin{eqnarray*}
{\rm HH}^{n, s}({\rm N}_m(R), {\rm M}_m(R)/{\rm N}_m(R)) & \cong & 
E_{\infty}^{n-s, s}({\rm N}_m(R), {\rm M}_m(R)/{\rm N}_m(R))  \\
& \cong & E_{2}^{n-s, s}({\rm N}_m(R), {\rm M}_m(R)/{\rm N}_m(R)). 
\end{eqnarray*} 
By the discussion above, we can verify the statement. 
\qed 

\bigskip 

As in  the proofs of Claim~\ref{claim:free}, Lemma~\ref{lemma:E2free}, and Theorem~\ref{th:HHMmNm}, 
we can show the following: 

\begin{proposition}
Let $m \ge 3$. 
The cohomology group ${\rm HH}^{n}({\rm N}_m(R), {\rm M}_m(R)/J({\rm N}_m(R)))$ is a free module over $R$ for $n \ge 0$. The rank of ${\rm HH}^{n}({\rm N}_m(R), {\rm M}_m(R)/J({\rm N}_m(R)))$ is given by  
\[
{\rm rank}_{R} {\rm HH}^{n}({\rm N}_m(R), {\rm M}_m(R)/J({\rm N}_m(R))) = \left\{ 
\begin{array}{cc}
m & (n=0), \\
(m-1)\varphi(n) & (n>0).  
\end{array}
\right. 
\]
For each $n \ge 0$ and $s \in {\Bbb Z}$, 
${\rm HH}^{n, s}({\rm N}_m(R), {\rm M}_m(R)/J({\rm N}_m(R)))$ is a free $R$-module. 
The rank is given by 
\[
{\rm rank}_{R} {\rm HH}^{n, s}({\rm N}_m(R), {\rm M}_m(R)/J({\rm N}_m(R))) = \left\{ 
\begin{array}{ll}
0 & (n\neq s), \\
m & (n=s=0), \\
(m-1)\varphi(n) & (n=s>0).   
\end{array}
\right. 
\]
\end{proposition}

\subsection{The Zariski tangent space of the moduli of subalgebras of ${\rm M}_m$ at ${\rm N}_m$} 
In the previous subsection, we have calculated ${\rm rank}_{R} {\rm HH}^{n}({\rm N}_m(R), {\rm M}_m(R)/{\rm N}_m(R))$.  
In this subsection, we calculate the dimension of the tangent space of the moduli of subalgebras of 
${\rm M}_m$ over ${\Bbb Z}$ at ${\rm N}_m$ for $m\ge 3$ 
by using ${\rm rank}_{R} {\rm HH}^{1}({\rm N}_m(R), {\rm M}_m(R)/{\rm N}_m(R))$. 

\begin{proposition}\label{prop:normalizerNm} 
Let $R$ be a commutative ring. 
Set 
\[
N({\rm N}_m(R)) = \{ A \in {\rm M}_m(R) \mid [A, B]:=AB-BA \in {\rm N}_m(R)  \mbox{ for any } B \in {\rm N}_m(R) \}. 
\]
For $m \ge 3$, $N({\rm N}_m(R)) = {\rm B}_m(R)$. 
\end{proposition}

\proof 
It is easy to see that $N({\rm N}_m(R)) \supseteq {\rm B}_m(R)$. 
Let us show that $N({\rm N}_m(R)) \subseteq {\rm B}_m(R)$. 
Recall the ${\Bbb Z}$-grading on ${\rm M} = {\rm M}_m(R)$ in \S\ref{subsection:sssubquotient}: 
\[ {\rm M}=\bigoplus_{r\in\mathbb{Z}}{\rm M}_r,  \mbox{ where } {\rm M}_r=\bigoplus_{j-i=r}R\{E_{i,j}\}.\]
Set $N({\rm N}_m(R))_r = N({\rm N}_m(R))\cap {\rm M}_r$. Since 
$N({\rm N}_m(R)) = \oplus_{r\in {\Bbb Z}} N({\rm N}_m(R))_r$ and ${\rm B}_m(R) = \oplus_{r \ge 0} {\rm M}_r$,  
it suffices to prove that $N({\rm N}_m(R))_r = 0$ for $-(m-1) \le r \le -1$. 
Suppose that there exists 
$x = a_{1}E_{1-r, 1} + a_{2} E_{2-r, 2} + \cdots + a_{m+r} E_{m, m+r} \in N({\rm N}_m(R))_r$ with 
$a_1=\cdots = a_{i-1} = 0$ and $a_i \neq 0$. 
If $-(m-1) \le r \le -2$, then 
\[
[x, E_{i, i+1} ] = a_i E_{i-r, i+1} \in {\rm N}_m(R).    
\]
This implies that $a_i = 0$, which is a contradiction. 
If $r=-1$, then 
\[
[x, E_{i, i+1} ] = a_i E_{i+1, i+1} -  a_{i} E_{i, i}  \in {\rm N}_m(R).      
\]
Since $m\ge 3$, we see that $a_i=0$, which is a contradiction. 
Hence,  $N({\rm N}_m(R)) = {\rm B}_m(R)$. 
\qed 

\bigskip 

Set $\displaystyle d = {\rm rank}_{R} {\rm N}_m(R) = \frac{m^2-m+2}{2}$. 
Recall the moduli of molds ${\rm Mold}_{m, d}$, in other words, 
the moduli of rank $d$ subalgebras of the full matrix ring ${\rm M}_m$ in \cite[\S3.1]{Nakamoto-Torii:Hochschild}. 
We can regard ${\rm N}_m$ as a point of ${\rm Mold}_{m, d}$. 
Let us consider the Zariski tangent space $T_{{\rm Mold}_{m, d}/{\Bbb Z}, {\rm N}_m}$ of  
${\rm Mold}_{m, d}$ over ${\Bbb Z}$ at ${\rm N}_m$ (for details, see  \cite[Definition~3.10]{Nakamoto-Torii:Hochschild}). 

\begin{theorem}\label{th:dimtangentspace} 
The dimension of the Zariski tangent space $T_{{\rm Mold}_{m, d}/{\Bbb Z}, {\rm N}_m}$ of  
${\rm Mold}_{m, d}$ over ${\Bbb Z}$ at ${\rm N}_m$ is 
\[
\dim T_{{\rm Mold}_{m, d}/{\Bbb Z}, {\rm N}_m} = \frac{3m^2 - 7m+4}{2} 
\]
for $m \ge 3$. 
\end{theorem}

\proof 
Let $m \ge 3$. For any field $k$,  
\[ 
\dim_{k} {\rm HH}^{1}({\rm N}_m(k), {\rm M}_m(k)/{\rm N}_m(k)) = (m-2)\varphi(1)=(m-2)(m-1) 
\]
by Corollaries~\ref{cor:varphi} and \ref{cor:HHMmNm}.   We also see that 
\[
\dim_{k} N({\rm N}_m(k)) = \dim_{k} {\rm B}_m(k) = \frac{m(m+1)}{2} 
\]
by Proposition~\ref{prop:normalizerNm}. 
Using \cite[Corollary~3.14]{Nakamoto-Torii:Hochschild}, we obtain 
\begin{eqnarray*} 
\dim T_{{\rm Mold}_{m, d}/{\Bbb Z}, {\rm N}_m} & = & 
 \dim_{k} {\rm HH}^{1}({\rm N}_m(k(x)), {\rm M}_m(k(x))/{\rm N}_m(k(x))) + m^2 - \dim_{k} N({\rm N}_m(k(x)))  \\
 & = & \frac{3m^2 - 7m+4}{2},  
\end{eqnarray*} 
where $k(x)$ is the residue field of $x={\rm N}_m$.  
\qed 

\begin{remark}\rm
By Theorem~\ref{th:dimtangentspace}, $\dim T_{{\rm Mold}_{3, 4}/{\Bbb Z}, {\rm N}_3} = 5$ for $m=3$. 
This result coincides with $\dim T_{{\rm Mold}_{3, 4}/{\Bbb Z}, {\rm N}_3}$ in \cite[Table~2]{Nakamoto-Torii:Hochschild}. 
\end{remark}

\begin{remark}\rm
In the case $m=2$, ${\rm N}_2(R)$ coincides with  ${\rm J}_2(R)$ defined in \cite[Definition~4.16]{Nakamoto-Torii:Hochschild}.  For any field $k$, we have obtained 
\[
\dim_{k} {\rm HH}^{1}({\rm N}_2(k), {\rm M}_2(k)/{\rm N}_2(k)) = \left\{ 
\begin{array}{cc}
1 &  ({\rm ch}(k) \neq 2), \\
2 &  ({\rm ch}(k) = 2) 
\end{array}
\right.
\]
by \cite[Corollary~4.20]{Nakamoto-Torii:Hochschild}. 
We also see that 
\[
N({\rm N}_2(k)) = \left\{ 
\begin{array}{cc}
{\rm B}_2(k) &  ({\rm ch}(k) \neq 2), \\
{\rm M}_2(k) &  ({\rm ch}(k) = 2) 
\end{array}
\right.
\]
by \cite[Proposition~4.21]{Nakamoto-Torii:Hochschild}. 
Using \cite[Corollary~3.14]{Nakamoto-Torii:Hochschild} or \cite[Example~4.22]{Nakamoto-Torii:Hochschild}, we obtain 
\[
\dim T_{{\rm Mold}_{2, 3}/{\Bbb Z}, {\rm N}_2} = 2,  
\]
while ${\rm Mold}_{2, 3} = {\Bbb P}^2_{{\Bbb Z}}$ (\cite[Example~3.6]{Nakamoto-Torii:Hochschild}).

\end{remark}

\section{The $R$-module structure of ${\rm HH}^{\ast}({\rm N}_m(R), {\rm N}_m(R))$}
\label{section:HH-N-N}

%\input{degeneration}
%\subsection{Degeneration of spectral sequences}

In this section, we determine the $R$-module structure of ${\rm HH}^{\ast}({\rm N}_m(R), {\rm N}_m(R))$ for 
$m\ge 3$.  (The case $m=2$ will be discussed in \S\ref{section:m=2}.) 
Throughout this section, we assume that $m \ge 3$. 
Set ${\rm N}={\rm N}_m(R)$, ${\rm B} = {\rm B}_m(R)$, and ${\rm J}=J({\rm N}_m(R))$. 
In \S\ref{subsection:Degss}, we consider a spectral sequence 
converging to the Hochschild cohomology
${\rm HH}^*({\rm N},{\rm N})$. In \S\ref{subsection:freenessE2pqB}, we show that $E_{2}^{p, q}({\rm B})=0$ unless $p=0, m-1$, where $E_1^{p, q}({\rm B}) \Longrightarrow {\rm HH}^{p+q}({\rm N}, {\rm B})$ is a spectral sequence converging to ${\rm HH}^{\ast}({\rm N}, {\rm B})$.   
We also show that $E_{2}^{p, q}({\rm B})$ is a finitely generated free module over $R$ for $p=0, m-1$. 
In \S\ref{subsection:rankE2pqB}, we calculate the rank of $E_{2}^{p, q}({\rm B})$ over $R$ for $p=0, m-1$. 
In \S\ref{subsection:freenessE21qN},  
we show that $E_{2}^{1, q}({\rm N})$ is a finitely generated free module over $R$ for 
the spectral sequence $E_{1}^{p, q}({\rm N}, {\rm N}) \Longrightarrow {\rm HH}^{p+q}({\rm N}, {\rm N})$. 
In \S\ref{subsection:rankE2pqN}, we calculate the rank of $E_{2}^{p, q}({\rm N})$ over $R$ for any $p$.  
As a result, we determine the $R$-module structure of ${\rm HH}^*({\rm N},{\rm N})$.

%In this section 
%we consider spectral sequences converging to
%the Hochschild cohomology
%${\rm HH}^*({\rm N},{\rm N})$.

%In this section, we calculate ${\rm HH}^{\ast}({\rm N}_m(R), {\rm N}_m(R))$. 
%Let $R$ be a commutative ring.
%We set $A={\rm N}_m(R)$.

%\input{degeneration}
\subsection{Degeneration of spectral sequences}\label{subsection:Degss}
Recall the filtration $\{\overline{\rm J}^p{\rm N}\}$ in (\ref{eq:JpM}) or \S\ref{subsection:sssubquotient}. 
Obviously, $\overline{\rm J}^p{\rm N} = {\rm J}^p$ as ideals of ${\rm N}$.  
By regarding ${\rm N}$ as a subobject of the $\mathbb{Z}$-graded
$R$-algebra ${\rm M}={\rm M}_m(R)$, we have 
$F^{p}{\rm N} = \overline{\rm J}^p{\rm N}$, where $F^p{\rm N} = {\rm N} \cap F^p{\rm M}$ has been defined in 
\S\ref{subsection:sssubquotient}. 
Using the filtration
$\{\overline{\rm J}^p{\rm N}\}$,
we have a spectral sequence
\[ {}^JE_1^{p,q}({\rm N},{\rm N})
   \Longrightarrow 
   {\rm HH}^{p+q}({\rm N},{\rm N})\]
of $R$-algebras
by Proposition~\ref{prop:spectral-sequence-A-A-version}.
By the argument in \S\ref{section:HH-N-M/N},
we see that it can be promoted to a spectral sequence
of $\mathbb{Z}$-graded $R$-algebras.

\begin{proposition}\label{prop:HHNNcollapses}
The spectral sequence 
\[ {}^JE_1^{p,q}({\rm N},{\rm N})
   \Longrightarrow 
   {\rm HH}^{p+q}({\rm N},{\rm N})\]
of $\mathbb{Z}$-graded $R$-algebras
collapses from the $E_2$-page.
There is an isomorphism
\begin{eqnarray}
 {\rm HH}^{n}({\rm N}, {\rm N}) = \bigoplus_{n,s}{\rm HH}^{n,s}({\rm N},{\rm N})\cong
   \bigoplus_{n,s}{}^JE_{\infty}^{n-s,s}({\rm N},{\rm N}) \label{eq:HHNNnoextprob} 
\end{eqnarray} 
of bigraded $R$-modules, where ${\rm HH}^{n,s}({\rm N},{\rm N})\cong {}^JE_{\infty}^{n-s,s}({\rm N},{\rm N})$.  
%$R$-algebras.
\end{proposition}

\proof
%We regard ${\rm N}$
%as a subquotient of the $\mathbb{Z}$-graded
%$R$-algebra ${\rm M}={\rm M}_m(R)$.
%Then we have the filtration $\{F^p{\rm N}\}$
%defined in \S\ref{section:HH-N-M/N}.
Since the filtration $\{\overline{\rm J}^p{\rm N}\}$
coincides with $\{F^p{\rm N}\}$,
the spectral sequence $\{{}^JE_r^{p,q}({\rm N},{\rm N})\}_{r\ge 1}$
is isomorphic to 
$\{{}^M{E}_r^{p,q}({\rm N},{\rm N})\}_{r\ge 1}$
given by Proposition~\ref{prop:spectral-sequence-matrix-filtration}.
Therefore, these spectral sequences
collapse from the $E_2$-pages, and 
there are
no extension problems as spectral sequences
of $\mathbb{Z}$-graded $R$-modules
by Theorem~\ref{theorem:collapse-E-2-term-no-extension}.
%Furthermore,
%by Lemma~\ref{lemma:graded_ss_vanishing-simple-version},
%we can see that 
%there are no extension problems
%as spectral sequences of $\mathbb{Z}$-graded
%$R$-algebras.
This completes the proof.
\qed

\begin{remark}\rm 
We can also prove that (\ref{eq:HHNNnoextprob}) is an isomorphism of bigraded $R$-algebras, which will be proved  in 
\S\ref{section:ProductHHNN}.  
\end{remark}

\bigskip 

In the sequel, we omit $J$ of ${}^J E_{1}^{p, q}({\rm N}, {\rm N})$. We also write $E_{r}^{p, q}({\rm N}) = E_{r}^{p, q}({\rm N}, {\rm N})$ and $E_{r}^{p, q, s}({\rm N}) = E_{r}^{p, q, s}({\rm N}, {\rm N})$.  
Here we rephrase Lemma~\ref{lemma:graded_ss_vanishing-simple-version}, which will be used later.  

\begin{lemma}\label{lemma:graded_ss_vanishing}
We have
\[ E_1^{p,q,s}({\rm N}) \cong
   \left\{\begin{array}{cl}
          {\rm N}^{!}_{p+q} \otimes_{R} \mathrm{Gr}^p({\rm N}) & (q=s),\\[2mm]
           0                  & (q\neq s).\\
          \end{array}\right.
\]
\end{lemma}

\bigskip 

When we consider ${\rm B}$ as a ${\Bbb Z}$-graded ${\rm N}$-bimodule, we can 
obtain a spectral sequence 
\[ E_1^{p,q}({\rm B})={\rm HH}^{p+q}({\rm N},\mathrm{Gr}^p({\rm B})) 
   \Longrightarrow {\rm HH}^{p+q}({\rm N},{\rm B}) \] 
by Proposition~\ref{prop:spectral-sequence-matrix-filtration}. 
%If we emphasize the commutative ring $R$, we write $E_{r}^{p,q}({\rm B}; R) = E_{r}^{p, q}({\rm B})$. 
For the ${\Bbb Z}$-graded ${\rm N}$-bimodule ${\rm B}/{\rm N}$, we also obtain a spectral sequence 
\[ E_1^{p,q}({\rm B}/{\rm N})={\rm HH}^{p+q}({\rm N},\mathrm{Gr}^p({\rm B}/{\rm N}))
   \Longrightarrow {\rm HH}^{p+q}({\rm N},{\rm B}/{\rm N}). \]
Then there exist morphisms of spectral sequences 
\[
E_{r}^{\ast, \ast}({\rm N}) \longrightarrow E_{r}^{\ast, \ast}({\rm B}) \longrightarrow E_{r}^{\ast, \ast}({\rm B}/{\rm N}). 
\]
Note that $E_{1}^{p, q}({\rm N}), E_{1}^{p, q}({\rm B})$, and $E_{1}^{p, q}({\rm B}/{\rm N})$ are finitely generated 
free modules over $R$. Indeed, for example, $E_{1}^{p, q}({\rm B}) \cong {\rm N}^{!}_{p+q} \otimes_{R} {\rm Gr}^{p}({\rm B})$ 
is a finitely generated free module over $R$. 
Unless $0 \le p  \le m-1$, $E_{r}^{p, q}({\rm N})=E_{r}^{p, q}({\rm B})=E_{r}^{p, q}({\rm B}/{\rm N})=0$ for $r \ge 1$, since  
${\rm Gr}^{p}({\rm N}) = {\rm Gr}^{p}({\rm B}) = {\rm Gr}^{p}({\rm B}/{\rm N})=0$.  
   
Let us describe results on $E_{r}^{p, q}({\rm B}/{\rm N})$ and ${\rm HH}^{\ast}({\rm N}_m(R),  {\rm B}_m(R)/{\rm N}_m(R))$.  Using 
\[
{\rm Gr}^{p}({\rm B}/{\rm N}) = 
\left\{ 
\begin{array}{cc}
{\rm B}/{\rm N} & (p = 0), \\
0 & (p \neq 0), 
\end{array}
\right.
\] 
we obtain the following theorem.

\begin{theorem}\label{th:HHB/N}
For $r \ge 1$,  
\[
E_{r}^{p, q}({\rm B}/{\rm N}) = \left\{ 
\begin{array}{cc}
{\rm N}^{!}_{q}\otimes_{R} ({\rm B}/{\rm N}) & (p=0), \\
0 & (p\neq 0).  
\end{array}
\right.
\]
For $n \ge 0$, 
\[
{\rm HH}^{n}({\rm N}_m(R), {\rm B}_m(R)/{\rm N}_m(R)) = {\rm N}^{!}_{n} \otimes_{R} ({\rm B}/{\rm N}) \cong 
{\rm N}^{!}_{n} \otimes_{R} R^{m-1}. 
\]

\end{theorem}    
   
In the sequel, if we emphasize the commutative ring $R$, then we write $E_{r}^{p,q}({\rm N}; R) = E_{r}^{p, q}({\rm N})$, 
$E_{r}^{p,q}({\rm B}; R) = E_{r}^{p, q}({\rm B})$,  
and $E_{r}^{p,q}({\rm B}/{\rm N}; R) = E_{r}^{p, q}({\rm B}/{\rm N})$, respectively.

\subsection{The freeness of $E_{2}^{p, q}({\rm B})$}\label{subsection:freenessE2pqB} 
In this subsection, we show that $E_{2}^{p, q}({\rm B})=0$ unless $p=0, m-1$. 
We also show that $E_{2}^{p, q}(B)$ is a finitely generated free module over $R$ for $p=0, m-1$.   
Note that $E_{1}^{p, q}({\rm B}) = {\rm N}^{!}_{p+q}\otimes_{R} {\rm Gr}^{p}({\rm B})$ is a finitely generated 
free module over $R$ for any $p$. 

%Let $B = {\rm B}_m(R)$. 
\bigskip

For $0 \le p \le m-2$,  
\[
\begin{array}{ccl}
%d_{1}^{p, q} : & 
E_{1}^{p, q}({\rm B})  = {\rm N}^{!}_{p+q}\otimes {\rm Gr}^{p}({\rm B}) & \stackrel{d_{1}^{p, q}}{\longrightarrow} & E_{1}^{p+1, q}({\rm B})  = {\rm N}^{!}_{p+q+1}\otimes {\rm Gr}^{p+1}({\rm B}) \\  [1mm]
% &  
f\otimes E_{i, i+p} & \longmapsto & \left\{ 
 	\begin{array}{lc} 
 	(-1)^{p+q+1} fy_{p+1}\otimes E_{1, p+2} & (i=1), \\ [2mm]
 	\begin{array}{l} 
 	  y_{i-1} f \otimes E_{i-1, i+p} \\ [2mm]
 	   +(-1)^{p+q+1} fy_{i+p} \otimes E_{i, i+p+1}
 	 \end{array}  & (2\le i \le m-p-1),  \\ [2mm]
 	y_{m-p-1}f \otimes E_{m-p-1, m} & (i=m-p) 
	\end{array}
\right. 
\end{array}
\]
for $1 \le i \le m-p$. 

\begin{definition}\rm 
For $1 \le p \le m-1$, we define an $R$-homomorphism 
$s^{p, q} : E_{1}^{p, q} \to E_{1}^{p-1, q}$ by 
\[
\begin{array}{ccl}
 E_{1}^{p, q}({\rm B})  = {\rm N}^{!}_{p+q}\otimes {\rm Gr}^{p}({\rm B}) & \stackrel{s^{p, q}}{\longrightarrow} & E_{1}^{p-1, q}({\rm B})  = {\rm N}^{!}_{p+q-1}\otimes {\rm Gr}^{p-1}({\rm B}) \\ 
   f\otimes E_{i, i+p} & \longmapsto & \left\{ 
 	\begin{array}{lc} 
 	f_{R} \otimes E_{i+1, i+p} & 
 		\begin{array}{r} 
 	  	(f=y_if_{R}, f_{R} \in {\rm N}^{!}_{p+q-1}, \\ 1 \le i \le m-p),
 	 	\end{array} \\ [4mm]
 	(-1)^{p+q}f_{L} \otimes E_{1, p}  & 
 		\begin{array}{l}
 		(i=1, f\notin y_1{\rm N}^{!}_{p+q-1}, \\ 
 		f=f_{L}y_{p}, f_{L} \in {\rm N}^{!}_{p+q-1}), 
 		\end{array} 
      \\ [4mm]
 	0 & (\mbox{otherwise})  
	\end{array} 
\right. 
\end{array}
\]
for $1 \le i \le m-p$ and $f \in {\mathcal B}({\rm N}^{!}_{p+q})$.  
\end{definition}

\begin{lemma}\label{lemma:homotopyE1pqR}
For $1\le p \le m-2$, $s^{p+1, q} \circ d_{1}^{p, q} + d_{1}^{p-1, q} \circ s^{p, q} = id_{E_{1}^{p, q}({\rm B})}$.  
\end{lemma}

\proof
Let $1 \le p \le m-2$. 
For $1 \le i \le m-p$ and $f \in {\mathcal B}({\rm N}^{!}_{p+q})$,  
\[
\begin{array}{ccl}
& & s^{p+1, q} \circ d_{1}^{p, q}(f\otimes E_{i, i+p}) \\ 
& = & 
\left\{ 
	\begin{array}{ll}
	(-1)^{p+q+1}f_{R}y_{p+1}\otimes E_{2, p+2} & (i=1, f=y_1f_{R}, f_{R} \in {\rm N}^{!}_{p+q-1}),  \\
	0 & (i=1, f\notin y_1{\rm N}^{!}_{p+q-1}, f=f_{L}y_p, f_{L} \in {\rm N}^{!}_{p+q-1}), \\
	f\otimes E_{1, p+1} & (i=1, f\notin y_1{\rm N}^{!}_{p+q-1}, f\notin {\rm N}^{!}_{p+q-1}y_{p}), \\
	(-1)^{p+q+1}f_{R}y_{i+p}\otimes E_{i+1, i+p+1} & (2 \le i \le m-p-1, f=y_{i}f_{R}, f_{R} \in {\rm N}^{!}_{p+q-1}), \\
	f\otimes E_{i, i+p} & (2 \le i \le m-p-1, f\notin y_{i}{\rm N}^{!}_{p+q-1}), \\ 
	0 & (i=m-p, f=y_{m-p}f_{R}, f_{R}\in {\rm N}^{!}_{p+q-1}), \\ 
	f\otimes E_{m-p, m} & (i=m-p, f \notin y_{m-p} {\rm N}^{!}_{p+q-1}).  
	\end{array}
\right. 
\end{array}
\]
For $1 \le i \le m-p$ and $f \in {\mathcal B}({\rm N}^{!}_{p+q})$,  
\[
\begin{array}{ccl}
& & d_{1}^{p-1, q} \circ s^{p, q} (f\otimes E_{i, i+p}) \\ 
& = & 
\left\{ 
	\begin{array}{ll}
	f\otimes E_{1, p+1} + (-1)^{p+q}f_{R}y_{p+1}\otimes E_{2, p+2} & (i=1, f=y_1f_{R}, f_{R} \in {\rm N}^{!}_{p+q-1}),  \\
	f\otimes E_{1, p+1} & (i=1, f\notin y_1{\rm N}^{!}_{p+q-1}, f=f_{L}y_p, f_{L} \in {\rm N}^{!}_{p+q-1}), \\
	0 & (i=1, f\notin y_1{\rm N}^{!}_{p+q-1}, f\notin {\rm N}^{!}_{p+q-1}y_{p}), \\
	f\otimes E_{i, i+p} +(-1)^{p+q}f_{R} y_{i+p} \otimes E_{i+1, i+p+1} & (2 \le i \le m-p-1, f=y_{i}f_{R}, f_{R} \in {\rm N}^{!}_{p+q-1}), \\
	0& (2 \le i \le m-p-1, f\notin y_{i}{\rm N}^{!}_{p+q-1}), \\ 
	f\otimes E_{m-p, m} & (i=m-p, f=y_{m-p}f_{R}, f_{R}\in {\rm N}^{!}_{p+q-1}), \\ 
	0 & (i=m-p, f \notin y_{m-p} {\rm N}^{!}_{p+q-1}).   
	\end{array}
\right. 
\end{array}
\]
Hence, we see that $s^{p+1, q} \circ d_{1}^{p, q} + d_{1}^{p-1, q} \circ s^{p, q} = id_{E_{1}^{p, q}({\rm B})}$. 
\qed 

\begin{proposition}\label{prop:E2pqB=0}
For $1 \le p \le m-2$, $E_{2}^{p, q}({\rm B}; R) =0$. 
\end{proposition}

\proof 
The statement follows from  Lemma~\ref{lemma:homotopyE1pqR}. 
\qed 

\bigskip 

By Proposition~\ref{prop:E2pqB=0}, $E_{2}^{p, q}({\rm B}; R) =0$ unless $p=0, m-1$. For $p=0, m-1$, we have the following proposition: 

\begin{proposition}\label{prop:E2m-1qfree}
For $p=0, m-1$, $E_{2}^{p, q}({\rm B}; R)$ is a finitely generated free $R$-module.  
\end{proposition}

\proof
Let $p=0$. For $R={\Bbb Z}$, $E_{2}^{0, q}({\rm B}; {\Bbb Z})$ is a finitely generated free ${\Bbb Z}$-module, since 
$E_{2}^{0, q}({\rm B}; {\Bbb Z}) \cong {\rm Ker} d_{1}^{0, q}$ is a submodule of the finitely generated free ${\Bbb Z}$-module $E_{1}^{0, q}({\rm B}; {\Bbb Z}) \cong {\rm N}^{!}_{q} 
\otimes_{{\Bbb Z}} {\rm Gr}^{0}({\rm B})$.   For an arbitrary commutative ring $R$, we have an exact sequence 
\[
0 \longrightarrow E_{2}^{0, q}({\rm B}; {\Bbb Z})\otimes_{{\Bbb Z}} R \longrightarrow E_{2}^{0, q}({\rm B}; R) \longrightarrow {\rm Tor}_{1}^{{\Bbb Z}}(E_{2}^{1, q}({\rm B}; {\Bbb Z}), R) 
\]
by the universal coefficient theorem. 
By Proposition~\ref{prop:E2pqB=0}, $E_{2}^{1, q}({\rm B}; {\Bbb Z})=0$ for $m \ge 3$. 
Hence, $E_{2}^{0, q}({\rm B}; {\Bbb Z})\otimes_{{\Bbb Z}} R \cong E_{2}^{0, q}({\rm B}; R)$. 
Therefore, $E_{2}^{0, q}({\rm B}; R)$ is a finitely generated free $R$-module. 

Let $p=m-1$. By Lemma~\ref{lemma:homotopyE1pqR} and the diagram 
\[
\xymatrix{
	E_{1}^{m-3, q}({\rm B}; R) \ar[r]^{d_{1}^{m-3, q}} \ar[d]_{id}& E_{1}^{m-2, q}({\rm B}; R) \ar[r]^{d_{1}^{m-2, q}}  \ar[d]_{id}\ar[ld]_{s^{m-2,q}}  &  E_{1}^{m-1, q}({\rm B}; R) \ar[ld]_{s^{m-1,q}}\ar[r] \ar[d]_{id}& 0 
	\\
	E_{1}^{m-3, q}({\rm B}; R) \ar[r]^{d_{1}^{m-3, q}} & E_{1}^{m-2, q}({\rm B}; R) \ar[r]^{d_{1}^{m-2, q}}  &  E_{1}^{m-1, q}({\rm B}; R) \ar[r] & 0,
%	0 \ar[r]  & C^{-(m-1), q\quad} \ar[r]^{\delta^{-(m-1),q}} &  C^{-(m-2), q\qquad } \ar[r]^{\qquad\delta^{-(m-2),q}} & &\!\!\!\!\!\!\cdots\!\!\!\!\!\!& \ar[r]^{\delta^{-2,q}} & C^{-1, q} \ar[r]^{\delta^{-1,q}} & C^{0, q} \ar[r]  & 0 
	}
\]
we obtain a short exact sequence 
\[
0 \longrightarrow {\rm Im} \: d_{1}^{m-2, q} \longrightarrow E_{1}^{m-1, q}({\rm B}; R) \longrightarrow E_{2}^{m-1, q}({\rm B}; R) \longrightarrow 0,  
\]
which is split by $s^{m-1, q}$.  Hence, $E_{1}^{m-1, q}({\rm B}; R) \cong {\rm Im} \: d_{1}^{m-2, q} \oplus E_{2}^{m-1, q}({\rm B}; R)$. 
For $R={\Bbb Z}$, $E_{2}^{m-1, q}({\rm B}; {\Bbb Z})$ is finitely generated free 
${\Bbb Z}$-module, since it is isomorphic to a submodule of 
the finitely generated free ${\Bbb Z}$-module $E_{1}^{m-1, q}({\rm B}; {\Bbb Z}) \cong {\rm N}^{!}_{m+q-1}\otimes_{{\Bbb Z}} {\rm Gr}^{m-1}({\rm B})$.  
For an arbitrary commutative ring $R$, we have an exact sequence 
\[
0 \longrightarrow E_{2}^{m-1, q}({\rm B}; {\Bbb Z})\otimes_{{\Bbb Z}} R \longrightarrow E_{2}^{m-1, q}({\rm B}; R)  \longrightarrow {\rm Tor}_{1}^{{\Bbb Z}}(E_{2}^{m, q}({\rm B}; {\Bbb Z}), R) 
\]
by the universal coefficient theorem. By $E_{2}^{m, q}({\rm B}; {\Bbb Z})=0$, $E_{2}^{m-1, q}({\rm B}; R) \cong E_{2}^{m-1, q}({\rm B}; {\Bbb Z})\otimes_{{\Bbb Z}} R$. Hence, $E_{2}^{m-1, q}({\rm B}; R)$ is a finitely generated free $R$-module. 
\qed

\begin{corollary}\label{cor:E2pqB=0}For $m \ge 4$ and $2 \le p \le m-2$,  
$E_{2}^{p, q}({\rm N}; R) =0$. For $m \ge 3$, $E_{2}^{m-1, q}({\rm N}; R)$ is isomorphic to 
the finitely generated free $R$-module $E_{2}^{m-1, q}({\rm B}; R)$. 
\end{corollary}

\proof 
Since $E_{1}^{1, q}({\rm N}) \to E_{1}^{2, q}({\rm N}) \to \cdots \to E_{1}^{m-1, q}({\rm N}) \to 0$ 
is isomorphic to $E_{1}^{1, q}({\rm B}) \to E_{1}^{2, q}({\rm B}) \to \cdots \to E_{1}^{m-1, q}({\rm B}) \to 0$,  
$E_{2}^{p, q}({\rm N}; R) \cong E_{2}^{p, q}({\rm B}; R)$ for $2 \le p \le m-1$. 
The statements follow from Propositions~\ref{prop:E2pqB=0} and \ref{prop:E2m-1qfree}. 
\qed

\begin{remark}\label{remark:E2pq=0unlessp=01m-1}\rm
Recall that $E_{2}^{p, q}({\rm N})=0$ unless $0 \le p \le m-1$. 
Furthermore, we also see that $E_{2}^{p, q}({\rm N})=0$ unless $p=0, 1, m-1$ for $m\ge 3$ 
by Corollary~\ref{cor:E2pqB=0}. 
\end{remark} 

\begin{definition}\label{def:notationxI}\rm 
Let $q > 0$.  For $I = (i_1, i_2, \ldots, i_q)$ with $1 \le i_1, i_2, \ldots, i_q \le m-1$, set $y_{I} = y_{i_1} y_{i_2} \cdots y_{i_q} \in {\rm N}^{!}_{q}$. For $I = (i_1, i_2, \ldots, i_q)$, we define the length $|I|$ of $I$ by $q$.  
We also define 
\[
z(i, I) = y_iy_{I} \otimes E_{i, i} + (-1)^{q+1} y_{I} y_i \otimes E_{i+1, i+1} \in {\rm N}^{!}_{q+1}\otimes_{R} {\rm Gr}^{0}({\rm B}) 
\cong E_{1}^{0, q+1}(B) 
\]
for $1 \le i \le m-1$ and $I$ with $|I|=q$. 

Let $q=0$. We define 
\[
z(i, \emptyset) = z(i) = y_i \otimes E_{i, i} -  y_i \otimes E_{i+1, i+1} \in {\rm N}^{!}_{1}\otimes_{R} {\rm Gr}^{0}({\rm B})
\cong E_{1}^{0, 1}(B) 
\]
for $1 \le i \le m-1$. If $I=\emptyset$, then put $y_{\emptyset} = 1$ and $| \emptyset | =0$. 
For $q\ge 0$, 
set 
\[ {\mathcal Z}(q) = \{ z(i, I) \mid i=1, 2, \ldots, m-1, |I|=q \}. \]
\end{definition}

\begin{proposition}\label{prop:E20qBgenerated}
As an $R$-submodule of $E_{1}^{0, 0}({\rm B}; R)$, 
$E_{2}^{0, 0}({\rm B}; R) = RI_m$.  For $q\ge 0$, 
$E_{2}^{0, q+1}({\rm B}; R) = R\{ z \mid  z\in {\mathcal Z}(q) \}$ as $R$-submodules of $E_{1}^{0, q+1}({\rm B}; R) \cong 
{\rm N}^{!}_{q+1}\otimes_{R} {\rm Gr}^{0}({\rm B})$.  
\end{proposition} 

\proof
Note that 
\[
\begin{array}{ccl}
E_{1}^{0, q+1}({\rm B}; R) \cong {\rm N}^{!}_{q+1}\otimes_{R} {\rm Gr}^{0}({\rm B})  & \stackrel{d_{1}^{0, q+1}}{\longrightarrow} & E_{1}^{1, q+1}({\rm B}; R)\cong {\rm N}^{!}_{q+2}\otimes_{R} {\rm Gr}^{1}({\rm B}) \\ [2mm]
  f\otimes E_{i,i} & \longmapsto & \left\{ 
 \begin{array}{ll}
 (-1)^{q}fy_1\otimes E_{1,2} & (i=1), \\ [2mm]
  	\begin{array}{l} 
  	y_{i-1}f\otimes E_{i-1, i} \\ 
  	+(-1)^{q} fy_i \otimes E_{i, i+1}
  	\end{array}
  	 & (2 \le i \le m-1), \\  [4mm]
  y_{m-1}f\otimes E_{m-1, m} & (i=m) 
 \end{array}
 \right. 
\end{array}
\]
for $f \in {\rm N}^{!}_{q+1}$. 
Let $z = \sum_{i=1}^{m} f_i \otimes E_{i,i} \in E_{1}^{0, q+1}({\rm B}; R)$ with $f_i \in {\rm N}^{!}_{q+1}$ ($1\le i \le m$). Since 
\[
d_{1}^{0, q+1}(z) = \{ (-1)^{q} f_1y_1+ y_1f_2\}\otimes E_{1,2} + %((-1)^{q}f_2x_2+x_2f_3)\otimes E_{2, 3} + 
\cdots + \{ (-1)^{q}f_{m-1}y_{m-1} + y_{m-1}f_{m}\}\otimes E_{m-1, m},  
\]
$z \in {\rm Ker} d_{1}^{0, q+1}$ if and only if 
\begin{eqnarray}
(-1)^{q+1}f_1y_1 = y_1f_2, (-1)^{q+1}f_2y_2=y_2f_3, %(-1)^{q+1} f_3x_3= x_3f_4, 
\ldots, (-1)^{q+1}f_{m-1}y_{m-1} = y_{m-1}f_{m}. \label{eq:condE20q} 
\end{eqnarray}
When $q=-1$, $z \in {\rm Ker} d_{1}^{0, 0}$ if and only if 
$f_1 = f_2 = \cdots = f_m \in R$. Hence, $E_{2}^{0, 0}({\rm B}; R) \cong {\rm Ker} d_{1}^{0, 0} = RI_m$. 

Let us consider $E_{2}^{0, q+1} = {\rm Ker} d_{1}^{0, q+1} \subseteq E_{1}^{0, q+1} \cong {\rm N}^{!}_{q+1}\otimes_{R} {\rm Gr}^{0}({\rm B})$ for $q\ge 0$. It is easy to verify that $z(i, I) \in {\rm Ker} d_{1}^{0, q+1}$. 
Conversely, suppose that $z = \sum_{i=1}^{m} f_i \otimes E_{i,i} \in {\rm Ker} d_1^{0, q+1} \subseteq E_{1}^{0, q+1}({\rm B}; R)$  with $f_i \in {\rm N}^{!}_{q+1}$ ($1\le i \le m$). By (\ref{eq:condE20q}), there exists $g_1, g_2, \ldots, g_{m-1} \in {\rm N}^{!}_{q}$ such that 
\begin{eqnarray*}
f_1 & = & y_1g_1, \\
f_2 & = & (-1)^{q+1}g_1y_1+y_2g_2, \\
\cdots & \cdots & \cdots \\
f_{m-1} & = & (-1)^{q+1}g_{m-2}y_{m-2}+y_{m-1}g_{m-1}, \\
f_{m} & = & (-1)^{q+1}g_{m-1}y_{m-1}.  
\end{eqnarray*}
Then 
\[
z=\{ y_1g_1\otimes E_{1,1}+(-1)^{q+1} g_1y_1\otimes E_{2,2} \}  + \cdots + \{ y_{m-1}g_{m-1}\otimes E_{m-1, m-1}+(-1)^{q+1} g_{m-1}y_{m-1}\otimes E_{m,m}\}. 
\]
Hence, $z \in R\{ z \mid z \in {\mathcal Z}(q) \}$. This completes the proof. 
\qed

\begin{corollary}\label{cor:E201B}
As $R$-modules, $E_{2}^{0, 1}({\rm B}; R) \cong R^{m-1}$. 
\end{corollary}

\proof
By Proposition~\ref{prop:E20qBgenerated}, 
\[ E_{2}^{0, 1}({\rm B}; R) = R\{ z(1), z(2), \ldots, z(m-1) \}. 
\]
It is easy to see that $z(1), z(2), \ldots, z(m-1)$ are linearly independent over $R$. 
Hence, $E_{2}^{0, 1}({\rm B}; R) \cong R^{m-1}$.  
\qed

\begin{proposition}\label{prop:E2m-1qBgenerated} 
As quotient modules of $E_1^{m-1, -(m-1)+q}({\rm B}) \cong {\rm N}^{!}_{q}\otimes_{R} {\rm Gr}^{m-1}({\rm B})$ over $R$, 
\[ E_2^{m-1,-(m-1)+q}({\rm B})=
   R\{y_I\otimes E_{1,m}|\ 
   I=(i_1,\ldots,i_q),i_1\neq 1, i_q\neq m-1\}\]
for $q \ge 0$. Here, we understand that $E_2^{m-1,-(m-1)}({\rm B})=R\{ 1\otimes E_{1, m} \}$ when $q=0$. 
Moreover, $\{y_I\otimes E_{1,m}|\  I=(i_1,\ldots,i_q),i_1\neq 1, i_q\neq m-1\}$ are linearly independent in 
$E_2^{m-1,-(m-1)+q}({\rm B})$ over $R$. 
\end{proposition}
\proof
Let us consider 
\[
d_{1}^{m-2, -(m-1)+q} : E_{1}^{m-2, -(m-1)+q}\cong {\rm N}_{q-1}\otimes_{R} {\rm Gr}^{m-2}({\rm B}) \longrightarrow E_{1}^{m-1, -(m-1)+q} \cong {\rm N}_{q}\otimes_{R} {\rm Gr}^{m-1}({\rm B}).
\] Since 
\begin{eqnarray*}
d_{1}^{m-2, -(m-1)+q}(y_{J}\otimes E_{1, m-1}) & = & (-1)^{q} y_{J}y_{m-1} \otimes E_{1, m}, \\
d_{1}^{m-2, -(m-1)+q}(y_{J}\otimes E_{2, m}) & = & y_{1}y_{J} \otimes E_{1, m}  
\end{eqnarray*} 
for $|J|=q-1$, we have 
\[
{\rm Im} \: d_{1}^{m-2, -(m-1)+q} = R\{y_I\otimes E_{1,m}|\ 
   I=(i_1,\ldots,i_q),i_1=1\mbox{ or }  i_q= m-1\}. 
\]
Hence, we can see that the statement is true. 
\qed

\subsection{The rank of $E_{2}^{p, q}({\rm B})$}\label{subsection:rankE2pqB} 

In this subsection, 
we calculate the rank
of the free $R$-module $E_2^{p,q}({\rm B})$.
We note that
${\rm rank}_RE_2^{p,q}({\rm B})=0$
unless $p=0$ or $m-1$ by Proposition~\ref{prop:E2pqB=0}. 
%For $q\in\mathbb{Z}$,
%we set 
Recall 
\[ \varphi(q)=\,{\rm rank}_R {\rm N}^!_q. \]
for $q\in\mathbb{Z}$. 
We shall show the following theorem: 

\begin{theorem}\label{th:rankE2pqB}
For each $q\in\mathbb{Z}$, 
we have
\[ \begin{array}{rcl}
    {\rm rank}_RE_2^{p,q}({\rm B})&=&
    \left\{\begin{array}{ll}
           \varphi(q)& (p=0), \\[2mm]
           \displaystyle
           (-1)^m\varphi(q)+\sum_{k=0}^{m-1}(-1)^k(k+1)
           \varphi(q+m-k-1)& (p=m-1).\\
           \end{array}\right.
   \end{array}\]
\end{theorem}

%\newpage

\subsubsection{The rank of $E_2^{0,q}({\rm B})$}

%In this subsection 
First, we calculate the rank of $E_2^{0,q}({\rm B})$. 
In this subsubsection, we show the following theorem, which claims the formula of ${\rm rank}_{R} E_2^{0, q}$   
in Theorem~\ref{th:rankE2pqB}. 

\begin{theorem}\label{theorem:rank-E-2-0-q-B}
For each $q\in\mathbb{Z}$,
we have
\[ {\rm rank}_RE_2^{0,q}({\rm B})=\varphi(q).\]
\end{theorem}

%\newpage
%\bigskip 

Recall $y_I$ in Definition~\ref{def:notationxI}. 
For $q\ge 0$,
we have $\varphi(q)=\sharp{{\mathcal B}}(q)$,
where
\[ \mathcal{B}(q)=\{ 0 \neq y_I\in {\rm N}^!_q|\ I=(i_1,\ldots,i_q)\}. \]
We put 
\begin{eqnarray*}
\varphi(q; i_1\neq a) & = &  \sharp{{\mathcal B}(q, i_1\neq a)}, \\ 
\varphi(q; i_1\neq a, i_q = b) & = & \sharp{{\mathcal B}(q, i_1\neq a, i_q = b)},
\end{eqnarray*} 
and so on, where 
\begin{eqnarray*} 
\mathcal{B}(q;i_1\neq a) & = & \{y_I\in \mathcal{B}(q) |\ i_1\neq a\}, \\ 
\mathcal{B}(q;i_1\neq a,i_q = b) & = & \{y_I\in \mathcal{B}(q) |\ i_1\neq a,i_q = b\}, 
\end{eqnarray*}  
and so on. For $\varphi(q;i_1\neq 1)$, we have the following lemma.

\begin{lemma}\label{lemma:varphi-neq-1}
For $q\ge 0$,
we have
\[ \varphi(q;i_1\neq 1)= 
   \sum_{r=0}^{m-1}(-1)^r\varphi(q-r),\]
where $\varphi(0;i_1\neq 1)=1$.
\end{lemma}

\proof
When $q=0$, we have
\[ \varphi(0;i_1\neq 1)=1=\varphi(0).\]

We consider the case when $0<q\le m-1$.
We have
\[ \begin{array}{cl}
    &\varphi(q;i_1\neq 1)\\[2mm]
   =&\varphi(q)-\varphi(q;i_1=1)\\[2mm]
   =&\varphi(q)-\varphi(q-1;i_1\neq 2)\\[2mm]
   =&\varphi(q)-\varphi(q-1)+\varphi(q-1;i_1=2)\\[2mm]
   =&\varphi(q)-\varphi(q-1)+\varphi(q-2;i_1\neq 3)\\[2mm]
   =&\cdots\\[2mm]
   =&\displaystyle\sum_{s=0}^{t-1}(-1)^s\varphi(q-s)
    +(-1)^t\varphi(q-t;i_1\neq t+1)\\[4mm]
   =&\displaystyle\sum_{s=0}^{q-2}(-1)^s\varphi(q-s)
    +(-1)^{q-1}\varphi(1;i_1\neq q)\\[4mm]
   =&\displaystyle\sum_{s=0}^{q}(-1)^s\varphi(q-s).\\
   \end{array}\]

We consider the case when $q\ge m$.
We have
\[ \begin{array}{cl}
    &\varphi(q;i_1\neq 1)\\[2mm]
   =&\displaystyle\sum_{s=0}^{t-1}(-1)^s\varphi(q-s)
    +(-1)^t\varphi(q-t;i_1\neq t+1)\\[4mm]
   =&\displaystyle\sum_{s=0}^{m-3}(-1)^s\varphi(q-s)
    +(-1)^{m-2}\varphi(q-(m-2);i_1\neq m-1)\\[4mm]
   =&\displaystyle\sum_{s=0}^{m-1}(-1)^s\varphi(q-s).\\
   \end{array}\]
This completes the proof. 
\qed

\bigskip 

By Lemma~\ref{lemma:varphi-relation}, we have another formula for $\varphi(q;i_1\neq 1)$. 

\begin{lemma}\label{lemma:varphi-neq-1-reformulation}
For $q>0$, we have
\[ \varphi(q;i_1\neq 1)=\sum_{r=1}^{m-2}(-1)^{r-1}(m-1-r)
   \varphi(q-r).\]
\end{lemma}

\proof
%By Remark~\ref{remark:phi-neq-i-independent}, we have 
By Lemma~\ref{lemma:varphi-neq-1}, we have 
\[ \varphi(q;i_1\neq 1)= 
   \varphi(q)+\sum_{r=1}^{m-1}(-1)^r\varphi(q-r).\]
Using Lemma~\ref{lemma:varphi-relation},
we obtain 
\[ \begin{array}{cl}
    &\varphi(q;i_1\neq 1)\\[2mm]
   =&\displaystyle\sum_{r=1}^{m-1}(-1)^{r-1}(m-r)\varphi(q-r)
   +\sum_{r=1}^{m-1}(-1)^r\varphi(q-r)\\[4mm]
   =&\displaystyle\sum_{r=1}^{m-1}(-1)^{r-1}(m-(r+1))\varphi(q-r)\\[4mm] 
   =&\displaystyle\sum_{r=1}^{m-2}(-1)^{r-1}(m-1-r)\varphi(q-r). 
   \end{array}\]
This completes the proof. 
\qed

%\bigskip 
%
%We put
%\[ \begin{array}{rcl}
%    \mathcal{B}(q;i_1=a,i_q=b)&=&
%    \{y_I\in \mathcal{B}(q) |\ i_1=a,i_q=b\}\\[2mm]
%    \mathcal{B}(q;i_1\neq a,i_q=b)&=&
%    \{y_I\in \mathcal{B}(q) |\ i_1\neq a,i_q=b\}\\[2mm]
%    \mathcal{B}(q;i_1=a,i_q\neq b)&=&
%    \{y_I\in \mathcal{B}(q) |\ i_1=a,i_q\neq b\}\\[2mm]
%    \mathcal{B}(q;i_1\neq a,i_q\neq b)&=&
%    \{y_I\in \mathcal{B}(q) |\ i_1\neq a,i_q\neq b\}. \\[2mm]
%   \end{array}\]
%%We define
%%\[ z(i,I)=y_iy_I\otimes E_{i,i}
%%   +(-1)^{q+1}y_Iy_i\otimes E_{i+1,i+1}
%%   \in E_1^{0,q+1}(B)\]
%%for $i=1,\ldots,m-1$, 
%%where $|I|=q$.

\bigskip 

Now, we calculate ${\rm rank}_RE_2^{0,q}({\rm B})$ for $q \ge 0$.   
Recall $z(i, I)$ in Definition~\ref{def:notationxI}. 
By Proposition~\ref{prop:E20qBgenerated},  
\[ E_2^{0,q+1}({\rm B})=R\{z|\ z\in\mathcal{Z}(q)\}\]
for $q\ge 0$, where 
\[ \mathcal{Z}(q)=\{z(i,I)|\ i=1,\ldots,m-1,|I|=q\}.\]
For $1\le i\le m-1$, we set
\[ \begin{array}{rcl}
    \mathcal{Z}(q;i)_1&=&
    \{z(i,I)\in\mathcal{Z}(q)|\ i_1\neq i+1,i_q\neq i-1\}, \\[2mm]
    \mathcal{Z}(q;i)_2&=&
    \{z(i,I)\in\mathcal{Z}(q)|\ i_1\neq i+1,i_q=i-1\}, \\[2mm]
    \mathcal{Z}(q;i)_3&=&
    \{z(i,I)\in\mathcal{Z}(q)|\ i_1=i+1,i_q\neq i-1\}, \\[2mm]
    \mathcal{Z}(q;i)_4&=&
    \{z(i,I)\in\mathcal{Z}(q)|\ i_1=i+1,i_q=i-1\}. \\[2mm]
   \end{array}\]
Here we understand 
$\mathcal{Z}(q;1)_2 = \mathcal{Z}(q;m)_2 = \mathcal{Z}(q;0)_3 = \mathcal{Z}(q;m-1)_3 = \mathcal{Z}(q;1)_4 = \mathcal{Z}(q;m-1)_4 = \emptyset$. 
We have a decomposition
\[ \mathcal{Z}(q)  = \bigcup_{r=1}^4
    \bigcup_{i=1}^{m-1}\mathcal{Z}(q;i)_r.\]
For any $z(i,I)\in \mathcal{Z}(q;i)_4$,
we have $z(i,I)=0$. 
%We can easily obtain the following lemma.
%
%\begin{lemma}
%For any $z(i,I)\in \mathcal{Z}(q;i)_4$,
%we have
%\[ z(i,I)=0.\]
%\end{lemma}
%
%We note that 
Note that 
\[ z(i,I)=y_iy_I\otimes E_{i,i} \]
for $z(i,I)\in \mathcal{Z}(q;i)_2$,
and
\[ z(i,I)= (-1)^{q+1}y_Iy_i\otimes E_{i+1,i+1}\]
for $z(i,I)\in \mathcal{Z}(q;i)_3$.
It is easy to see that 
\[
R\{  z(i, I)  \mid z(i, I) \in \mathcal{Z}(q;i)_2 \} = R\{  z(i-1, I)  \mid z(i-1, I) \in \mathcal{Z}(q;i)_3 \}
\]  
for $2 \le i  \le m-1$ and $q \ge 0$. Hence, we have 
\[
E_2^{0,q+1}({\rm B})=R\{z|\ z\in\bigcup_{i=1}^{m-1}\mathcal{Z}(q; i)_1\cup \bigcup_{i=2}^{m-1}\mathcal{Z}(q; i)_2  \}.  
\]

%\bigskip 

%%%%%%%%%%%%%%%%   the beginning of if0 
\if0 

We shall construct a bijection $\Phi$
between $\mathcal{Z}(q;i)_2$ and $\mathcal{Z}(q;i-1)_3$.
We write
\[ (i,I)=(I',i-1) \]
for $I=(i_1,\ldots,i_q)$ with $i_1\neq i+1$ and $i_q=i-1$.
For $1\le i\le m$,
we define a map
$\Phi: \mathcal{Z}(q;i)_2\to\mathcal{Z}(q;i-1)_3$ by 
\[ \Phi(z(i,I))=z(i-1,I').\]
Note that we understand
$\mathcal{Z}(q;m)_2=\mathcal{Z}(q;0)_3=\emptyset$.

We can easily obtain the following lemma.

\begin{lemma}
For $1\le i\le m$,
the map $\Phi: \mathcal{Z}(q;i)_2\to\mathcal{Z}(q;i-1)_3$
is a bijection.
\end{lemma}

\begin{corollary}
For $1\le i\le m$, we have
\[ \sharp\mathcal{Z}(q;i)_2=\sharp\mathcal{Z}(q,i-1)_3.\]
\end{corollary}

\if0
\begin{corollary}
We have
\[ {\rm rank}_RE_2^{0,q+1}(B)=\sum_{r=1}^2\sum_{i=1}^{m-1}
   \sharp\mathcal{Z}(q;i)_r.\]
\end{corollary}
\fi

Note that
\[ \begin{array}{rcl}
   \sharp\mathcal{Z}(q;i)_1=
   \sharp\mathcal{B}(q;i_1\neq i+1,i_q\neq i-1),\\[2mm]   
   \sharp\mathcal{Z}(q;i)_2=
   \sharp\mathcal{B}(q;i_1\neq i+1,i_q=i-1),\\[2mm]   
   \sharp\mathcal{Z}(q;i)_3=
   \sharp\mathcal{B}(q;i_1= i+1,i_q\neq i-1).\\[2mm]   
    \end{array}\]
We easily obtain the following lemma.

\begin{lemma}
For $1\le i\le m-1$, we have
\[ \varphi(q)=\sum_{r=1}^3\mathcal{Z}(q;i)_r
   +\sharp\mathcal{B}(q;i_1=i+1,i_q=i-1).\]
\end{lemma}

\fi
%%%%%%%%%%%%%%%%%%%%%   End of \if0 

\begin{proposition}\label{prop:rank-E-2-0-q-i-neq-description}
For $q>0$,
we have
\[ {\rm rank}_RE_2^{0,q+1}({\rm B})=
   (m-1)\varphi(q)-\varphi(q;i_1\neq 1). \]
\end{proposition}

%Note that we have
%\[ {\rm rank}_RE_2^{0,1}(B)=m-1.\]

\proof
We can easily see that 
\[ 
\{ z \mid 0 \neq z \in \bigcup_{i=1}^{m-1}\mathcal{Z}(q; i)_1\cup \bigcup_{i=2}^{m-1}\mathcal{Z}(q; i)_2 \}
\] 
is linearly independent over $R$. Note that $\mathcal{Z}(q; 1)_2=\emptyset$. 
Hence, we have
\[ \begin{array}{rcl}
   {\rm rank}_RE_2^{0,q+1}({\rm B})&=& \displaystyle 
   \sum_{r=1}^2\sum_{i=1}^{m-1}
   \sharp\{ z \mid 0 \neq z \in \mathcal{Z}(q;i)_r\} \\[4mm] 
   &=& 
   \displaystyle
   \sum_{r=1}^3\sum_{i=1}^{m-1}
   \sharp\{ z \mid 0 \neq  z \in \mathcal{Z}(q;i)_r\} 
    -\sum_{i=1}^{m-1}\sharp\{ z \mid 0 \neq z \in \mathcal{Z}(q;i)_3\} \\[4mm]
   &=&\displaystyle
   \sum_{i=1}^{m-1}(\varphi(q)
   -\sharp\mathcal{B}(q;i_1=i+1,i_q=i-1))
   -\sum_{i=1}^{m-1}\sharp\mathcal{B}(q;i_1=i+1,i_q\neq i-1) \\[4mm] 
   &=&\displaystyle
   (m-1)\varphi(q)
   -\sum_{i=1}^{m-1}\sharp\mathcal{B}(q;i_1=i+1,i_q=i-1)\\[4mm]
   &&\displaystyle\phantom{(m-1)\varphi(q)}
   -\sum_{i=1}^{m-1}\sharp\mathcal{B}(q;i_1=i+1,i_q\neq i-1)\\[4mm]
   &=&\displaystyle
   (m-1)\varphi(q)-\sum_{i=1}^{m-1}\sharp\mathcal{B}(q;i_1=i+1)\\[4mm]
   &=&
   (m-1)\varphi(q)-\varphi(q;i_1\neq 1).\\
   \end{array}\]
This completes the proof. 
\qed

\proof[Proof of Theorem~\ref{theorem:rank-E-2-0-q-B}]
For $q<0$,
we have
\[ {\rm rank}_RE_2^{0,q}({\rm B})=0=\varphi(q). \]
For $q=0$,
we have
\[ {\rm rank}_RE_2^{0,0}({\rm B})=1=\varphi(0) \]
by Proposition~\ref{prop:E20qBgenerated}. For $q=1$,
we have
\[ {\rm rank}_RE_2^{0,1}({\rm B})=m-1=\varphi(1) \]
by Corollaries~\ref{cor:varphi} and \ref{cor:E201B}. %Here we used Corollary~\ref{cor:varphi}. 

We assume that $q>1$.
By Proposition~\ref{prop:rank-E-2-0-q-i-neq-description},
Lemmas~\ref{lemma:varphi-relation}
and \ref{lemma:varphi-neq-1-reformulation},
we obtain
\[ \begin{array}{rcl}
    {\rm rank}_RE_2^{0,q}({\rm B})&=&(m-1)\varphi(q-1)
     -\varphi(q-1;i_1\neq 1)\\[2mm]
    &=&\displaystyle (m-1)\varphi(q-1)+
     \sum_{r=1}^{m-2}(-1)^r(m-1-r)\varphi(q-1-r)\\[2mm]
    &=&\varphi(q).\\
  \end{array} \]
%Note that
%\[ \begin{array}{rcl}
%    \varphi(q)&=&
%    \displaystyle\sum_{i=1}^{m-1}\varphi(q;i_1=i)\\[4mm]
%    &=&\displaystyle\sum_{i=1}^{m-2}\varphi(q-1;i_1\neq i+1)
%        +\varphi(q-1)\\[4mm]
%    &=&\displaystyle\sum_{i=1}^{m-2}(\varphi(q-1)-\varphi(q-1;i_1=i+1))
%       +\varphi(q-1)\\[4mm]
%    &=&\displaystyle (m-1)\varphi(q-1)-\varphi(q-1;i_1\neq 1).\\ 
%   \end{array}\]
%Hence, 
%by Proposition~\ref{prop:rank-E-2-0-q-i-neq-description},
%we obtain 
%\[ {\rm rank}_RE_2^{0,q}({\rm B})=\varphi(q).\]
Therefore, we have proved 
Theorem~\ref{theorem:rank-E-2-0-q-B}.  
\qed

\subsubsection{The rank of $E_2^{m-1,q}({\rm B})$}

%In this subsection 
Next, we calculate the rank of the $E_2^{m-1,q}({\rm B})$. 
In this subsubsection, we show the following theorem, which claims the equivalent 
formula of ${\rm rank}_{R} E_{2}^{m-1, q}$ in  
Theorem~\ref{th:rankE2pqB}. 

\begin{theorem}\label{th:E2m-1,q}
For each $q\in\mathbb{Z}$,
we have
\[ {\rm rank}_RE_2^{m-1,-(m-1)+q}({\rm B})=
   \sum_{k=0}^{m-1}(-1)^k (k+1)\varphi(q-k)
   +(-1)^m\varphi(q-m+1).\]
\end{theorem}

\if0
\begin{remark}\rm
Since $E_2^{m-1,q}(A)\cong E_2^{m-1,q}({\rm B})$,
we have
\[ {\rm rank}_RE_2^{m-1,-(m-1)+q}(A)=
   \sum_{k=0}^{m-1}(-1)^k (k+1)\varphi(q-k)
   +(-1)^m\varphi(q-m+1).\]
\end{remark}
\fi

%\newpage

%\subsubsection{The case when $0\le q<m-1$}

For $q<0$, Theorem~\ref{th:E2m-1,q} is true since 
the both sides of the formula are $0$. It suffices to prove Theorem~\ref{th:E2m-1,q} for $q\ge 0$. 
Let us consider the case when $0\le q<m-1$. 

\begin{proposition}\label{prop:rankE2-B-m-1<0}
For $0\le q<m-1$,
we have
\[ {\rm rank}_RE_2^{m-1,-(m-1)+q}({\rm B})=
   \sum_{r=0}^q(-1)^r \varphi(q-r;i_1\neq 1),\]
where $\varphi(0;i_1\neq 1)=1$.
\end{proposition}

\proof
By Proposition~\ref{prop:E2m-1qBgenerated},  
\[ E_2^{m-1,-(m-1)+q}({\rm B})=
   R\{y_I\otimes E_{1,m}|\ 
   I=(i_1,\ldots,i_q),i_1\neq 1, i_q\neq m-1\}.\]
Since $\{ y_I\otimes E_{1,m} |\ I=(i_1,\ldots,i_q),i_1\neq 1, i_q\neq m-1, y_{I} \neq 0 \}$ is linearly independent over $R$,   
we have
\[ \begin{array}{cl}
    &{\rm rank}_RE_2^{m-1,-(m-1)+q}({\rm B})\\[2mm]
   =&\varphi(q;i_1\neq 1,i_q\neq m-1)\\[2mm]
   =&\varphi(q;i_1\neq 1)-\varphi(q; i_1\neq 1,i_q=m-1)\\[2mm]
   =&\varphi(q;i_1\neq 1)-\varphi(q-1;i_1\neq 1;i_{q-1}\neq m-2)\\[2mm]
   =&\varphi(q;i_1\neq 1)-
    \varphi(q-1;i_1\neq 1)+\varphi(q-1;i_1\neq 1,i_{q-1}=m-2)\\[2mm]
   =&\varphi(q;i_1\neq 1)-\varphi(q-1;i_1\neq 1)
    +\varphi(q-2;i_1\neq 1,i_{q-2}\neq m-3)\\[2mm]
   =&\cdots\\[2mm]
   =&\displaystyle\sum_{r=0}^{s-1}(-1)^r\varphi(q-r;i_1\neq 1)
    +(-1)^s\varphi(q-s;i_1\neq 1,i_{q-s}\neq m-(s+1))\\[4mm]
   =&\displaystyle\sum_{r=0}^{q-2}(-1)^r\varphi(q-r;i_1\neq 1)
    +(-1)^{q-1}\varphi(1;i_1\neq 1,i_1\neq m-q)\\[4mm]
   =&\displaystyle\sum_{r=0}^q (-1)^r\varphi(q-r;i_1\neq 1),  
   \end{array}\]
which is what we wanted. 
\qed

\begin{proposition}\label{prop:E2pqBq<m-1}
For $0\le q< m-1$,
we have
\[ {\rm rank}_RE_2^{m-1,-(m-1)+q}({\rm B})=
   \sum_{k=0}^q(-1)^k(k+1)\varphi(q-k).\]
\end{proposition}

\proof
By Lemma~\ref{lemma:varphi-neq-1} and 
Proposition~\ref{prop:rankE2-B-m-1<0},  
we have 
\[ \begin{array}{rcl}
    {\rm rank}_RE_2^{m-1,-(m-1)+q}({\rm B})&=&\displaystyle
    \sum_{r=0}^q(-1)^r\sum_{s=0}^{q-r}(-1)^s\varphi(q-r-s)\\[4mm]
    &=&\displaystyle
    \sum_{k=0}^q(-1)^k(k+1)\varphi(q-k).\\
   \end{array}\]
This completes the proof. 
\qed

\bigskip 

%\newpage

%\subsubsection{The case when $q=m-1$}

Let us  consider the case when $q=m-1$. 

\begin{proposition}\label{prop:rank-E-2-B-m-1-0}
We have
\[ {\rm rank}_RE_2^{m-1,0}({\rm B})=
   \sum_{r=0}^{m-2}(-1)^r\varphi(m-1-r;i_1\neq 1).\]
\end{proposition}

\proof
By the proof of Proposition~\ref{prop:rankE2-B-m-1<0},
we have
\[ \begin{array}{cl}
   &{\rm rank}_RE_2^{m-1,0}({\rm B})\\[2mm]
  =&\varphi(m-1;i_1\neq 1,i_{m-1}\neq m-1)\\[2mm]
  =&\displaystyle\sum_{r=0}^{m-3}(-1)^r\varphi(m-1-r;i_1\neq 1)
    +(-1)^{m-2}\varphi(1;i_1\neq 1,i_1\neq 1)\\[4mm]
  =&\displaystyle\sum_{r=0}^{m-2}(-1)^r\varphi(m-1-r;i_1\neq 1), \\
 \end{array}\]
which is what we wanted. 
\qed

\begin{proposition}\label{prop:E2pqBq=m-1}
We have
\[ {\rm rank}_RE_2^{m-1,0}({\rm B})=
   (-1)^m+\sum_{k=0}^{m-1}(-1)^k(k+1)\varphi(m-1-k).\]
\end{proposition}

\proof
%By Lemma~\ref{lemma:phi-q-r-for-q-le-m-1}
By Lemma~\ref{lemma:varphi-neq-1}
and Proposition~\ref{prop:rank-E-2-B-m-1-0},
we have
\[ \begin{array}{cl}
     &{\rm rank}_RE_2^{m-1,0}({\rm B})\\[2mm]
    =&\displaystyle\sum_{r=0}^{m-2}(-1)^r\varphi(m-1-r;i_1\neq 1)\\[4mm]
    =&\displaystyle\sum_{r=0}^{m-2}(-1)^r
      \sum_{s=0}^{m-1-r}(-1)^s\varphi(m-1-r-s)\\[4mm]
    =&\displaystyle\sum_{k=0}^{m-1}(-1)^k(k+1)\varphi(m-1-k)
      -(-1)^{m-1}\varphi(0).\\
   \end{array}\]
This completes the proof. 
\qed

%\newpage

%\subsubsection{The case when $q\ge m$}

\bigskip 

Finally, let us consider 
the case when $q\ge m$. 

\begin{proposition}\label{rank-e-2-B-m-1-positive}
For $q\ge m$,
we have
\[ {\rm rank}_RE_2^{m-1,-(m-1)+q}({\rm B})=
   \sum_{r=0}^{m-1}(-1)^r\varphi(q-r;i_1\neq 1).\]
\end{proposition}

\proof
%Since
%\[ E_2^{m-1,-(m-1)+q}({\rm B})=
%   R\{y_I\otimes E_{1,m}|\ 
%   I=(i_1,\ldots,i_q),i_1\neq 1, i_q\neq m-1\},\] 
As in the proof of Proposition~\ref{prop:rankE2-B-m-1<0}, 
we have
\[ \begin{array}{cl}
    &{\rm rank}_RE_2^{m-1,-(m-1)+q}({\rm B})\\[2mm]
   =&\varphi(q;i_1\neq 1,i_q\neq m-1)\\[2mm]
   =&\displaystyle\sum_{r=0}^{s-1}(-1)^r\varphi(q-r;i_1\neq 1)
    +(-1)^s\varphi(q-s;i_1\neq 1,i_{q-s}\neq m-(s+1))\\[4mm]
   =&\displaystyle\sum_{r=0}^{m-3}(-1)^r\varphi(q-r;i_1\neq 1)
    +(-1)^{m-2}\varphi(q-(m-2);i_1\neq 1,i_{q-(m-2)}\neq 1)\\[4mm]
   =&\displaystyle\sum_{r=0}^{m-2}(-1)^r\varphi(q-r;i_1\neq 1)
    +(-1)^{m-1}\varphi(q-(m-2);i_1\neq 1,i_{q-(m-2)}=1)\\[4mm]
   =&\displaystyle\sum_{r=0}^{m-1} (-1)^r\varphi(q-r;i_1\neq 1),  
   \end{array}\]
which is what we wanted. 
\qed

\begin{proposition}\label{prop:E2pqBq>m-1} 
For $q\ge m$,
we have
\[ {\rm rank}_RE_2^{m-1,-(m-1)+q}({\rm B})=
   \sum_{k=0}^{m-1}(-1)^k (k+1)\varphi(q-k)
   +(-1)^m\varphi(q-m+1).\]
\end{proposition}

\proof
By Lemma~\ref{lemma:varphi-neq-1}
and %Remark~\ref{remark:phi-neq-i-independent},
Proposition~\ref{rank-e-2-B-m-1-positive},
we have
\[ \begin{array}{cl}
     &{\rm rank}_RE_2^{m-1,-(m-1)+q}({\rm B})\\[2mm]
    =&\displaystyle\sum_{r=0}^{m-1}(-1)^r\varphi(q-r;i_1\neq 1)\\[4mm]
    =&\displaystyle\sum_{r=0}^{m-1}(-1)^r
      \sum_{s=0}^{m-1}(-1)^s\varphi(q-r-s)\\[4mm]
    =&\displaystyle\sum_{k=0}^{2(m-1)}(-1)^kC_k\varphi(q-k),\\
   \end{array}\]
where 
\[ C_k=\left\{ \begin{array}{ll}
                k+1& (0\le k\le m-1), \\[2mm]
                2(m-1)-k+1 & (m\le k\le 2(m-1)).\\
               \end{array}\right. \]
By Lemma~\ref{lemma:varphi-relation},
we have 
\[ \begin{array}{cl}
    &\displaystyle\sum_{k=m}^{2(m-1)}(-1)^kC_k\varphi(q-k)\\[2mm]
   =&
    \displaystyle
    (-1)^m\sum_{k=1}^{m-1}(-1)^{k-1}(m-k)\varphi((q-m+1)-k)\\[4mm]
   =&(-1)^m\varphi(q-m+1).\\
   \end{array}\]
This completes the proof. 
%\[ {\rm rank}_RE_2^{m-1,-(m-1)+q}(B)
%   =\displaystyle\sum_{k=0}^{m-1}(-1)^k(k+1)\varphi(q-k)
%   +(-1)^m\varphi(q-m+1).\]
\qed

%\newpage

\bigskip 

Hence, we have proved Theorem~\ref{th:E2m-1,q} by Propositions~\ref{prop:E2pqBq<m-1}, \ref{prop:E2pqBq=m-1}, and \ref{prop:E2pqBq>m-1}. Therefore, we have finished the proof of Theorem~\ref{th:rankE2pqB}. 

\begin{remark}\rm
We consider the cochain complex
\[ (E_1^{*,-(m-1)+q}({\rm B}),d_1),\]
where
\[ E_1^{(m-1)-r,-(m-1)+q}({\rm B})\cong
   {\rm N}^!_{q-r}\otimes_R {\rm Gr}^{(m-1)-r}({\rm B}). \]
Since
\[ {\rm rank}_RE_1^{(m-1)-r,-(m-1)+q}({\rm B})=
    (r+1)\varphi(q-r)\]
for $0\le r\le m-1$,
the Euler characteristic
of $(E_1^{*,-(m-1)+q}({\rm B}),d_1)$ is given by
\[ \chi(E_1^{*,-(m-1)+q}({\rm B}))=
   (-1)^{m-1}\sum_{r=0}^{m-1} (-1)^r (r+1)\varphi(q-r).\]
For $q \in {\mathbb Z}$,  
we can directly verify that 
   \begin{eqnarray*}
%   (-1)^{m-1}\sum_{r=0}^q(-1)^r(r+1)\varphi(q-r)=
   \chi(E_2^{*,-(m-1)+q}({\rm B})) 
   & = & {\rm rank}_{R}E_{2}^{0, -(m-1)+q}({\rm B}) +(-1)^{m-1}{\rm rank}_RE_2^{m-1,-(m-1)+q}({\rm B}) \\
   & = &  \chi(E_1^{*,-(m-1)+q}({\rm B}))
   \end{eqnarray*}
by Theorems~\ref{theorem:rank-E-2-0-q-B} and \ref{th:E2m-1,q}.  
\end{remark}

\bigskip 

Summarizing the discussion in \S\ref{subsection:rankE2pqB}, we have the following theorems by Theorems~\ref{theorem:collapse-E-2-term-no-extension} and \ref{th:rankE2pqB}. 

\begin{theorem}
Let $m \ge 3$. The Hochschild cohomology ${\rm HH}^{n}({\rm N}_m(R), {\rm B}_m(R))$ is a free $R$-module for $n\ge 0$. For $n \ge 0$, the rank of ${\rm HH}^{n}({\rm N}_m(R), {\rm B}_m(R))$ is given by 
\begin{eqnarray*} 
& & {\rm rank}_{R} {\rm HH}^{n}({\rm N}_m(R), {\rm B}_m(R)) \\ 
& = & 2\varphi(n)+(-1)^{m-1}(m-1)\varphi(n-m+1) + \sum_{k=1}^{m-2} (-1)^k (k+1)\varphi(n-k).  
\end{eqnarray*} 
\end{theorem}

\begin{theorem}
Let $m \ge 3$. For $n \ge 0$ and $s \in {\Bbb Z}$, ${\rm HH}^{n, s}({\rm N}_m(R), {\rm B}_m(R))$ is a free $R$-module. 
Furthermore, 
\[ {\rm HH}^{n,s}({\rm N}_m(R), {\rm B}_m(R))\cong E_2^{n-s,s}({\rm B}_m(R))\] 
as $R$-modules and 
\begin{eqnarray*} 
& & {\rm rank}_{R} {\rm HH}^{n, s}({\rm N}_m(R), {\rm B}_m(R)) \\ 
& = & \left\{
\begin{array}{ll}
\varphi(s) &  (n=s), \\
\displaystyle (-1)^m\varphi(s)+\sum_{k=0}^{m-1} (-1)^k (k+1) \varphi(s+m-k-1) & (n=s+m-1), \\
0 & (\mbox{\rm otherwise}).  
\end{array}
\right.
\end{eqnarray*} 
\end{theorem}

%%%%%%%%%%%%%%%%  the beginning of \if0 
\if0 

\subsubsection{Double check}

We consider the cochain complex
\[ (E_1^{*,-(m-1)+q}({\rm B}),d_1),\]
where
\[ E_1^{(m-1)-r,-(m-1)+q}({\rm B})\cong
   {\rm N}^!_{q-r}\otimes_R {\rm Gr}^{(m-1)-r}({\rm B}). \]
Since
\[ {\rm rank}_RE_1^{(m-1)-r,-(m-1)+q}({\rm B})=
    (r+1)\varphi(q-r)\]
for $0\le r\le m-1$,
the Euler characteristic
of $(E_1^{*,-(m-1)+q}({\rm B}),d_1)$ is given by
\[ \chi(E_1^{*,-(m-1)+q}({\rm B}))=
   (-1)^{m-1}\sum_{r=0}^{m-1} (-1)^r (r+1)\varphi(q-r).\]

\begin{remark}\rm
For $0\le q<m-1$,
we have
\[ \chi(E_1^{*,-(m-1)+q}({\rm B}))=
   (-1)^{m-1}\sum_{r=0}^q(-1)^r(r+1)\varphi(q-r)=
   \chi(E_2^{*,-(m-1)+q}({\rm B})).\]
\end{remark}

\begin{remark}\rm
We have
\[ \begin{array}{rcl}
    \chi(E_2^{*,0}({\rm B}))&=&
    {\rm rank}_RE_2^{0,0}({\rm B})+(-1)^{m-1}{\rm rank}_RE_2^{m-1,0}({\rm B})\\[2mm]
    &=&\displaystyle
    (-1)^{m-1}\sum_{k=0}^{m-1}(-1)^k(k+1)\varphi(m-1-k)\\[4mm]
    &=&\chi(E_1^{*,0}({\rm B})).\\
   \end{array}\]
\end{remark}

\begin{remark}\rm
For $q\ge m$,
we have
\[ \begin{array}{cl}
     &\chi(E_2^{*,-(m-1)+q}({\rm B}))\\[2mm]
    =&{\rm rank}_RE_2^{0,-(m-1)+q}({\rm B})
      +(-1)^{m-1}{\rm rank}_RE_2^{m-1,-(m-1)+q}({\rm B})\\[2mm]
    =&\displaystyle\varphi(q-(m-1))
     +(-1)^{m-1}\sum_{k=0}^{m-1}(-1)^k(k+1)\varphi(q-k)
      -\varphi(q-m+1)\\[4mm]
    =&\displaystyle(-1)^{m-1}\sum_{k=0}^{m-1}(-1)^k(k+1)\varphi(q-k)\\[4mm]
    =&\chi(E_1^{*,-(m-1)+q}({\rm B})).\\
   \end{array}\]
\end{remark}

\fi
%%%%%%%%%%%%%%%%%%   the end of \if0

\subsection{Freeness of $E_2^{1,q}({\rm N})$}\label{subsection:freenessE21qN} 
We have shown that $E_{2}^{p, q}({\rm N}) = 0$ unless $p=0, 1, m-1$ (Remark~\ref{remark:E2pq=0unlessp=01m-1}) 
and that $E_{2}^{m-1, q}({\rm N})$ 
is a finitely generated free module over $R$ by Corollary~\ref{cor:E2pqB=0}. 
In this subsection, 
we show that $E_2^{1,q}({\rm N})$ is a free $R$-module. We also show that 
$E_2^{0, 0}({\rm N})\cong R$ and $E_2^{0, q}({\rm N})=0$ for $q \neq 0$. 

\begin{proposition}\label{prop:E200A-E200Bisom} 
For $q = 0$, $E_2^{0, 0}({\rm N})\cong R$ and $E_2^{0, 0}({\rm N}) \to E_2^{0, 0}({\rm B})$ is an isomorphism. 
\end{proposition}

\proof
Let us consider $d_{1}^{0, 0} : E_1^{0, 0}({\rm N}) \cong  {\rm N}^{!}_{0}\otimes_{R} {\rm Gr}^{0}({\rm N}) \to E_1^{1, 0}({\rm N}) \cong {\rm N}^{!}_{1}\otimes_{R} {\rm Gr}^{1}({\rm N})$. For $c I_m \in {\rm N}^{!}_{0}\otimes_{R} {\rm Gr}^{0}({\rm N}) = RI_m$ with $c \in R$,  
\[
d_{1}^{0, 0}(c I_m) = \sum_{i=1}^{m-1} (y_ic-cy_i)\otimes E_{i, i+1} = 0. 
\]
Hence, $E_{2}^{0, 0}({\rm N}) \cong {\rm Ker} d_{1}^{0, 0} = RI_m \cong R$. Since $E_{2}^{0, 0}({\rm B}) =RI_m$ by 
Proposition~\ref{prop:E20qBgenerated}, $E_2^{0, 0}({\rm N}) \to E_2^{0, 0}({\rm B})$ is an isomorphism.  
\qed 

\begin{proposition}\label{prop:E20qA=0}
For $q \neq 0$, $E_2^{0, q}({\rm N})=0$. 
\end{proposition}

\proof
Obviously, $E_{2}^{0, q}({\rm N}) \cong {\rm Ker} d_{1}^{0, q} = E_{1}^{0, q}({\rm N}) = 0$ if $q<0$. 
We only need to consider the case that $q>0$. 
Let us consider $d_{1}^{0, q} : E_1^{0, q}({\rm N}) \cong  {\rm N}^{!}_{q}\otimes_{R} {\rm Gr}^{0}({\rm N}) \to E_1^{1, q}({\rm N}) \cong {\rm N}^{!}_{q+1}\otimes_{R} {\rm Gr}^{1}({\rm N})$ for $q>0$. 
Let $z \in {\rm Ker} d_{1}^{0, q}$. 
Note that $d_{1}^{0, q}$ can be regarded as a restriction of $d_{1}^{0, q} : E_1^{0, q}({\rm B}) \to E_1^{1, q}({\rm B}) \cong E_1^{1, q}({\rm N})$. 
As in the proof of Proposition~\ref{prop:E20qBgenerated}, there exists $g_1, g_2, \ldots, g_{m-1} \in {\rm N}^{!}_{q-1}$ 
such that 
\begin{multline*}
z = y_1g_1\otimes E_{1,1} + \{ y_2g_2 +(-1)^{q} g_1y_1\}\otimes E_{2,2} + \cdots \\ 
 + \{ y_{m-1}g_{m-1} +(-1)^{q} g_{m-2}y_{m-2} \}\otimes E_{m-1,m-1} 
+(-1)^{q} g_{m-1}y_{m-1}\otimes E_{m,m}. 
\end{multline*} 
Using $z \in E_1^{0, q}({\rm N}) \cong  {\rm N}^{!}_{q}\otimes_{R} RI_m$, we have 
\begin{eqnarray}
y_1g_1 =  y_{2}g_{2} +(-1)^{q} g_{1}y_{1} = 
\cdots = (-1)^{q} g_{m-1}y_{m-1}. \label{eq:condKerd1pqN}
\end{eqnarray}
%When $q=1$, the condition (\ref{eq:condKerd1pqN}) implies that $z=0$. 
The left hand side and the right hand side of  (\ref{eq:condKerd1pqN}) are contained in $y_1{\rm N}^{!}_{q-1}$ and 
${\rm N}^{!}_{q-1}y_{m-1}$, respectively.  Hence, the both sides of (\ref{eq:condKerd1pqN}) are contained 
in $y_1{\rm N}^{!}_{q-2}y_{m-1}$, while $y_{2}g_{2} +(-1)^{q} g_{1}y_{1} \in y_2{\rm N}^{!}_{q-1} + {\rm N}^{!}_{q-1}y_1$.  
This implies that $y_1g_1 = \cdots = (-1)^{q+1} g_{m-1}y_{m-1}=0$ and that $z=0$.   
Thus, we have $E_{2}^{0, q}({\rm N}) = {\rm Ker} d_1^{0, q} = 0$ for $q > 0$.  
\qed 

\begin{theorem}\label{th:E21qNisfree}
For $q \in {\mathbb Z}$, $E_2^{1,q}({\rm N})$ is a finitely generated free module over $R$.  
%for any commutative ring $R$ and any integer $q$.
\end{theorem}

\proof
We have an exact sequence of cochain complexes
\[ 0 \longrightarrow E_1^{*,q}({\rm N};R) \longrightarrow 
   E_1^{*,q}({\rm B};R) \longrightarrow E_1^{*,q}({\rm B}/{\rm N};R) \longrightarrow 0. \]
This induces a long exact sequence
\begin{eqnarray} 
\cdots\longrightarrow E_2^{*,q}({\rm N};R)\longrightarrow E_2^{*,q}({\rm B};R)
   \longrightarrow E_2^{*,q}({\rm B}/{\rm N};R)\longrightarrow E_2^{*+1,q}({\rm N};R) \longrightarrow \cdots.  \label{eq:longexactseqE2} 
\end{eqnarray}

Let $q=0$. 
The map $E_2^{0,0}({\rm N};R)\to E_2^{0,0}({\rm B}; R)$ is an isomorphism by Proposition~\ref{prop:E200A-E200Bisom}. 
Since $E_2^{1,0}({\rm B};R)=0$ by Proposition~\ref{prop:E2pqB=0}, 
we obtain that $E_2^{1,0}({\rm N};R)\cong E_2^{0,0}({\rm B}/{\rm N};R)$,  
which is isomorphic to 
the finitely generated free $R$-module ${\rm B}/{\rm N}$ by Theorem~\ref{th:HHB/N}.   

Let $q\neq 0$. 
Since $E_2^{0,q}({\rm N};R)=0$ and $E_2^{1,q}({\rm B};R)=0$ by Propositions~\ref{prop:E20qA=0} and \ref{prop:E2pqB=0},
we obtain an exact sequence
\begin{equation}\label{equation:exac-sequence-E2-A-B-B/A}
   0 \longrightarrow E_2^{0,q}({\rm B};R) \longrightarrow 
   E_2^{0,q}({\rm B}/{\rm N};R) \longrightarrow E_2^{1,q}({\rm N};R) \longrightarrow 0.
\end{equation}
By the universal coefficient theorem,
we have an exact sequence
\[ 0 \longrightarrow E_2^{1,q}({\rm N};\mathbb{Z})\otimes R
   \longrightarrow E_2^{1,q}({\rm N};R)\longrightarrow
   {\rm Tor}^{\mathbb{Z}}_1(E_2^{2,q}({\rm N}; {\Bbb Z}),R) \longrightarrow 0.\] 
Since $E_2^{2,q}({\rm N};{\Bbb Z})$ is a free ${\Bbb Z}$-module by Corollary~\ref{cor:E2pqB=0},  
we obtain an isomorphism
\begin{equation}\label{equation:e2-1-q-A-Z-R}
   E_2^{1,q}({\rm N};\mathbb{Z})\otimes R
   \stackrel{\cong}{\longrightarrow}
   E_2^{1,q}({\rm N};R). 
\end{equation}

Let $k$ be a field.
Note that
$E_2^{0,q}({\rm B};\mathbb{Z})$ and
$E_2^{0,q}({\rm B}/{\rm N};\mathbb{Z})$
are finitely generated free $\mathbb{Z}$-modules by Proposition~\ref{prop:E2m-1qfree} and Theorem~\ref{th:HHB/N}, 
and hence that
$\dim_kE_2^{0,q}({\rm B};k)=
{\rm rank}_{\mathbb{Z}}E_2^{0,q}({\rm B};\mathbb{Z})$ 
and 
$\dim_kE_2^{0,q}({\rm B}/{\rm N};k)=
{\rm rank}_{\mathbb{Z}}E_2^{0,q}({\rm B}/{\rm N};\mathbb{Z})$ by $E_{2}^{1, q}({\rm B}; {\Bbb Z}) = 
E_{2}^{1, q}({\rm B}/{\rm N}; {\Bbb Z}) = 0$ (Proposition~\ref{prop:E2pqB=0} and Theorem~\ref{th:HHB/N})  and 
the universal coefficient theorem.  
By (\ref{equation:exac-sequence-E2-A-B-B/A})
and
(\ref{equation:e2-1-q-A-Z-R}),
we obtain 
\[ \begin{array}{rcl}
    \dim_k E_2^{1,q}({\rm N};\mathbb{Z})\otimes k&=&
    \dim_k E_2^{1,q}({\rm N};k)\\[2mm]
    &=&\dim_kE_2^{0,q}({\rm B}/{\rm N};k)-\dim_kE_2^{0,q}({\rm B};k)\\[2mm]
    &=&{\rm rank}_{\mathbb{Z}}E_2^{0,q}({\rm B}/{\rm N};\mathbb{Z})
      -{\rm rank}_{\mathbb{Z}}E_2^{0,q}({\rm B};\mathbb{Z}).\\[2mm] 
   \end{array}\]
This shows that
$\dim_kE_2^{1,q}({\rm N}; {\Bbb Z})\otimes k$ is independent
from the field $k$. 
Since each $E_{1}^{p, q}({\rm N};\mathbb{Z})$ is a finitely generated free ${\Bbb Z}$-module, 
$E_2^{1,q}({\rm N};\mathbb{Z})$ is finitely generated ${\Bbb Z}$-module. 
Thus,
we see that
$E_2^{1,q}({\rm N};\mathbb{Z})$
is a finitely generated free $\mathbb{Z}$-module by the fundamental theorem of finitely generated abelian groups.
Hence 
$E_2^{1,q}({\rm N};R)\cong E_2^{1,q}({\rm N};\mathbb{Z})\otimes R$
is a finitely generated free $R$-module.
\qed

\subsection{The rank of $E_2^{p,q}({\rm N})$}\label{subsection:rankE2pqN}
Note that  
${\rm rank}_RE_2^{p,q}({\rm N})=0$
unless $p=0,1$ or $m-1$.
In this subsection, 
we calculate the rank of $E_2^{p,q}({\rm N})$ for $p=0, 1, m-1$, 
which is a finitely generated free $R$-module. 
As a result, we can determine the $R$-module structure of ${\rm HH}^{n}({\rm N}, {\rm N})$.  

%For $q\in\mathbb{Z}$,
%we recall that  
%\[ \varphi(q)=\,{\rm rank}_R {\rm N}^!_q .\]

\begin{theorem}\label{th:E2pqA} 
We have
\[ \begin{array}{rcl}
   {\rm rank}_RE_2^{0,q}({\rm N})&=&
    \left\{\begin{array}{ll}
            1 & (q=0),\\[2mm]
            0 & (q\neq 0), \\
           \end{array}\right.\\[7mm]
    {\rm rank}_RE_2^{1,q}({\rm N})&=&
     \left\{\begin{array}{ll}
             m-1 & (q=0), \\[2mm]
             (m-2)\varphi(q) & (q\neq 0),\\
            \end{array}\right.\\[7mm]
    {\rm rank}_RE_2^{m-1,q}({\rm N})&=&\displaystyle
           (-1)^m\varphi(q)+\sum_{k=0}^{m-1}(-1)^k(k+1)
           \varphi(q+m-k-1).\\
   \end{array} \]
\end{theorem}

%\newpage

\proof
When $p=0$,
recall that 
$E_2^{0,0}({\rm N})\cong R$ and
$E_2^{0,q}({\rm N})=0$ for $q\neq 0$ by Propositions~\ref{prop:E200A-E200Bisom} and \ref{prop:E20qA=0}. 
When $p=m-1$,
since $E_2^{m-1,q}({\rm N})\cong E_2^{m-1,q}({\rm B})$ by Corollary~\ref{cor:E2pqB=0}, 
we have
${\rm rank}_RE_2^{m-1,q}({\rm N})={\rm rank}_RE_2^{m-1,q}({\rm B})$, which can be calculated by  
Theorem~\ref{th:rankE2pqB}.  

We consider the case when $p=1$. 
Recall the proof of Theorem~\ref{th:E21qNisfree}.  
We have an exact sequence
\[ 0 \longrightarrow E_2^{0,q}({\rm N}) \longrightarrow E_2^{0,q}({\rm B}) \longrightarrow 
        E_2^{0,q}({\rm B}/{\rm N}) \longrightarrow E_2^{1,q}({\rm N}) \longrightarrow 0. \]
When $q=0$, 
%the map 
%$E_2^{0,0}({\rm N})\to E_2^{0,0}({\rm B})$ is an isomorphism.
%Thus,
%$E_2^{1,0}({\rm N})\cong E_2^{0,0}({\rm B}/{\rm N})$ and
%\[ {\rm rank}_RE_2^{1,0}({\rm N})={\rm rank}_RE_2^{0,0}({\rm B}/{\rm N})
%   = m-1.\]
we have seen that $E_2^{1,0}({\rm N};R)\cong E_2^{0,0}({\rm B}/{\rm N};R) \cong {\rm B}/{\rm N}$.    Hence, we obtain 
\[ {\rm rank}_RE_2^{1,0}({\rm N})  = m-1.\]
When $q\neq 0$, since
$E_2^{0,q}({\rm N})=0$ by Proposition~\ref{prop:E20qA=0}, 
we have
\[ \begin{array}{rcl}
    {\rm rank}_RE_2^{1,q}({\rm N})&=&
    {\rm rank}_RE_2^{0,q}({\rm B}/{\rm N})-{\rm rank}_RE_2^{0,q}({\rm B})\\[2mm]
    &=& (m-1)\varphi(q)-\varphi(q)\\[2mm]
    &=& (m-2)\varphi(q).\\
   \end{array}\]
Here we used Theorems~\ref{th:HHB/N} and \ref{th:rankE2pqB}. 
\qed

\bigskip 

Summarizing the discussions above, we have the following theorems by Proposition~\ref{prop:HHNNcollapses} and 
Theorem~\ref{th:E2pqA}.

\begin{theorem}\label{th:mainthm} 
Let $m \ge 3$. The Hochschild cohomology ${\rm HH}^{n}({\rm N}_m(R), {\rm N}_m(R))$ is a free $R$-module for $n\ge 0$. The rank of ${\rm HH}^{n}({\rm N}_m(R), {\rm N}_m(R))$ is given by 
\begin{eqnarray*} 
& & {\rm rank}_{R} {\rm HH}^{n}({\rm N}_m(R), {\rm N}_m(R)) \\ 
& = & \left\{
\begin{array}{lc}
2 &  (n=0), \\
2m-4 & (n=1), \\
\displaystyle \varphi(n)+(m-4)\varphi(n-1) +(-1)^m\varphi(n-m+1) +\sum_{k=2}^{m-1} (-1)^k (k+1)\varphi(n-k) & (n\ge 2).  
\end{array}
\right.
\end{eqnarray*} 
\end{theorem}

\begin{theorem}\label{th:mainthmHHns} 
Let $m \ge 3$. For $n \ge 0$ and $s \in {\Bbb Z}$, ${\rm HH}^{n, s}({\rm N}_m(R), {\rm N}_m(R))$ is a free $R$-module. 
Furthermore, 
\[ {\rm HH}^{n,s}({\rm N}_m(R), {\rm N}_m(R))\cong E_2^{n-s,s}({\rm N}_m(R))\] 
as $R$-modules and 
\begin{eqnarray*} 
& & {\rm rank}_{R} {\rm HH}^{n, s}({\rm N}_m(R), {\rm N}_m(R)) \\ 
& = & \left\{
\begin{array}{ll}
1 &  (n=0, s=0), \\
m-1 & (n=1, s=0), \\
(m-2)\varphi(s) & (n=s+1, s\neq 0), \\ 
\displaystyle (-1)^m\varphi(s)+\sum_{k=0}^{m-1} (-1)^k (k+1) \varphi(s+m-k-1) & (n=s+m-1), \\
0 & (\mbox{\rm otherwise}).    
\end{array}
\right.
\end{eqnarray*} 
\end{theorem}

%\subsection{Construction of an $R$-basis of ${\rm HH}^{\ast}({\rm N}_m(R), {\rm M}_m(R)/{\rm N}_m(R))$} 
%In this subsection, we construct an $R$-basis of ${\rm HH}^{\ast}({\rm N}_m(R), {\rm M}_m(R)/{\rm N}_m(R))$. 

%
%\begin{remark}\rm
%Note that we have an isomorphism
%\[ {\rm HH}^{n,s}({\rm N}, {\rm N})\cong E_2^{n-s,s}({\rm N}).\]
%\end{remark}
%}

%\input{N3product_structure}
\section{Product structure on ${\rm HH}^*({\rm N}_m(R),{\rm N}_m(R))$}\label{section:ProductHHNN} 
In this section, we describe the product structure on ${\rm HH}^*({\rm N}_m(R),{\rm N}_m(R))$ for $m \ge 3$. 
In \S\ref{subsection:productm=3}, we deal with the case $m=3$ explicitly, which is different from the case 
$m \ge 4$.  
In \S\ref{subsection:productm>=4}, we deal with the case $m \ge 4$ in general.  
In any case, there exists an augmentation map $\epsilon : {\rm HH}^{\ast}({\rm N}_m(R), {\rm N}_m(R)) \to R$ as an $R$-algebra homomorphism such that the Kernel $\overline{{\rm HH}^{\ast}}({\rm N}_m(R), {\rm N}_m(R))$ of $\epsilon$ satisfies 
\[
\overline{{\rm HH}^{\ast}}({\rm N}_m(R), {\rm N}_m(R))  \cdot \overline{{\rm HH}^{\ast}}({\rm N}_m(R), {\rm N}_m(R)) =0.
\] 
In particular, we see that ${\rm HH}^*({\rm N}_m(R),{\rm N}_m(R))$ is an infinitely generated algebra over $R$. 

%\bigskip 

\subsection{The case $m=3$}\label{subsection:productm=3} 
In this subsection, we set $m=3$ and ${\rm N} = {\rm N}_3(R)$. 
Recall ${\rm N}^{!} = R\langle y_1, y_2 \rangle/\langle y_1y_2 \rangle$ in \S\ref{subsection:HHNmR}.  
We define $c(i, j) \in {\rm N}^{!}$ by 
\[
c(i, j) = y_2^i y_1^j \in {\rm N}^{!}
\]
for $i, j \ge 0$. (Set $c(0, 0)=1$.) Then we can describe the homogeneous part ${\rm N}^{!}_{n}$ of ${\rm N}^{!}$ of degree $n$ by 
\[
{\rm N}^{!}_{n} = R\{ c(i, j) \mid i, j \ge 0, i+j = n \}
\]
for $n \ge 0$. Note that $\varphi(n) = {\rm rank}_{R} {\rm N}^{!}_{n} = n+1$.

Let us consider the spectral sequence
\[ E_1^{p,q}={\rm HH}^{p+q}({\rm N},{\rm Gr}^p({\rm N}))
   \Longrightarrow {\rm HH}^{p+q}({\rm N}, {\rm N}).\]
By the discussions in \S\ref{section:HH-N-N},    
we have $E_2^{p,q}\cong E_{\infty}^{p,q}$ and
\[ \begin{array}{rcl}
     E_{\infty}^{0,q}&\cong&
     \left\{\begin{array}{cl}
            R\{c(0,0)\otimes I_3\}& (q=0),\\[2mm]
            0                     & (q\neq 0),\\  
            \end{array}\right.\\[5mm]
     E_{\infty}^{1,q}&\cong&
     \left\{\begin{array}{cl}
            R\{c(i,q-i+1)\otimes E_{1,2},
               c(q+1,0)\otimes E_{2,3}|\ 0\le i< q\}& (q\ge 1),\\[2mm]
            R\{c(0,1)\otimes E_{1,2},c(1,0)\otimes E_{2,3}\}
                                       & (q=0),\\[2mm]  
            0                          & (q<0),\\[2mm]
            \end{array}\right.\\[9mm]
     E_{\infty}^{2,q}&\cong&
     \left\{\begin{array}{cl}
            R\{c(i,j)\otimes E_{1,3}|\ i+j=q+2,i>0,j>0\}& (q\ge 0),\\[2mm]
            R\{c(0,0)\otimes E_{1,3}\}& (q=-2),\\[2mm] 
            0                     & (\mbox{\rm otherwise}). \\  
            \end{array}\right.\\[2mm]

   \end{array}\] 

By direct inspection,
we obtain the following lemma.

\begin{lemma}\label{lemma:product-E-infinity-term}
The element $c(0,0)\otimes I_3\in E_{\infty}^{0,0}$
is a unit of the bigraded algebra $E_{\infty}^{*,*}$.
For any $r,s\ge 1$,
the product map 
$E_{\infty}^{r,q}\otimes_R E_{\infty}^{s,q'}\to E_{\infty}^{r+s,q+q'}$
is a zero map.
\end{lemma}

We have
\[ F^r{\rm HH}^*({\rm N}, {\rm N})=\ \mbox{\rm Im}
   ({\rm HH}^*({\rm N},F^r{\rm N})\longrightarrow {\rm HH}^*({\rm N}, {\rm N})) \]
and
\[ E_{\infty}^{p,q}\cong F^p{\rm HH}^{p+q}({\rm N}, {\rm N})/
                       F^{p+1}{\rm HH}^{p+q}({\rm N}, {\rm N}).\]

The map
\[ {\rm HH}^0({\rm N}, {\rm N})=F^0{\rm HH}^0({\rm N}, {\rm N})
   \longrightarrow F^0{\rm HH}^0({\rm N}, {\rm N})/
   F^1{\rm HH}^0({\rm N}, {\rm N})\cong E_{\infty}^{0,0} \]
gives an augmentation map
$\epsilon: {\rm HH}^*({\rm N}, {\rm N})\to R$.
Note that
\[ {\rm HH}^0({\rm N}, {\rm N})=R\{c(0,0)\otimes I_3,c(0,0)\otimes E_{1,3}\} \]
and 
\[ \begin{array}{rcl}
    \epsilon(c(0,0)\otimes I_3)&=&1, \\[2mm] 
    \epsilon(c(0,0)\otimes E_{1,3})&=&0. \\
   \end{array}\]

We can identify $F^1{\rm HH}^*({\rm N}, {\rm N})$  
with the kernel of $\epsilon$.
By Lemma~\ref{lemma:product-E-infinity-term},
the product map
\[ F^1{\rm HH}^*({\rm N}, {\rm N})\otimes_R F^1{\rm HH}^*({\rm N}, {\rm N})
   \longrightarrow F^2{\rm HH}^*({\rm N}, {\rm N}) \]
is trivial.
Hence we obtain the following theorem.

\begin{theorem}\label{th:productzerom=3} 
There is an augmentation map $\epsilon:{\rm HH}^*({\rm N}, {\rm N})\to R$ 
such that 
$\epsilon(c(0,0)\otimes I_3)=1$ 
and
$\epsilon(c(0,0)\otimes E_{1,3})=0$.
Let $\overline{\rm HH}^*({\rm N}, {\rm N})$ be the kernel
of $\epsilon$.
Then we have 
\[ \overline{\rm HH}^*({\rm N}, {\rm N})\cdot \overline{\rm HH}^*({\rm N}, {\rm N})=0. \]
\end{theorem}

\subsection{The case $m\ge4$}\label{subsection:productm>=4} 

%First, we recall that a spectral sequence can be
%formulated in an abelian category. 
%Let $\mathcal{A}$ be an abelian symmetric monoidal  category
%in which the tensor product
%$\otimes: \mathcal{A}\times\mathcal{A}\to\mathcal{A}$
%is right exact separately in each variable.
%
%\begin{definition}\rm
%Let $(A^*,d)$ be a differential graded algebra in $\mathcal{A}$.
%Suppose that we have a filtration
%\[ A^*=F^0A^*\supset F^1A^*\supset\cdots\supset
%       F^nA^*\supset\cdots.\]
%A triple $(A^*,d,\{F^rA^*\}_{r\ge 0})$ is said to be
%a filtered differential graded algebra
%if it satisfies the following two conditions:
%\begin{enumerate}
%\item[(1)]
%For any $n\ge 0$, $d(F^nA^*)\subset F^nA^*$.
%\item[(2)]
%For any $r,s\ge 0$, $F^rA^*\cdot F^sA^*\subset F^{r+s}A^*$.
%\end{enumerate}
%\end{definition}
%
%By \cite[Theorem~2.14]{McCleary},
%there is a spectral sequence 
%\[ E_1^{p,q}=H^{p+q}(F^pA/F^{p+1}A)
%   \Longrightarrow H^{p+q}(A) \]
%of algebras in $\mathcal{A}$,
%which converges to $H^{p+q}(A)$
%as an algebra.

Let ${\rm N}={\rm N}_m(R)$ for $m\ge 4$.
Recall that we have a decomposition
\[ C^*({\rm N}, {\rm N})=\bigoplus_{s\in\mathbb{Z}}C^{*,s}({\rm N}, {\rm N}) \]
which is compatible with the filtration.
We regard 
\[ C^p({\rm N}, {\rm N})=\bigoplus_{s\in\mathbb{Z}}C^{p,s}({\rm N}, {\rm N}) \]
as a $\mathbb{Z}$-graded $R$-module.
Then the triple 
$(C^*({\rm N}, {\rm N}),d,\{F^rC^*({\rm N}, {\rm N})\}_{r\ge 0})$
is a filtered differential graded algebra
in the category of $\mathbb{Z}$-graded $R$-modules.
Thus, we obtain a multiplicative spectral sequence 
\[ E_1^{p,q}({\rm N})\Longrightarrow {\rm HH}^{p+q}({\rm N}, {\rm N}) \]
in the abelian category of $\mathbb{Z}$-graded $R$-modules (for details, see \S\ref{subsection:reviewss} and \S\ref{subsection:gradingonss}).  

\begin{lemma}\label{lemma:productzerom>=4}
Let $m\ge 4$.  
If $a\in {\rm HH}^{1+q,q}({\rm N}, {\rm N})$ and
$b\in {\rm HH}^{1+q',q'}({\rm N}, {\rm N})$, then $ab=0$ in ${\rm HH}^{2+q+q',q+q'}({\rm N}, {\rm N})$.
\end{lemma}

\proof
By Theorem~\ref{th:mainthmHHns}, we may assume that $q, q' \ge 0$. 
Let $x\in E_{\infty}^{1,q,q}$ and $y\in E_{\infty}^{1,q',q'}$ 
be elements which represent $a$ and $b$, respectively. 
Since $E_{2}^{2,q+q',q+q'}({\rm N}) =0$ for $m \ge 4$ by Theorem~\ref{th:E2pqA},  
$E_{\infty}^{2,q+q',q+q'}({\rm N}) =0$. Hence $xy=0$, which implies that 
$ab$ is represented by an element in $E_{\infty}^{m-1,q+q'-m+3,q+q'}({\rm N})$. 
By Lemma~\ref{lemma:graded_ss_vanishing}, if $m \ge 4$, then 
$E_{\infty}^{m-1,q+q'-m+3,q+q'}({\rm N})=E_{1}^{m-1,q+q'-m+3,q+q'}({\rm N})=0$. 
Therefore $ab=0$. 
%We have
%$xy\in E_{\infty}^{2,q+q',q+q'}(A)$.
%Since $m\ge 4$,
%we have $E_{\infty}^{2,q+q',q+q'}(A)=0$.
%Thus, $xy=0$.
%This implies that $ab$
%is represented by an element
%in $E_{\infty}^{m-1,q+q'-m+3,q+q'}(A)$.
%Since $m\ge 4$,
%$E_{\infty}^{m-1,q+q'-m+3,q+q'}(A)=0$.
%Hence $ab=0$.
\qed

\bigskip 

Recall $I_m \in C^{0, 0}({\rm N}, {\rm N})$ %(in \S~\ref{subsection:sssubquotient}) 
is a generator of ${\rm HH}^{0, 0}({\rm N}, {\rm N})$ ({\it cf}. Proposition~\ref{prop:E200A-E200Bisom} and Theorem~\ref{th:mainthmHHns}).  
%%%%%%%%%%%%%%%%%%%%%%%%%%%%%%   the beginning of \if0 
\if0 
In the same way of the case $m=3$, 
the map
\[ {\rm HH}^0({\rm N}, {\rm N})=F^0{\rm HH}^0({\rm N}, {\rm N})
   \longrightarrow F^0{\rm HH}^0({\rm N}, {\rm N})/
   F^1{\rm HH}^0({\rm N}, {\rm N})\cong E_{\infty}^{0,0} \]
gives an augmentation map
$\epsilon: {\rm HH}^*({\rm N}, {\rm N})\to R$.
Note that
\[ {\rm HH}^0({\rm N}, {\rm N})={\rm HH}^{0, 0}({\rm N}, {\rm N})\oplus {\rm HH}^{m-1, -(m-1)}({\rm N}, {\rm N}) \]
and 
\[ \begin{array}{rcl}
    \epsilon(I_m)&=&1\\[2mm] 
    \epsilon( {\rm HH}^{m-1, -(m-1)}({\rm N}, {\rm N}))&=&0. \\
   \end{array}\]
We can identify $F^1{\rm HH}^*({\rm N}, {\rm N})$  
with the kernel of $\epsilon$. 
By Lemma~\ref{lemma:productzerom>=4}, we have the following theorem. 
\fi
%%%%%%%%%%%%%%%%%%%%%%%%%%%%%  the end of \if0 
By the decomposition ${\rm HH}^{\ast}({\rm N}, {\rm N}) = 
\bigoplus_{n\ge 0, s\in {\Bbb Z}} {\rm HH}^{n, s}({\rm N}, {\rm N})$,  
we have an augmentation map $\epsilon: {\rm HH}^*({\rm N}, {\rm N})\to R$ as an $R$-algebra homomorphism 
such that $\epsilon(I_m) = 1$ and $\epsilon({\rm HH}^{n, s}({\rm N}, {\rm N})) = 0$ for 
$(n, s)\neq (0, 0)$. We can identify $F^1{\rm HH}^*({\rm N}, {\rm N})$  
with the kernel of $\epsilon$. Using Lemma~\ref{lemma:productzerom>=4}, we see that 
${\rm HH}^{n, s}({\rm N}, {\rm N})\cdot {\rm HH}^{n', s'}({\rm N}, {\rm N}) = 0$ if 
$(n, s)\neq (0, 0)$ and $(n', s')\neq (0, 0)$. Hence, we have the following theorem. 

\begin{theorem}\label{th:productzerom>=4}
Let $m \ge 4$. 
There is an augmentation map $\epsilon:{\rm HH}^*({\rm N}, {\rm N})\to R$ 
such that 
$\epsilon(I_m)=1$ for $I_m \in {\rm HH}^{0, 0}({\rm N}, {\rm N})$  
and
$\epsilon({\rm HH}^{m-1, -(m-1)}({\rm N}, {\rm N}))=0$. 
Let $\overline{\rm HH}^*({\rm N}, {\rm N})$ be the kernel
of $\epsilon$.
Then we have 
\[ \overline{\rm HH}^*({\rm N}, {\rm N})\cdot \overline{\rm HH}^*({\rm N}, {\rm N})=0. \]
\end{theorem} 

\begin{corollary}\label{cor:infinitelygenerated} 
Let $m \ge 3$. The Hochschild cohomology algebra ${\rm HH}^{\ast}({\rm N}, {\rm N})$ is an infinitely generated algebra over $R$. 
\end{corollary} 

\proof
Suppose that there exists a finite set $G = \{ x_i \mid 1 \le i \le l \}$ 
of generators of ${\rm HH}^{\ast}({\rm N}, {\rm N})$ 
as an $R$-algebra. 
We may assume that $x_i$ is contained in ${\rm HH}^{n_i, s_i}({\rm N}, {\rm N})$ for each $i$. 
By Theorems~\ref{th:productzerom=3} and \ref{th:productzerom>=4}, $x_ix_j=0$ if $(n_i, s_i)\neq (0, 0)$ and 
$(n_j, s_j)\neq (0, 0)$. However, ${\rm rank}_{R} {\rm HH}^{s+1, s}({\rm N}, {\rm N}) = (m-2)\varphi(s) >0$ for $s> 0$ by Theorem~\ref{th:mainthmHHns}.   This implies that $G$ can not generate ${\rm HH}^{\ast}({\rm N}, {\rm N})$, 
which is a contradiction. Hence,  ${\rm HH}^{\ast}({\rm N}, {\rm N})$ is an infinitely generated algebra over $R$. 
\qed

\begin{remark}\rm 
By \cite[Theorem~7.3]{GSS}, if $\Lambda = KQ/I$ is an indecomposable monomial algebra over a field $K$, then ${\rm HH}^{\ast}(\Lambda)/{\mathcal N}$ is a commutative finitely generated $K$-algebra of Krull dimension at most one, where ${\mathcal N}$ is the ideal of ${\rm HH}^{\ast}(\Lambda)$ generated by the homogeneous nilpotent elements. In the ${\rm N}_m(K)$ case for $m \ge 3$, ${\mathcal N} = \overline{\rm HH}^*({\rm N}_m(K), {\rm N}_m(K)) = {\rm Ker} \epsilon$, and ${\rm HH}^{\ast}(\Lambda)/{\mathcal N} \cong K$ has Krull dimension zero.  
\end{remark}

\section{Gerstenhaber bracket on ${\rm HH}^*({\rm N}_m(R), {\rm N}_m(R))$}
In this section, we describe the Gerstenhaber bracket on ${\rm HH}^*({\rm N}_m(R), {\rm N}_m(R))$. 

\subsection{Cocycle representatives}
%Let $R$ be a commutative ring.
Set ${\rm N}={\rm N}_m(R)$ for $m \ge 3$. 

\begin{proposition}\label{prop:E21qgenerator}
For $q\ge 0$, 
\[ 
E_2^{1,q}({\rm N}, {\rm N}) = R\{ y_{i-1}y_I\otimes E_{i-1, i}-(-1)^qy_Iy_i\otimes E_{i, i+1}  \mid 1 \le i  \le m, |I|=q \} 
\]
as $R$-subquotients of $E_{1}^{1, q}({\rm N}, {\rm N}) \cong {\rm N}^{!}_{q+1}\otimes_{R} {\rm Gr}^{1}({\rm N})$. 
Here, $y_{i-1}y_I\otimes E_{i-1, i}-(-1)^qy_Iy_i\otimes E_{i, i+1}$ is regarded as 
$-(-1)^qy_Iy_{1}\otimes E_{1, 2}$ if $i=1$ and $y_{m-1}y_{I}\otimes E_{m-1, m}$ if $i=m$, respectively. 
\end{proposition}

\proof
Let ${\rm M} = {\rm M}_m(R)$. 
We have an exact sequence of cochain complexes
\[ 0 \longrightarrow E_1^{*,q}({\rm N}, {\rm N}) \longrightarrow 
   E_1^{*,q}({\rm N}, {\rm M}) \longrightarrow E_1^{*,q}({\rm N}, {\rm M}/{\rm N}) \longrightarrow 0. \]
This induces a long exact sequence
\begin{eqnarray} 
\cdots \longrightarrow E_2^{*,q}({\rm N}, {\rm N})\longrightarrow E_2^{*,q}({\rm N}, {\rm M})
   \longrightarrow E_2^{*,q}({\rm N}, {\rm M}/{\rm N}) \longrightarrow E_2^{*+1,q}({\rm N}, {\rm N}) \longrightarrow \cdots.  \label{eq:longexactseqE2MN} 
\end{eqnarray}
Since 
\[
\begin{array}{ccccccccccccc}
 \cdots & \longrightarrow &  0  & \longrightarrow & 0 & \longrightarrow & E_1^{0,q}({\rm N}, {\rm B}) & \longrightarrow  & E_1^{1,q}({\rm N}, {\rm B}) & \longrightarrow &  E_1^{2,q}({\rm N}, {\rm B}) & \longrightarrow & \cdots \\  
 & & \downarrow & & \downarrow & & \downarrow  & & \downarrow  & & \downarrow  & &  \\
 \cdots & \longrightarrow &  0  & \longrightarrow & 0 & \longrightarrow & E_1^{0,q}({\rm N}, {\rm M}) & \longrightarrow  & E_1^{1,q}({\rm N}, {\rm M}) & \longrightarrow &  E_1^{2,q}({\rm N}, {\rm M}) & \longrightarrow & \cdots 
\end{array}
\]
is an isomorphism of cochain complexes, $E_{2}^{1, q}({\rm N}, {\rm M}) \cong E_{2}^{1, q}({\rm N}, {\rm B}) = 0$ by 
Proposition~\ref{prop:E2pqB=0}. 
%As in the proof of Theorem~\ref{th:E21qNisfree}, 
Then there is an surjection 
\begin{eqnarray}
\delta : E_{2}^{0, q}({\rm N}, {\rm M}/{\rm N}) \longrightarrow E_{2}^{1, q}({\rm N}, {\rm N}), \label{eq:connectinghom} 
\end{eqnarray}
%which is the connecting homomorphism associated to 
%\[  
%0 \to E_{1}^{\ast, q}({\rm N}, {\rm N})  \to E_{1}^{\ast, q}({\rm N}, {\rm M}) \to E_{1}^{\ast, q}({\rm N}, {\rm M}/{\rm N}) \to 0. 
%\] 
Using $E_{2}^{0, q}({\rm N}, {\rm M}/{\rm N}) \cong {\rm N}^{!}_{q}\otimes_{R} ((\oplus_{i=1}^{m} RE_{i, i})/RI_m)$, 
we obtain a set of generators 
$\{ y_{I}\otimes E_{ii} \mid  1\le i \le m, |I|= q \}$ of $E_{2}^{0, q}({\rm N}, {\rm M}/{\rm N})$.   
Since $\delta(y_I\otimes E_{ii}) =  y_{i-1}y_I\otimes E_{i-1, i}-(-1)^qy_Iy_i\otimes E_{i, i+1}$, we can verify the statement.   
\qed 

\bigskip
 
By Proposition~\ref{prop:E21qgenerator}, 
the $R$-module $E_2^{1,|I|}({\rm N}, {\rm N})$
is generated by
\[ y_{i-1}y_I\otimes E_{i-1, i}-(-1)^{|I|}y_Iy_i\otimes E_{i, i+1} \]
for $1\le i\le m$. 
%where $y_I=y_{i_1}y_{i_2}\cdots y_{i_{|I|}}$
%for $I=(i_1,\ldots, i_{|I|})$.
Let $\{ I_m^{\ast} \} \cup \{ E_{ij}^{\ast} \mid 1\le i < j \le m \}$ be the dual basis of 
$\{ I_m \} \cup \{ E_{ij} \mid 1\le i < j \le m \}$ of ${\rm N}_m(R)$ over $R$. 
For $I=(i_1,\ldots, i_{|I|})$, 
we set
\[ E_I^*=E_{i_1,i_1+1}^*E_{i_2,i_2+1}^*\cdots E_{i_{|I|},i_{|I|+1}}^*.\]
If $|I|=0$, then set $E^{\ast}_{\emptyset} = 1$.

\begin{lemma}
In the cochain complex $C^*({\rm N}, {\rm N})$,
the cochain
\[ \sum_{1\le k<i}E_{k,i}^*E_I^*\otimes E_{k,i}
   -(-1)^{|I|}
   \sum_{i<k\le m}E_I^*E_{i,k}^*\otimes E_{i,k}\]
is a cocycle.
\end{lemma}

\proof   
The lemma follows from the following calculations:
\[ \begin{array}{cl}
    &d(E_{i,j}^*E_I^*\otimes E_{i,j})\\[2mm]
    =&\displaystyle
    \sum_{k<i}E_{k,i}^*E_{i,j}^*E_I^*\otimes E_{k,j}%\\[2mm]
    %&&\displaystyle
    -\sum_{i<k<j}E_{i,k}^*E_{k,j}^*E_I^*\otimes E_{i,j}%\\[2mm]
    %%&\displaystyle
    +(-1)^{|I|}\sum_{j<k}E_{i,j}^*E_I^*E_{j,k}^*\otimes E_{i,k},\\[8mm]
    &d(E_I^*E_{i,j}^*\otimes E_{i,j})\\[2mm]
    =&\displaystyle
    \sum_{k<i}E_{k,i}^*E_I^*E_{i,j}^*\otimes E_{k,j}%\\[2mm]  
    %&&\displaystyle
    -(-1)^{|I|}\sum_{i<k<j}E_I^*E_{i,k}^*E_{k,j}^*\otimes E_{i,j}%\\[2mm]
    %&&\displaystyle
    +(-1)^{|I|}\sum_{j<k}E_I^*E_{i,j}^*E_{j,k}^*\otimes E_{i,k}.\\
   \end{array}\]
\qed

We define
\[ a(i,I)\in {\rm HH}^{|I|+1,|I|}({\rm N}, {\rm N}) \]
to be the cohomology class represented by
the cocycle
\begin{eqnarray} \sum_{1\le k<i}E_{k,i}^*E_I^*\otimes E_{k,i}
   -(-1)^{|I|}
   \sum_{i<k\le m}E_I^*E_{i,k}^*\otimes E_{i,k}.  \label{eq:aiI} 
\end{eqnarray}

\begin{lemma}\label{lemma:aiI}
The cohomology class
$a(i,I)$ corresponds to
\[ y_{i-1}y_I\otimes E_{i-1, i}-(-1)^{|I|}y_Iy_i\otimes E_{i, i+1} \] 
under the isomorphism
\[ {\rm HH}^{|I|+1,|I|}({\rm N}, {\rm N})\cong
   E_2^{1,|I|}({\rm N}, {\rm N}). \]
\end{lemma}

\proof
This follows from the fact that
\[ \begin{array}{cl}
     &\sum_{1\le k<i}E_{k,i}^*E_I^*\otimes E_{k,i}
     -(-1)^{|I|}
     \sum_{i<k\le m}E_I^*E_{i,k}^*\otimes E_{i,k}\\[2mm]
     \equiv &
     E_{i-1,i}^*E_I^*\otimes E_{i-1,i}
     -(-1)^{|I|}E_I^*E_{i,i+1}^*\otimes E_{i,i+1}\\
   \end{array}\]
in $C^*({\rm N}, {\rm Gr}^1{\rm N})$.
\qed

\bigskip

Proposition~\ref{prop:E2m-1qBgenerated} shows that 
the $R$-module $E_2^{m-1,|J|-(m-1)}({\rm N}, {\rm N})$ is generated by
\[ y_J\otimes E_{1,m}\]
over $R$. 
By the direct calculation,
we obtain the following lemma.

\begin{lemma}\label{lemma:dJ}
In the cochain complex $C^*({\rm N}, {\rm N})$,
the cochain
\[ E_J^*\otimes E_{1,m} \]
is a cocycle.
The cohomology class represented by $E_J^*\otimes E_{1,m}$
corresponds to $y_J\otimes E_{1,m}$
under the isomorphism
${\rm HH}^{|J|,|J|-(m-1)}({\rm N}, {\rm N})\cong 
E_2^{m-1,|J|-(m-1)}({\rm N}, {\rm N})$.
\end{lemma}

We define
\[ d(J)\in {\rm HH}^{|J|,|J|-(m-1)}({\rm N}, {\rm N})\]
to be the cohomology class represented
by $E_J^*\otimes E_{1,m}$.

%\newpage

\subsection{Construction of an $R$-basis of ${\rm HH}^{\ast}({\rm N}_m(R), {\rm N}_m(R))$} 
In this subsection, we construct an $R$-basis of ${\rm HH}^{\ast}({\rm N}_m(R), {\rm N}_m(R))$. 
We set
\[ \mathbf{1}= [1\otimes 1] \in {\rm HH}^{0,0}({\rm N}, {\rm N}).\]
By Theorem~\ref{th:mainthmHHns},  Propositions~\ref{prop:E2m-1qBgenerated} and \ref{prop:E21qgenerator}, and Lemmas~\ref{lemma:aiI} and \ref{lemma:dJ}, ${\rm HH}^{\ast}({\rm N}, {\rm N})$ is generated by 
\[
%{\mathcal S} = 
\{ \mathbf{1} \} \cup \{ a(i, I) \mid 1\le i \le m, |I| \ge 0 \} \cup \{ d(J) \mid |J| \ge 0 \}
\] as $R$-modules.  Since   
\[
\{y_I\otimes E_{1,m}|\  I=(i_1,\ldots,i_q),i_1\neq 1, i_q\neq m-1\}
\]
is an $R$-basis of $E_{2}^{m-1, q-(m-1)}({\rm B})\cong E_{2}^{m-1, q-(m-1)}({\rm N})$ by Proposition~\ref{prop:E2m-1qBgenerated}, 
\[
\{ d(J) \mid  J=(j_1,\ldots,j_q),j_1\neq 1, j_q\neq m-1 \}
\]
is an $R$-basis of $R\{ d(J) \mid |J| =q \}$. 

Let us consider ${\rm HH}^{q+1,q}({\rm N}_m(R), {\rm N}_m(R))$ for $q\ge 0$.  
%Let us choose an $R$-basis of ${\rm HH}^{q+1,q}({\rm N}_m(R), {\rm N}_m(R)) = R\{ a(i, I) \mid 1\le i \le m, |I| = q \}$. 
%$R\{ a(i, I) \mid 1\le i \le m, |I| \ge 0 \}$.  
Set 
\[
{\mathcal S}(q) = \{ I=(i_1,\ldots,i_q) \mid 1 \le i_1, \ldots i_q \le m-1, y_{I}\neq 0 \}
\]
for $q > 0$ and ${\mathcal S}(0) = \{ \emptyset \}$.  Note that $\sharp {\mathcal S}(q) = {\rm rank}_{R} {\rm N}_m(R)^{!}_{q} = \varphi(q)$ for $q \ge 0$. By Proposition~\ref{prop:E21qgenerator} and Lemma~\ref{lemma:aiI}, 
${\rm HH}^{q+1,q}({\rm N}_m(R), {\rm N}_m(R)) \cong E_{2}^{1, q}({\rm N}, {\rm N})$ is generated by $\{ a(i, I) \mid 1 \le i \le m, I \in {\mathcal S}(q) \}$ for $q\ge 0$. 
Recall the long exact sequence (\ref{eq:longexactseqE2}) in Theorem~\ref{th:E21qNisfree}: 
\begin{eqnarray}   
\cdots \longrightarrow E_2^{0,q}({\rm N}, {\rm B})\stackrel{\pi}{\longrightarrow} E_2^{0,q}({\rm N}, {\rm B}/{\rm N}) \stackrel{\delta}{\longrightarrow} E_2^{1,q}({\rm N}, {\rm N})\longrightarrow \cdots . \label{eq:longexact2} 
\end{eqnarray}   
For $1 \le i  \le m$, let $b(i, I) = E^{\ast}_{I}\otimes E_{i, i} \in E_{2}^{0, q}({\rm N}, {\rm B}/{\rm N})$ for $q>0$ and $b(i, \emptyset) = 1\otimes E_{i, i} \in E_{2}^{0, 0}({\rm N}, {\rm B}/{\rm N})$ for $q=0$.  We see that 
$\delta(b(i, I)) = a(i, I)$. %, where $\delta$ is the connecting homomorphism in (\ref{eq:connectinghom}).  
Since $\sum_{i=1}^{m} b(i, I) = E^{\ast}_{I}\otimes I_m =  0 \in E_{2}^{0, q}({\rm N}, {\rm B}/{\rm N})$, 
\begin{eqnarray}
a(1, I)+ a(2, I) +\cdots + a(m, I) = 0. \label{eq:aiIrelation1} 
\end{eqnarray} 
In Definition~\ref{def:notationxI}, we have defined 
\[
z(i, I) = y_iy_{I} \otimes E_{i, i} + (-1)^{q+1} y_{I} y_i \otimes E_{i+1, i+1} \in E_{2}^{0, q+1}({\rm N}, {\rm B}) 
\]
for $I \in {\mathcal S}(q)$ and $1 \le i \le m-1$ (see also Proposition~\ref{prop:E20qBgenerated}).  
Note that $\pi(z(i, I))=b(i, (i, I))+(-1)^{q+1}b(i+1, (I, i))$, where $\pi$ is the $R$-homomorphism 
$\pi : E_2^{0,q}({\rm N}, {\rm B}) 
\to E_2^{0,q}({\rm N}, {\rm B}/{\rm N})$ in (\ref{eq:longexact2}). 
Since $\delta(\pi(z(i, I)))=\delta(b(i, (i, I)))+(-1)^{q+1}\delta(b(i+1, (I, i)))=0$,  
\begin{eqnarray}
a(i, (i, I)) + (-1)^{q+1} a(i+1, (I, i)) = 0 \label{eq:aiIrelation2}  
\end{eqnarray}
for $I \in {\mathcal S}(q)$ and $1 \le i \le m-1$. 

Let us construct an $R$-basis of ${\rm HH}^{q+1,q}({\rm N}, {\rm N}) \cong E_{2}^{1, q}({\rm N}, {\rm N})$. 
Let $q=0$. By Theorem~\ref{th:mainthmHHns}, ${\rm rank}_{R} {\rm HH}^{1,0}({\rm N}, {\rm N}) = m-1$. The set 
$\{ a(i, \emptyset) \mid 1 \le i \le m-1 \}$ is an $R$-basis of ${\rm HH}^{1,0}({\rm N}, {\rm N})$, since 
\[
a(m, \emptyset) = - a(1, \emptyset) - a(2, \emptyset) - \cdots - a(m-1, \emptyset) 
\]
by (\ref{eq:aiIrelation1}).  

Let $q>0$. By Theorem~\ref{th:mainthmHHns}, ${\rm rank}_{R} {\rm HH}^{q+1,q}({\rm N}, {\rm N}) = (m-2)\varphi(q)$. 
Set 
\[{\mathcal T}(q) = \{ (i, I) \mid 1 \le i \le m, I \in {\mathcal S}(q) \}.\]  
Note that $\{ a(i, I) \mid (i, I) \in {\mathcal T}(q) \}$ generates ${\rm HH}^{q+1,q}({\rm N}, {\rm N})$ 
as an $R$-module.  

\begin{definition}\rm 
%Let $\Omega(q) = \{ (i, I) \mid 1\le i \le m, I \in {\mathcal S}(q) \}$ for $q >0$. 
For $q>0$, set
\[
{\mathcal T}(q)_{i} = \{ (i, I) \in {\mathcal T}(q) \mid I = (i, J) \mbox{ for some } J \in {\mathcal S}(q-1) \}   
\]
for $1 \le i \le m-1$. 
We also define 
\begin{eqnarray*}
{\mathcal T}(q)^{0}_{i} & = & \{ (i, (i, J)) \in {\mathcal T}(q)_i \mid y_{J}y_i = 0 \}, \\ 
{\mathcal T}(q)^{1}_{i} & = & \{ (i, (i, J)) \in {\mathcal T}(q)_i \mid y_{J}y_i \neq 0 \}. 
\end{eqnarray*} 
Note that ${\mathcal T}(q)_i = {\mathcal T}(q)^{0}_i \coprod {\mathcal T}(q)^{1}_i$ for $q>0$ and $1 \le i \le m-1$.  
\end{definition}

\begin{lemma}\label{lemma:aiIreduction} 
Let $1 \le i \le m-1$. 
For $(i, (i, J)) \in {\mathcal T}(q)^{0}_i$, $a(i, (i, J)) = 0$. For $(i, (i, J)) \in {\mathcal T}(q)^{1}_i$, 
\begin{eqnarray}
a(i, (i, J))+(-1)^q a(i+1, (J, i))=0. \label{eq:Omegaq1} 
\end{eqnarray}
\end{lemma}
\proof
By direct calculation and (\ref{eq:aiIrelation2}), we can verify the statement. 
\qed 

\begin{definition}\rm 
For $q>0$, set 
\begin{eqnarray*}
{\mathcal T}^{-}(q) & = &  \left( \bigcup_{i=1}^{m-1} {\mathcal T}(q)^{0}_i \right)  \bigcup \left( \bigcup_{i=1}^{m-2} \left\{ (i+1, (J, i)) \mid (i, (i, J)) \in {\mathcal T}(q)^{1}_i \right\} \right) \bigcup  {\mathcal T}(q)^{1}_{m-1} \\ 
& & \bigcup \left(\bigcup_{i=1}^{m-2} \left\{ (m, (i, J)) \mid (i, (i, J)) \in {\mathcal T}(q)_i \right\} \right) \\ 
& & \bigcup \left\{ (m, (m-1, J) \mid (m-1, (m-1, J)) \in {\mathcal T}(q)^{0}_{m-1} \right\} \\
& & \bigcup \left\{ (1, (m-1, J) \mid (m-1, (m-1, J)) \in {\mathcal T}(q)^{1}_{m-1} \right\} 
\end{eqnarray*} 
and 
\[
{\mathcal T}(q)^{+} = {\mathcal T}(q)\setminus {\mathcal T}(q)^{-}. 
\]
Note that $\sharp {\mathcal T}(q) = m\varphi(q)$, $\sharp {\mathcal T}(q)^{-} = 2\sum_{i=1}^{m-1} \sharp {\mathcal T}(q)_i = 2\varphi(q)$, and $\sharp {\mathcal T}(q)^{+} = (m-2)\varphi(q)$ for $q>0$. 
\end{definition}

\begin{proposition}\label{prop:basisofHHq+1q}
For $q>0$, $\{ a(i, I) \mid (i, I) \in {\mathcal T}(q)^{+} \}$ is an $R$-basis of ${\rm HH}^{q+1,q}({\rm N}, {\rm N}) \cong E_{2}^{1, q}({\rm N}, {\rm N})$.
\end{proposition}

\proof 
Set $T=R\{ a(i, I) \mid (i, I) \in {\mathcal T}(q)^{+} \}$. 
Let us show that $T={\rm HH}^{q+1,q}({\rm N}, {\rm N})$. 
It suffices to prove that $a(i, I) \in T$ for any $(i, I) \in  {\mathcal T}(q)^{-}$. 
%Assume that $(i, I) \in  {\mathcal T}(q)^{-}$. 
If $(i, I) \in \left( \bigcup_{i=1}^{m-1} {\mathcal T}(q)^{0}_i \right)$, then $a(i, I) = 0 \in T$ by Lemma~\ref{lemma:aiIreduction}. 
We easily see that ${\mathcal T}(q)^{1}_i \subset {\mathcal T}(q)^{+}$ for $1 \le i \le m-2$.  
Put ${\mathcal T}^{-}(q)' = \bigcup_{i=1}^{m-2} \left\{ (i+1, (J, i)) \mid (i, (i, J)) \in {\mathcal T}(q)^{1}_i \right\}$. 
If $(i, I) \in {\mathcal T}^{-}(q)' $, then $a(i, I) \in T$ 
by (\ref{eq:Omegaq1}). 
%Hence, $a(i, I) \in R\{ a(i, I) \mid (i, I) \in {\mathcal T}(q)^{+} \}$ for any $(i, I) \in {\mathcal T}(q)_i = {\mathcal T}(q)^{0}_i \cup {\mathcal T}(q)^{1}_i$. 
If $(i, I) \in \bigcup_{i=1}^{m-2} \left\{ (m, (i, J)) \mid (i, (i, J)) \in {\mathcal T}(q)_i \right\}$, then   $a(i, I) \in T$, since 
\[
R\left\{ 
\begin{array}{c|c}
a(i, (i, J)) & \displaystyle 
(i, (i, J)) \in \bigcup_{i=1}^{m-2} {\mathcal T}(q)_i 
\end{array}
\right\} + R\{ a(i, I) \mid (i, I) \in {\mathcal T}^{-}(q)' \} \subseteq T
\] 
and 
\begin{eqnarray}
a(i, I) = - \sum_{j\neq i} a(j, I)  \label{eq:aiIrelation3}
\end{eqnarray}
by (\ref{eq:aiIrelation1}).  
For $(i, I) \in \left\{ (m, (m-1, J) \mid (m-1, (m-1, J)) \in {\mathcal T}(q)^{0}_{m-1} \right\}$, we see that $a(i, I) \in T$ by 
using $a(m-1,(m-1, J))=0$, $a(j, (m-1, J)) \in T$ for $1 \le j \le m-2$ and  (\ref{eq:aiIrelation3}). 
Summarizing the discussion above, we have $a(i, I) \in T$ for any $(i, I) \in {\mathcal T}(q)$ with $i\neq 1, m-1$.

For $(m-1, (m-1, J)) \in {\mathcal T}(q)^{1}_{m-1}$, let us prove that $a(1, (m-1, J)), a(m-1, (m-1, J)) \in T$.  
By (\ref{eq:Omegaq1}) and  (\ref{eq:aiIrelation3}), 
\begin{eqnarray*}
a(m-1, (m-1, J)) + (-1)^{q}a(m, (J, m-1))  & = &  0, \\
a(1, (m-1, J)) + a(2, (m-1, J))+ \cdots + a(m-1, (m-1, J)) + a(m, (m-1, J)) & = & 0. 
\end{eqnarray*}
Using $a(j, (m-1, J)), a(m, (J, m-1)) \in T$ for $j \neq 1, m-1$, we obtain $a(1, (m-1, J)), a(m-1, (m-1, J)) \in T$.  

Hence, $a(i, I) \in T$ for any $(i, I) \in  {\mathcal T}(q)^{-}$ and $T={\rm HH}^{q+1,q}({\rm N}, {\rm N})$. 
%
%We also see that $a(i, I) \in R\{ a(i, I) \mid (i, I) \in {\mathcal T}(q)^{+} \}$ for 
%$(i, I) \in \left\{ (1, (m-1, J) \mid (m-1, (m-1, J)) \in {\mathcal T}(q)^{1}_{m-1} \right\}$.  
%Thus, we have shown that $\{ a(i, I) \mid (i, I) \in {\mathcal T}(q)^{+} \}$ generates ${\rm HH}^{q+1,q}({\rm N}, {\rm N})$. 
Since ${\rm rank}_{R}{\rm HH}^{q+1,q}({\rm N}, {\rm N}) = \sharp {\mathcal T}(q)^{+} = (m-2)\varphi(q)$,  $\{ a(i, I) \mid (i, I) \in {\mathcal T}(q)^{+} \}$ is an $R$-basis of ${\rm HH}^{q+1,q}({\rm N}, {\rm N})$. 
\qed 

%\bigskip 
%
%By Proposition~\ref{prop:E2m-1qBgenerated}, 
%\[ {\rm HH}^{q, -(m-1)+q}({\rm N}, {\rm N}) \cong E_2^{m-1,-(m-1)+q}({\rm B})=
%   R\{y_I\otimes E_{1,m}|\ 
%   I=(i_1,\ldots,i_q),i_1\neq 1, i_q\neq m-1\}. \]
%The set $\{ d(J) \mid J=(j_1, \ldots, j_q), j_1 \neq 1, j_{q} \neq m-1 \}$ is an $R$-basis of ${\rm HH}^{q, -(m-1)+q}({\rm N}, {\rm N})$. 

\bigskip 

By the discussion above, we obtain the following corollary. 

\begin{corollary}\label{cor:basisof HHNN}
We have an $R$-basis 
\begin{eqnarray*}
{\mathcal S} & = & \{ \mathbf{1} \} \bigcup \{ a(1, \emptyset), \ldots, a(m-1, \emptyset) \} \bigcup \left(\bigcup_{q > 0} \{ a(i, I) \mid (i, I) \in {\mathcal T}(q)^{+} \} \right) \\
& & \hspace*{5ex} \bigcup \left(\bigcup_{q\ge 0} \{ d(J) \mid  J=(j_1,\ldots,j_q),j_1\neq 1, j_q\neq m-1 \} \right) 
\end{eqnarray*} 
of ${\rm HH}^{\ast}({\rm N}_m(R), {\rm N}_m(R))$. 
\end{corollary}

\subsection{Gerstenhaber bracket}
In this subsection, we calculate the Gerstenhaber bracket $[\; ,\; ]$ of  ${\rm HH}^{\ast}({\rm N}_m(R), {\rm N}_m(R))$.    
By Corollary~\ref{cor:basisof HHNN}, we have the $R$-basis ${\mathcal S}$ of ${\rm HH}^{\ast}({\rm N}_m(R), {\rm N}_m(R))$.  
%We set
%\[ \mathbf{1}= [1\otimes 1] \in {\rm HH}^{0,0}({\rm N}, {\rm N}).\]
%By Lemmas~\ref{lemma:aiI} and \ref{lemma:dJ}, ${\rm HH}^{\ast}({\rm N}, {\rm N})$ is generated by 
%\[
%{\mathcal S} = \{ \mathbf{1} \} \cup \{ a(i, I) \mid 1\le i \le m, |I| \ge 0 \} \cup \{ d(J) \mid |J| \ge 0 \}
%\] as $R$-modules.  
For investigating the Gerstenhaber bracket on ${\rm HH}^{\ast}({\rm N}, {\rm N})$, 
we only need to calculate $[x, y]$ for $x, y \in {\mathcal S}$. 

\bigskip

We easily obtain the following lemmas.

\begin{lemma}
For any $z\in {\rm HH}^{*,*}({\rm N}, {\rm N})$,
we have
\[ [\mathbf{1},z]=0. \]
\end{lemma}

\begin{lemma}
For any $J,J'$, we have
\[ [d(J),d(J')]=0. \]
\end{lemma}

\bigskip

Next, we will calculate
$[d(J),a(i,I)]$.

For $I=(i_1,\ldots,i_{|I|})$,
$J=(j_1,\ldots,j_{|J|})$,
and $1\le k\le |I|$,
we set
\[ J\subrel{k}{\circ}I=    
   (j_1,\ldots,j_{k-1},i_1,\ldots,i_{|I|},j_{k+1},\ldots,j_{|J|}).
\]
We also set
\[ J(r)=\{k\in \{1,2,\ldots,|J|\}|\ j_k=r\}. \]

\begin{proposition}
We have 
\[ \begin{array}{cl}
   &[d(J),a(i,I)]\\[5mm]
   =&
   \left\{
   \begin{array}{ll}
   \displaystyle
   (-1)^{|I|}d(I,J)   
   -\sum_{k\in J(1)} (-1)^{k|I|}
   d(J\subrel{k}{\circ}(I,1))& (i=1),\\[10mm]

   \displaystyle 
   \sum_{k\in J(i-1)}(-1)^{(k-1)|I|} 
   d(J\subrel{k}{\circ}(i-1,I))%\\[5mm]
   %&&\displaystyle 
   -\sum_{k\in J(i)}(-1)^{k|I|}
   d(J\subrel{k}{\circ}(I,i)) & (1<i<m),\\[10mm]

   \displaystyle 
   \sum_{k\in J(m-1)}(-1)^{(k-1)|I|}   
   d(J\subrel{k}{\circ}(m-1,I))
   -(-1)^{|I|(|J|-1)} d(J,I) & (i=m).\\%[10mm]
   \end{array}\right.\\
   \end{array}
\]
\end{proposition}

%\newpage

Next, we will calculate
$[a(i,I),a(i',I')]$.

\begin{lemma}
\label{lemma:circle-product-truncation}
Let $x=x^1+x^2+\cdots +x^{m-1} \in C^{p+1, p}({\rm N}, {\rm N})\cap \overline{C}^{p+1}({\rm N}, {\rm N})$ and
$y=y^1+y^2+\cdots +y^{m-1} \in C^{q+1, q}({\rm N}, {\rm N}) \cap \overline{C}^{q+1}({\rm N}, {\rm N})$ be cocycles, 
where $x^i \in C^{p+1}({\rm N},F^i{\rm N})$ and $y^i\in C^{q+1}({\rm N},F^i{\rm N})$.
%Suppose that
%$[x],[y]\in {\rm HH}^{*+1,*}({\rm N}, {\rm N})$.
Then we have 
\[ x\circ y\equiv x^1\circ y^1 \mod C^{p+q+1}({\rm N},F^2{\rm N}).\]
%Let $x=x^1+x^2+\cdots +x^{m-1}$ and
%$y=y^1+y^2+\cdots +y^{m-1}$ be cocycles in $\overline{C}^*({\rm N}, {\rm N}) \subset C^*({\rm N}, {\rm N})$,
%where $x^i,y^i\in C^*({\rm N},F^i{\rm N})$.
%Suppose that
%$[x],[y]\in {\rm HH}^{*+1,*}({\rm N}, {\rm N})$.
%Then we have 
%\[ x\circ y\equiv x^1\circ y^1 \mod C^*({\rm N},F^2{\rm N}).\]
%Let $x=x^1+x^2+\cdots +x^{m-1}$ and
%$y=y^1+y^2+\cdots +y^{m-1}$ be cocycles in $C^*({\rm N}, {\rm N})$,
%where $x^i,y^i\in C^*({\rm N},F^i{\rm N})$.
%Suppose that
%$[x],[y]\in {\rm HH}^{*+1,*}({\rm N}, {\rm N})$.
%Then we have 
%\[ x\circ y\equiv x^1\circ y^1 \mod C^*({\rm N},F^2{\rm N}).\]
\end{lemma}

\proof
Since $x^1 \in  C^{p+1, p}({\rm N}, {\rm N})\cap \overline{C}^{p+1}({\rm N}, {\rm N})$, 
$x^1$ is a linear combination of 
\[\{ E^{\ast}_{i_1, i_1+1}E^{\ast}_{i_2, i_2+1}\cdots E^{\ast}_{i_{p+1}, i_{p+1}+1}\otimes E_{i, i+1} \mid 1 \le i_1, i_2, \ldots, i_{p+1}, i \le m-1 \}
\] 
modulo $C^{p+1}({\rm N},F^2{\rm N})$.   
It is easy to see that $x^1\circ (y^2+\cdots +y^{m-1})\equiv 0 \mod C^{p+q+1}({\rm N},F^2{\rm N})$. Hence, we can verify the statement. 
%The lemma follows from the fact that
%\[ x^1\circ y^i\in C^*({\rm N},F^2{\rm N})\quad (i>1) \]
%by degree reasons.
\qed

\begin{lemma}
\label{lemma:braket-representetive}
Let $x=x^1+x^2+\cdots+x^{m-1} \in C^{p+1, p}({\rm N}, {\rm N})\cap \overline{C}^{p+1}({\rm N}, {\rm N})$ and
$y=y^1+y^2+\cdots+y^{m-1} \in C^{q+1, q}({\rm N}, {\rm N}) \cap \overline{C}^{q+1}({\rm N}, {\rm N})$ be cocycles,  
%in $\overline{C}^*({\rm N}, {\rm N}) \subset C^*({\rm N}, {\rm N})$,
where $x^i \in C^{p+1}({\rm N},F^i{\rm N})$ and $y^i\in C^{q+1}({\rm N},F^i{\rm N})$. 
%where $x^i,y^i\in C^*({\rm N},F^i{\rm N})$.
%Suppose that
Let 
$[x],[y]\in {\rm HH}^{*+1,*}({\rm N},{\rm N})$ be the cohomology classes represented by $x, y$, respectively.
Then the element of
$E_2^{1,p+q}({\rm N}, {\rm N})$ that corresponds to 
the Gerstenhaber bracket $[[x],[y]]\in {\rm HH}^{p+q+1,p+q}({\rm N}, {\rm N})$
is represented by
$x^1\circ y^1-(-1)^{(|x|-1)(|y|-1)}y^1\circ x^1$.
\end{lemma}

\proof
By Lemma~\ref{lemma:circle-product-truncation},
$[x,y]\equiv x^1\circ y^1-(-1)^{(|x|-1)(|y|-1)}y^1\circ x^1$
mod $C^{p+q+1}({\rm N},F^2{\rm N})$.
The proposition follows from
the isomorphism
${\rm HH}^{p+q+1, p+q}({\rm N}, {\rm N})\cong E_2^{1,p+q}({\rm N}, {\rm N})$
since $[[x],[y]] \in {\rm HH}^{p+q+1,p+q}({\rm N},{\rm N})$.
\qed

\bigskip

We set
\[ \begin{array}{cl}
    &A(i,I;i',I')\\[5mm]
   =&\displaystyle
    \sum_{k\in I(i'-1)}(-1)^{k|I'|} 
    a(i,I\subrel{k}{\circ}(i'-1,I'))
   -\sum_{k\in I(i')}(-1)^{(k+1)|I'|}
   a(i,I\subrel{k}{\circ}(I',i'))\\[5mm]
   &\displaystyle
    -(-1)^{|I||I'|}
    \left(\sum_{k\in I'(i-1)}(-1)^{k|I|} 
    a(i',I'\subrel{k}{\circ}(i-1,I))
   -\sum_{k\in I'(i)}(-1)^{(k+1)|I|}
   a(i',I'\subrel{k}{\circ}(I,i))
   \right).\\[5mm]
   \end{array}\]
   
By direct calculation, we have the following proposition. 

\begin{proposition}\label{prop:Gerstenhaberaa} 
We have
\[ [a(i,I),a(i',I')]=
   \left\{\begin{array}{ll}
   A(i,I;i',I') & (i\neq i'),\\[2mm]
   A(i,I;i,I')
   +a(i,(I',I))-(-1)^{|I||I'|}a(i,(I,I')) &
                   (i=i').\\     
          \end{array}\right.
   \]
\end{proposition}

%
%\begin{remark}\rm
%We can delete 
%Lemmas~\ref{lemma:circle-product-truncation} 
%and \ref{lemma:braket-representetive}.
%\end{remark}}

\subsection{Batalin-Vilkovisky structure on ${\rm HH}^{\ast}({\rm N}_m(R), {\rm N}_m(R))$}

Recall that a Batalin-Vilkovisky algebra is a Gerstenhaber algebra $(A^{\ast}, \cup, [\;, \;])$ with an operator 
$\Delta : A^{\ast} \to A^{\ast -1}$ of degree $-1$ such that $\Delta\circ \Delta = 0$ and 
\begin{eqnarray}
[a, b] & = & (-1)^{|a|}\{ \Delta(a\cup b) - \Delta(a)\cup b -(-1)^{|a|}a\cup \Delta(b) \} 
\label{eq:BV}  
\end{eqnarray}
for homogeneous elements $a, b \in A^{\ast}$ (see, for example, \cite[Definition~3.6]{BKL}).   In this subsection,  
we consider the question whether the Hochschild cohomology ${\rm HH}^{\ast}({\rm N}_m(R), {\rm N}_m(R))$ has a Batalin-Vilkovisky algebra structure over $R$ which gives the Gerstenhaber bracket $[\;, \;]$ or not. 

\begin{lemma}\label{lemma:noBV}
Let $A$ be an associated algebra over $R$ such that $A$ is a projective module over $R$.  
Assume that ${\rm HH}^{k}(A, A) \cup {\rm HH}^{l}(A, A)  = 0$ for any $k, l > 0$.  If there exist $a \in {\rm HH}^{k}(A, A)$ and $b \in {\rm HH}^{l}(A, A)$ with $k, l \ge 2$ such that $[a, b]\neq 0$, then ${\rm HH}^{\ast}(A, A)$ has no Batalin-Vilkovisky algebra structure over $R$ which gives the Gerstenhaber bracket  $[\;, \;]$. 
\end{lemma}

\proof 
Suppose that ${\rm HH}^{\ast}(A, A)$ has a Batalin-Vilkovisky algebra structure which gives the Gerstenhaber bracket $[\;, \;]$. By (\ref{eq:BV}) and $a\cup b=\Delta(a)\cup b = a\cup \Delta(b)=0$, we obtain $[a, b] = 0$, which is a contradiction.  Hence, ${\rm HH}^{\ast}(A, A)$ has no Batalin-Vilkovisky algebra structure which gives $[\;, \;]$. 
\qed 

\bigskip 

Let us show that ${\rm HH}^{\ast}({\rm N}_m(R), {\rm N}_m(R))$ has no Batalin-Vilkovisky algebra structure over $R$ giving $[\;, \;]$ for $m\ge 3$. %For this purpose, we make some preparations. 

\begin{lemma}\label{lemma:exampleGerstenhaber}
Let $m \ge 3$.  For $a(1, (1, 1)), a(1, (2, 1)) \in {\rm HH}^{3, 2}({\rm N}_m(R), {\rm N}_m(R))$, 
% $\cong E_{2}^{1, 2}({\rm N}_m(R), {\rm N}_m(R))$, 
\[[a(1, (1, 1)), a(1, (2, 1))] = a(1, (2, 1, 1, 1)) \neq 0.
\]
\end{lemma}

\proof
By Proposition~\ref{prop:Gerstenhaberaa}, $[a(1, (1, 1)), a(1, (2, 1))] = a(1, (2, 1, 1, 1))$. 
Since $(1, (2, 1, 1, 1)) \in {\mathcal T}(4)^{+}$, $a(1, (2, 1, 1, 1)) \neq 0$ by Proposition~\ref{prop:basisofHHq+1q}.  
\qed 

\begin{theorem}\label{th:noBVstructure} 
For $m \ge 3$, ${\rm HH}^{\ast}({\rm N}_m(R), {\rm N}_m(R))$ has no Batalin-Vilkovisky algebra structure over $R$ which gives the Gerstenhaber bracket  $[\;, \;]$. 
\end{theorem}

\proof
The statement follows from Theorems~\ref{th:productzerom=3} and \ref{th:productzerom>=4}, and Lemmas \ref{lemma:noBV} and \ref{lemma:exampleGerstenhaber}.  
\qed 

\section{Appendix: the case $m=2$}\label{section:m=2} 
In this appendix, we deal with ${\rm N}_2(R)$ for a commutative ring $R$. Set ${\rm N} = {\rm N}_2(R)$. 
Putting $x = E_{1,2} \in {\rm N}$, we see that ${\rm N} \cong R[x]/(x^2)$. 
Throughout this section, we set ${\rm Ann}(2)=\{ a \in R \mid 2a=0 \}$. 
We introduce the following proposition without proof, which gives a projective resolution of ${\rm N}$ over ${\rm N}^{e} = {\rm N}\otimes_{R}{\rm N}^{op}$ over $R$.  

\begin{proposition}[{\cite[Proposition~1.3]{Buenos}, \cite[Example~2.6]{Redondo}}]\label{prop:projresonegen}   
The following complex gives a projective resolution of ${\rm N}$ over ${\rm N}^{e}$: 
\begin{eqnarray} 
\cdots \longrightarrow {\rm N}^{e} \stackrel{d_n}{\longrightarrow} {\rm N}^{e} \longrightarrow \cdots \longrightarrow {\rm N}^{e} \stackrel{d_1}{\longrightarrow} {\rm N}^{e} \stackrel{\mu}{\longrightarrow} {\rm N} \longrightarrow 0, \label{eq:projresolN2}  
\end{eqnarray}
where 
\[ d_i(a\otimes b) = \left\{ 
\begin{array}{cc}
(x\otimes 1 + 1 \otimes x) (a\otimes b) & (i : \mbox{\rm even}),  \\
(x \otimes 1-1\otimes x) (a\otimes b) & (i : \mbox{\rm odd}) \\
\end{array} 
\right. 
\]  
and $\mu(a\otimes b) = ab$. 
\end{proposition}  

In \cite{Nakamoto-Torii:Hochschild}, 
we have calculated ${\rm HH}^{n}({\rm N}_2(R), {\rm M}_2(R)/{\rm N}_2(R))$ by using the projective resolution above.   

\begin{theorem}[{\cite[Proposition~4.19]{Nakamoto-Torii:Hochschild}}]
We have 
\[ 
{\rm HH}^{n}({\rm N}_2(R), {\rm M}_2(R)/{\rm N}_2(R)) \cong  
\left\{ 
\begin{array}{cc}
R \oplus {\rm Ann}(2) & (n:\mbox{\rm even}),  \\
R \oplus (R/2R) & (n:\mbox{\rm odd}). \\
\end{array} 
\right. 
\] 
\end{theorem} 

\begin{corollary}[{\cite[Corollary~4.20]{Nakamoto-Torii:Hochschild}}]\label{cor:hochschildjn}
Let $k$ be a field. For each $n \ge 0$,  
\[ 
{\rm HH}^{n}({\rm N}_2(k), {\rm M}_2(k)/{\rm N}_2(k)) \cong  
\left\{ 
\begin{array}{cc}
k & ({\rm ch}(k) \neq 2),  \\
k^2 & ({\rm ch}(k) = 2).  \\
\end{array} 
\right. 
\] 
\end{corollary} 

\bigskip 

By using the same discussions in \S\ref{subsection:sssubquotient} and \S\ref{subsection:calHHM/N}, we also have the following result. 

\begin{theorem}
For each $n \ge 0$ and $s \in {\Bbb Z}$, 
\[
{\rm HH}^{n, s}({\rm N}_2(R), {\rm M}_2(R)/{\rm N}_2(R)) = \left\{ 
\begin{array}{ll}
R & (n: \mbox{\rm even}, s=n), \\
{\rm Ann}(2) & (n: \mbox{\rm even}, s=n+1), \\
R/2R & (n: \mbox{\rm odd}, s=n), \\    
R & (n: \mbox{\rm odd}, s=n+1), \\    
0 & (\mbox{\rm otherwise}).  
\end{array}
\right. 
\] 
\end{theorem}

Next, let us consider ${\rm HH}^{n}({\rm N}_2(R), {\rm N}_2(R))$. 
By taking ${\rm Hom}_{{\rm N}^e}(-, {\rm N})$ of (\ref{eq:projresolN2}), we obtain the following complex 
\[
0 \longrightarrow {\rm Hom}_{{\rm N}^e}({\rm N}^e, {\rm N}) \stackrel{d_1^{\ast}}{\longrightarrow} 
{\rm Hom}_{{\rm N}^e}({\rm N}^e, {\rm N}) \stackrel{d_2^{\ast}}{\longrightarrow} 
{\rm Hom}_{{\rm N}^e}({\rm N}^e, {\rm N}) \stackrel{d_3^{\ast}}{\longrightarrow} \cdots,  
\]
which is isomorphic to 
\[
0 \longrightarrow {\rm N} \stackrel{\delta'^1}{\longrightarrow} {\rm N} \stackrel{\delta'^2}{\longrightarrow} {\rm N} \stackrel{\delta'^3}{\longrightarrow} \cdots,  
\]
where $\delta'^i : {\rm N} \to {\rm N}$ is defined by 
\[ \delta'^i(a) = \left\{ 
\begin{array}{cc}
2xa & (i : \mbox{even}), \\
0 & (i : \mbox{odd}). 
\end{array} 
\right. 
\]  
Thus, we obtain 
\begin{theorem}\label{th:HHN2N2}
We have 
\[ 
{\rm HH}^{n}({\rm N}_2(R), {\rm N}_2(R)) \cong  
\left\{ 
\begin{array}{ll}
{\rm N}_2(R) & (n=0), \\
{\rm N}_2(R)/(2E_{1, 2}{\rm N}_2(R)) \cong R\oplus (R/2R) & (n:\mbox{\rm even}, n>0), \\
RE_{1, 2} \oplus {\rm Ann}(2)I_2 \cong R\oplus {\rm Ann}(2) & (n:\mbox{\rm odd}). 
\end{array} 
\right. 
\] 
\end{theorem}

Notice that ${\rm HH}^{n}({\rm N}_2(R), {\rm N}_2(R))$ is not a free $R$-module in general, which is different from the case ${\rm N}_m(R)$ for $m \ge 3$.  

Third, let us consider the product structure on ${\rm HH}^{n}({\rm N}_2(R), {\rm N}_2(R))$. 
Set $\overline{\rm N} = {\rm N}/RI_2 \cong Rx$. 
Recall the reduced bar complex $\overline{B}_p({\rm N}, {\rm N}, {\rm N}) = {\rm N}\otimes_{R} \overline{\rm N}^{\otimes p}\otimes_{R} {\rm N}$.   
Let us consider a homomorphism of chain complexes 
\begin{eqnarray}
\begin{array}{ccccccccccc}
\cdots &  \longrightarrow &  {\rm N}\otimes_{R}\overline{\rm N}^{\otimes 2}\otimes_{R}{\rm N}& \longrightarrow &  {\rm N}\otimes_{R}\overline{\rm N}\otimes_{R}{\rm N} & \longrightarrow &  {\rm N}\otimes_{R}{\rm N} & \stackrel{\mu}{\longrightarrow} & {\rm N} & \longrightarrow & 0 \\ 
& & f_2\downarrow\hspace*{3ex}  & & f_1\downarrow\hspace*{3ex}  & & f_0\downarrow\hspace*{3ex}  & & \parallel & &   \\ 
\cdots &  \stackrel{d_3}{\longrightarrow} & {\rm N}^{e} & \stackrel{d_2}{\longrightarrow} &  {\rm N}^{e} & \stackrel{d_1}{\longrightarrow} &  {\rm N}^{e} & \stackrel{\mu}{\longrightarrow} & {\rm N} & \longrightarrow & 0, \\ 
\end{array} \label{eq:chainmaptoN}
\end{eqnarray} 
where $f_p : \overline{B}_p({\rm N}, {\rm N}, {\rm N}) = {\rm N}\otimes_{R} \overline{\rm N}^{\otimes p}\otimes_{R} {\rm N} \to {\rm N}^e$ is the ${\rm N}^e$-homomorphism defined by $f_p(x^{\otimes p}) = I_2\otimes I_2$ for 
$p \ge 0$. 
By taking ${\rm Hom}_{{\rm N}^e}(-, {\rm N})$ of (\ref{eq:chainmaptoN}), we have a quasi-isomorphism of cochain complexes 
\[
\begin{array}{ccccccccc}
0 & \longrightarrow & {\rm N} & \stackrel{\delta'^1}{\longrightarrow} & {\rm N}  & \stackrel{\delta'^2}{\longrightarrow} & {\rm N}   & \stackrel{\delta'^3}{\longrightarrow} & \cdots \\
 & & \downarrow & & \downarrow & & \downarrow & &  \\
0 & \longrightarrow & \overline{C}^{0}({\rm N}, {\rm N}) & \stackrel{\delta^1}{\longrightarrow} & \overline{C}^{1}({\rm N}, {\rm N})    & \stackrel{\delta^2}{\longrightarrow} & \overline{C}^{2}({\rm N}, {\rm N}) & \stackrel{\delta^3}{\longrightarrow} & \cdots,   
\end{array} 
\]
where 
\[\overline{C}^{p}({\rm N}, {\rm N}) = {\rm Hom}_{{\rm N}^e}({\rm N}\otimes_{R} \overline{\rm N}^{\otimes p}\otimes_{R} {\rm N}, {\rm N}) \cong {\rm Hom}_{R}(\overline{\rm N}^{\otimes p}, {\rm N}). \]
For $p \ge 0$, we define $f_p, g_p \in \overline{C}^{p}({\rm N}, {\rm N})$ by 
\begin{eqnarray*}
f_p(x^{\otimes p}) & = &  I_2, \\
g_p(x^{\otimes p}) & = &  E_{1,2},  
\end{eqnarray*}
respectively. Then $\overline{C}^{p}({\rm N}, {\rm N}) = Rf_p \oplus Rg_p$. 
By rephrasing Theorem~\ref{th:HHN2N2}, we obtain 
\begin{theorem}
For $n \ge 0$, we have 
\[
{\rm HH}^{n}({\rm N}_2(R), {\rm N}_2(R)) \cong 
\left\{
\begin{array}{ll}
Rf_0 \oplus Rg_0 & (n=0), \\
Rf_n \oplus (R/2R)g_n & (n:\mbox{\rm even}, n>0), \\
{\rm Ann}(2)f_n \oplus Rg_n & (n:\mbox{\rm odd}).  
\end{array}
\right.
\]
\end{theorem}

By direct calculation, we obtain the following theorem. 
\begin{theorem}\label{th:productHHN2N2}
For any $a \in {\rm HH}^{\ast}({\rm N}_2(R), {\rm N}_2(R))$, $f_0 a = a f_0 =a$. 
For $i, j \ge 0$, we have 
\[
\begin{array}{ccccc}
f_if_j & = & f_jf_i  & = &  f_{i+j}, \\ 
f_ig_j & = & g_j f_i & = & g_{i+j}, \\ 
& & g_i g_j & = & 0.
\end{array}
\]
\end{theorem}

\begin{remark}\rm
For an odd integer $i > 0$, if $a_i \in {\rm Ann}(2)$, then $a_i f_i \in {\rm HH}^{i}({\rm N}_2(R), {\rm N}_2(R)) \cong {\rm Ann}(2)f_i \oplus Rg_i$.   
For an even integer $j >0$, $(a_i f_i) g_j = a_i g_{i+j} \in {\rm HH}^{i+j}({\rm N}_2(R), {\rm N}_2(R)) \cong {\rm Ann}(2) f_{i+1}\oplus R g_{i+j}$ is well-defined.  
\end{remark} 

\begin{remark}\rm
Theorem~\ref{th:productHHN2N2} is compatible with the result in \cite[Theorem~7.1]{Holm}: 
Let $k$ be a commutative ring with ${\rm ch}(k) = p$, where $p$ is a prime number or $0$.  
For $A_2 =k[X]/(X^2)$,  the Hochschild cohomology ring of $A_2$ has the following structure
\[
{\rm HH}^{*}(A_2) \cong 
\left\{
\begin{array}{ll}
 k[x, y, z]/(x^2, y^2-z) & \mbox{ if } p=2 \\
 k[x, y, z]/(x^2, 2xz, yx, y^2) & \mbox{ if }  p\neq 2 \mbox{ and } 2 \in k^{\times},    
\end{array}
\right.
\] 
where $\deg x=0, \deg y=1$, and $\deg z=2$. 
(Note that $t$ in \cite[Theorem~7.1]{Holm} is needed to be regarded as $n$.)  
\end{remark}

\begin{remark}\rm
We can show that ${\rm HH}^{\ast}({\rm N}_2(R), {\rm N}_2(R))$ is a finitely generated algebra over $R$ if and only if 
${\rm Ann}(2)$ is a finitely generated ideal of $R$. Indeed, the ``only if'' part follows from that ${\rm HH}^{1}({\rm N}_2(R), {\rm N}_2(R)) = R\oplus {\rm Ann}(2)$. If ${\rm Ann}(2) = Ra_1+\cdots + Ra_s$, then  ${\rm HH}^{\ast}({\rm N}_2(R), {\rm N}_2(R))$ is generated by 
\[
\{ f_0, g_0, a_1f_1, \ldots, a_sf_1, g_1, f_2 \}
\]
as an $R$-algebra. In particular, if $R$ is a noetherian ring, then ${\rm HH}^{\ast}({\rm N}_2(R), {\rm N}_2(R))$ is a finitely generated algebra over $R$.  
\end{remark} 

By calculating $E_{2}^{p, q, s}({\rm N}_2(R), {\rm N}_2(R))$ in \S\ref{subsection:Degss} directly, 
%Using the same discussions in %\S\ref{subsection:sssubquotient} 
%\S\ref{subsection:Degss}, %and \S\ref{subsection:calHHM/N}, 
we also have the following result. 

\begin{theorem}
For each $n \ge 0$ and $s \in {\Bbb Z}$, 
\[
{\rm HH}^{n, s}({\rm N}_2(R), {\rm N}_2(R)) = \left\{ 
\begin{array}{ll}
R & (n=0, s=-1), \\    
R & (n: \mbox{\rm even}, n=s\ge 0), \\
R/2R & (n: \mbox{\rm even}, n=s+1\ge 2), \\    
{\rm Ann}(2) & (n: \mbox{\rm odd}, n=s\ge 1), \\
R & (n: \mbox{\rm odd}, n=s+1\ge 1), \\  
0 & (\mbox{\rm otherwise}).  
\end{array}
\right. 
\] 
\end{theorem}

\bigskip 

Finally, let us consider the Gerstenhaber bracket on ${\rm HH}^{n}({\rm N}_2(R), {\rm N}_2(R))$. 
We can easily verify the following result. 

\begin{theorem}
We have 
\[
\begin{array}{ccl}
[f_i, f_j] & = & 0 \hspace*{3ex} (i, j \ge 0), \\ 
{[f_i,g_j]} &= & 
\left\{
\begin{array}{ll}
0 & (i:\mbox{\rm even}, \; j:\mbox{\rm even}), \\
f_{i+j-1} & (i:\mbox{\rm odd}, \; j:\mbox{\rm even}), \\
if_{i+j-1} & (j:\mbox{\rm odd}), 
\end{array}
\right. \\
{[g_i,g_j]} &= & 
\left\{
\begin{array}{ll}
0 & (i:\mbox{even}, \; j:\mbox{\rm even}), \\
-(j-1)g_{i+j-1} & (i:\mbox{\rm odd}, \; j:\mbox{\rm even}), \\
(i-1)g_{i+j-1} & (i:\mbox{\rm even}, \; j:\mbox{\rm odd}), \\
(i-j)g_{i+j-1} &(i:\mbox{\rm odd}, \; j:\mbox{\rm odd}). 
\end{array}
\right. 
\end{array}
\]
\end{theorem}

Suppose that  $R$ is a field $k$ of characteristic ${\rm ch}(k) \neq 2$. Then 
\[
{\rm HH}^{n}({\rm N}_2(k), {\rm N}_2(k)) \cong 
\left\{
\begin{array}{ll}
kf_0 \oplus kg_0 & (n=0), \\
kf_n  & (n:\mbox{even}, n>0), \\
kg_n & (n:\mbox{odd}),   
\end{array}
\right.
\]
where $f_0$ is the unit and $f_{2n}=f_2^n$, $g_{2n+1}=f_2^n g_1=g_1f_2^n$, $f_{2}g_0 = g_0f_{2}=g_{0}g_{1}=g_{1}g_{0}=0$ for $n\ge 0$. 
In particular, ${\rm HH}^{\ast}({\rm N}_2(k), {\rm N}_2(k))$ is generated by $g_0, g_1, f_2$ as a $k$-algebra.  

\begin{theorem}
Let $k$ be a field of characteristic ${\rm ch}(k) \neq 2$. 
For $c \in k$, define an operator $\Delta_c : {\rm HH}^{\ast}({\rm N}_2(k), {\rm N}_2(k)) \to {\rm HH}^{\ast-1}({\rm N}_2(k), {\rm N}_2(k))$ by 
\begin{eqnarray*}
\Delta_c(f_0)=\Delta_c(g_0)=0, &\\
\Delta_c(g_1) = f_0+cg_0, & \\
\Delta_c(f_{2n}) = \Delta_c(f_2^n) = 0 & \hspace*{2ex} (n > 0), \\
\Delta_c(g_{2n+1}) = \Delta_c(f_2^ng_1) = (2n+1)f_2^n & \hspace*{2ex} (n \ge 0). 
\end{eqnarray*}
Then $\Delta_c$ gives ${\rm HH}^{\ast}({\rm N}_2(k), {\rm N}_2(k))$ a Batalin-Vilkovisky algebra structure which induces 
$[\;, \;]$.  
In particular, Batalin-Vilkovisky algebra structures on ${\rm HH}^{\ast}({\rm N}_2(k), {\rm N}_2(k))$ giving 
$[\;, \;]$ are not unique.   
\end{theorem}

\proof
By direct calculation, we can verify (\ref{eq:BV}). 
\qed 

\bigskip 

Suppose that  $R$ is a field $k$ of characteristic ${\rm ch}(k) = 2$. Then 
\[
{\rm HH}^{n}({\rm N}_2(k), {\rm N}_2(k)) \cong kf_n \oplus kg_n \hspace*{5ex} (n \ge 0), 
%\left\{
%\begin{array}{ll}
%kf_n \oplus kg_n & (n=0) \\
%kf_n  & (n:\mbox{even}, n>0) \\
%kg_n & (n:\mbox{odd}),   
%\end{array}
%\right.
\]
where $f_0$ is the unit and $f_{n}=f_1^n$, $g_{n}=f_1^n g_0=g_0f_1^n$, $g_{0}^2 =0$ for $n\ge 0$. 
In particular, ${\rm HH}^{\ast}({\rm N}_2(k), {\rm N}_2(k))$ is generated by $g_0, f_1$ as a $k$-algebra.  

\begin{theorem}
Let $k$ be a field of characteristic ${\rm ch}(k) = 2$. 
For $c, c' \in k$, define an operator $\Delta_{c, c'} : {\rm HH}^{\ast}({\rm N}_2(k), {\rm N}_2(k)) \to {\rm HH}^{\ast-1}({\rm N}_2(k), {\rm N}_2(k))$ by 
\begin{eqnarray*}
\Delta_{c, c'}(f_{2n})=\Delta_{c, c'}(g_{2n})=0 & \hspace*{2ex} (n \ge 0), \\
\Delta_{c, c'}(f_{2n+1})=cf_{2n}+c'g_{2n} & \hspace*{2ex} (n \ge 0), \\
\Delta_{c, c'}(g_{2n+1})=-f_{2n}+cg_{2n} & \hspace*{2ex} (n \ge 0).  
\end{eqnarray*}
Then $\Delta_{c, c'}$ gives ${\rm HH}^{\ast}({\rm N}_2(k), {\rm N}_2(k))$ a Batalin-Vilkovisky algebra structure which induces $[\;, \;]$.  
In particular, Batalin-Vilkovisky algebra structures on ${\rm HH}^{\ast}({\rm N}_2(k), {\rm N}_2(k))$ giving 
$[\;, \;]$ are not unique.   
\end{theorem}

\proof
By direct calculation, we can verify (\ref{eq:BV}). 
\qed

\end{document}